\documentclass[3p,times]{elsarticle}

\usepackage{ecrc}
\usepackage{url}
\usepackage{todonotes}

\volume{00}

\firstpage{1}

\journalname{Journal of Computational and Applied Mathematics}

\runauth{L. Fang and L. Stals}

\jid{JCAM}

\jnltitlelogo{JCAM}

\CopyrightLine{2024}{This manuscript version is made available under the CC-BY-NC-ND 4.0 license https://creativecommons.org/licenses/by-nc-nd/4.0/}

\usepackage{amssymb}

\usepackage[figuresright]{rotating}

\usepackage{amsmath}
\usepackage{algorithm}
\usepackage{algpseudocode}
\usepackage{bm}
\usepackage{subcaption}
\usepackage{dsfont}

\begin{document}

\begin{frontmatter}

\dochead{}

\title{Data-based Adaptive Refinement of Finite Element Thin Plate Spline}

\author[1]{Lishan Fang\corref{cor1}}
\ead{fanglishan@hqu.edu.cn}
\address[1]{School of Mathematical Sciences, Huaqiao University, Xiamen, China}

\author[2]{Linda Stals}
\ead{Linda.Stals@anu.edu.au}
\address[2]{Mathematical Sciences Institute, Australian National University, ACT, Australia}

\cortext[cor1]{Corresponding author.}

\begin{abstract}
	The thin plate spline, as introduced by Duchon, interpolates a smooth surface through scattered data. It is computationally expensive when there are many data points. The finite element thin plate spline (TPSFEM) possesses similar smoothing properties and is efficient for large data sets. Its efficiency is further improved by adaptive refinement that adapts the precision of the finite element grid. Adaptive refinement processes and error indicators developed for partial differential equations may not apply to the TPSFEM as it incorporates information about the scattered data. This additional information results in features not evident in partial differential equations. An iterative adaptive refinement process and five error indicators were adapted for the TPSFEM. We give comprehensive depictions of the process in this article and evaluate the error indicators through a numerical experiment with a model problem and two bathymetric surveys in square and L-shaped domains.
\end{abstract}


\end{frontmatter}

\section{Introduction}\label{intro}

	Scattered data from real-world applications often consist of noise-perturbed data points that need to be modelled and analysed~\cite{majdiscova2017big}. The finite element thin plate spline (TPSFEM) combines the favourable properties of the thin plate spline (TPS) and finite element surface fitting, which allows it to approximate and smooth large data sets efficiently~\cite{roberts2003approximation}. The accuracy of the TPSFEM is improved by iteratively refining its finite element grid, and the efficiency is further enhanced by using adaptive refinement that only refines certain sensitive regions marked by error indicators.
 To construct the error indicators, we combine techniques from adaptive refinement of the finite element method (FEM) with techniques from large data approximation. To the best of the author's knowledge, this idea has not been previously explored.
 
	Many error indicators have been developed for the FEM to approximate partial differential equations (PDEs) and many of them also estimate errors in the solution.
 The formulation of the TPSFEM incorporates the data and consequently, these PDE-based error indicators may not be directly applicable. Alternatively, the accuracy of data approximation may be evaluated using regression metrics such as the root-mean-square error (RMSE). 
 These regression metrics can estimate the difference between the smoother and the observed data in each element of the grid. However, their performance deteriorates in the presence of noise as discussed below in Section~\ref{sec:regression}. Moreover, they do not work for regions without data points, which may be of interest in applications like surface reconstruction~\cite{carr2001reconstruction}. We have developed a technique that combines the PDE-based error indicators and regression metrics in a modified iterative adaptive refinement process.

	A numerical experiment was conducted using the peaks function from MATLAB and two bathymetric surveys, which are public data sets released from the U.S. Geological Survey\footnote{U.S. Geological Survey, https://www.usgs.gov/}, to validate their performance. The peaks function contains three local maxima and three local minima, which makes it ideal for testing adaptive refinement. The two surveys were conducted in a crater lake and a coastal region, respectively. The Crater Lake data portrays the bottom of the lake and was collected using a state-of-the-art multibeam sonar system from a survey vessel\footnote{2000 Multibeam sonar survey of Crater Lake, Oregon, https://pubs.usgs.gov/dds/dds-72/index.htm}. It aims to support both biological and geological research in the area, for example, aquatic biology and volcanic processes. The Coastal Region data contains sounding (water depth) measurements of coastal areas and may be converted to a digital terrain model to interpret the historic seafloor elevation\footnote{Historical bathymetry in the Mississippi-Alabama Coastal Region, https://coastal.er.usgs.gov/data-release/doi-P9GRUK4B//}. In contrast to previous TPSFEM studies, these two real-world data sets have more complex surfaces and boundaries. For example, the Coastal Region data's distribution fits in an L-shaped domain. See Section~\ref{sec:experiment}.

	The remainder of this article is organised as follows. A short description and formulation of the TPSFEM are provided in Section~\ref{sec:tpsfem}. The iterative adaptive refinement process of the TPSFEM is presented in Section~\ref{sec:adaptive}. The four error indicators of the TPSFEM and an alternate approach that uses regression metrics are described in Section~\ref{sec:indicator}. Results and discussions of numerical experiments are shown in Section~\ref{sec:experiment}. And lastly, a conclusion of findings is given in Section~\ref{sec:conclusion}.

\section{Finite element thin plate spline}
\label{sec:tpsfem}

	The TPS is a spline-based interpolation and smoothing technique introduced by Duchon~\cite{duchon1977splines}. Given a data set~$\left\{(\bm{x}_{i}, y_{i}):i=1,2,\ldots,n\right\}$ of size~$n$ and dimension~$d$ on a bounded domain~$\Omega$, where~$\bm{x}_{i} \in \mathbb{R}^{d}$ and~$y_{i} \in \mathbb{R}$ are~$i$-th predictor value and response value, respectively, the TPS is the function $t$ that minimises the functional
     \begin{equation} \label{eqn:tps_minimiser}
        J_{\alpha}(t) = \frac{1}{n}\sum^{n}_{i=1}\left(t(\bm{x}_{i})-y_{i}\right)^{2}+\alpha \int_{\Omega}\sum_{|\bm{v}|=2}\begin{pmatrix}2\\\bm{v}\end{pmatrix}\left(D^{\bm{v}}t(\bm{x})\right)^{2}\,d\bm{x}
    \end{equation}
    over~$H^{2}(\Omega)$, where~$\bm{v}=(v_{1},\ldots,v_{d})$, $|v|=\sum_{j=1}^{d}v_{j}$ and~$v_{j}\in\mathbb{N}$ for~$j=1,\ldots,d$. Smoothing parameter~$\alpha$ balances the goodness of fit and smoothness of~$t$.
    
    The TPS possesses many favourable properties including smoothness and rotation invariance and has been widely applied in fields including geological modelling~\cite{cowan2003practical} and rainfall interpolation~\cite{hutchinson1998interpolation}. However, it is computationally expensive and has a high memory requirement for large data sets. A number of techniques were developed to improve the efficiency of the traditional TPS. These techniques include local basis functions, some of which are discussed in Section~\ref{sec:cost}, as well as finite element-based techniques~\cite{ramsay2002spline,chen2018stochastic}.

	In this article, we focus on a mixed finite element approximation to the TPS known as the TPSFEM. The TPSFEM was presented by Roberts et al.~\cite{roberts2003approximation} to achieve smoothing properties similar to the TPS. 
 Unlike other FEM approximations that use second or higher-order finite elements, the TPSFEM is the first technique that is based on first-order elements. Specifically, it uses mixed finite elements to obtain a sparse system of equations. Further studies of the TPSFEM are done by several authors, including Stals and Roberts~\cite{stals2006smoothing} for three-dimensional data sets and Stals~\cite{stals2015efficient} for efficient solvers. A short description of the TPSFEM is given in Section~\ref{sec:discrete}. We compare the performance of the TPSFEM and four other radial basis smoothers in Section~\ref{sec:cost}.

\subsection{Discrete formulation}
\label{sec:discrete}

\begin{figure}
		\centering
		\includegraphics[width=0.5\textwidth]{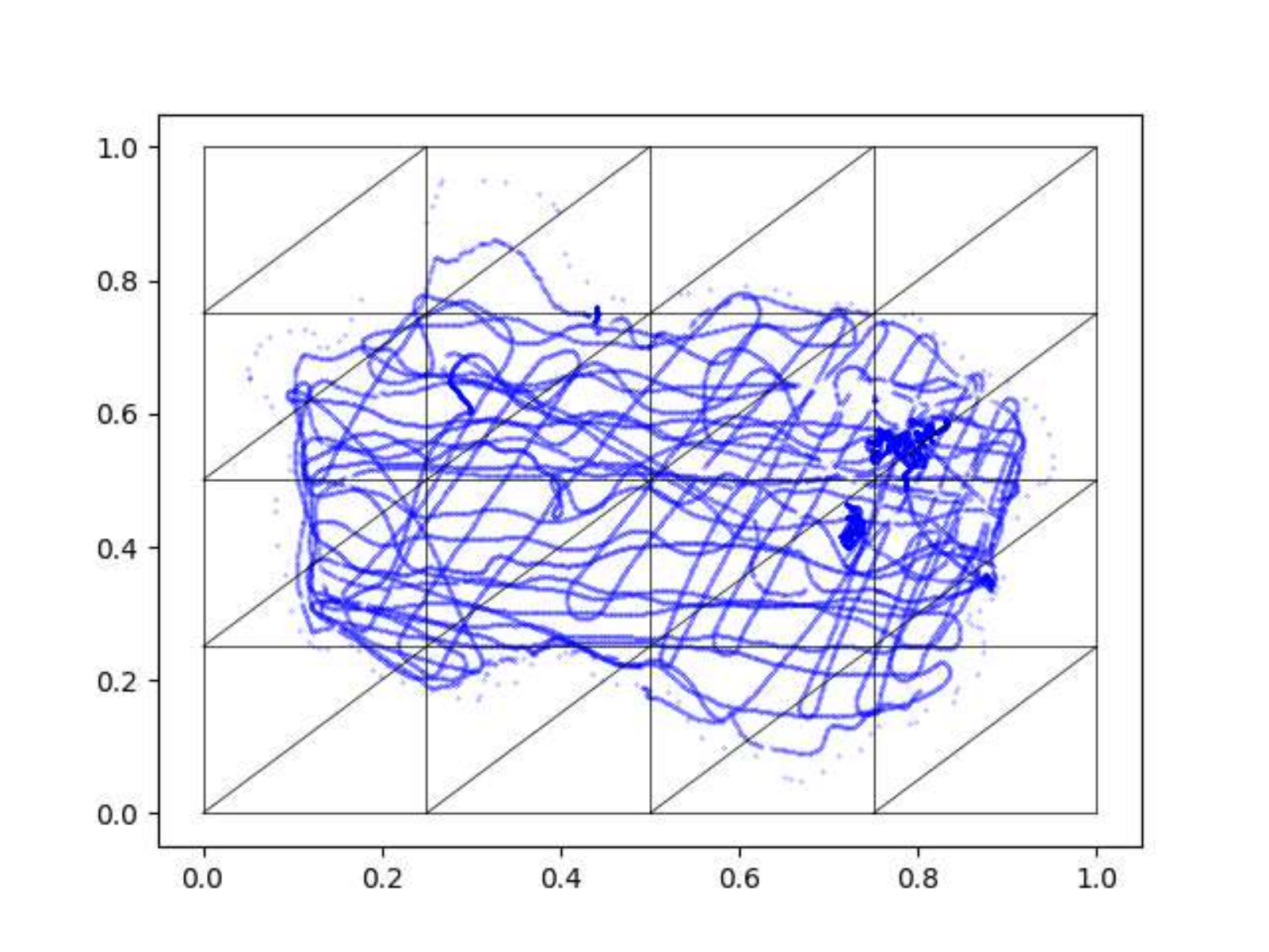}
		\caption{Data from a bathymetric survey in a triangular FEM grid. Data points are represented as blue dots.}
		\label{fig:grid_data}
	\end{figure}

	The TPSFEM is discretised over a continuous piecewise polynomial space~$\mathds{V}_{h} \in H^{1}(\Omega)$ parameterised by mesh size~$h$. A more precise definition of $h$ is given below. The TPSFEM smoothers take the form of~$s(\bm{x}) = \bm{b}(\bm{x})^{T}\bm{c}$, where~$\bm{b}=[b_{1},\ldots,b_{m}]^{T}$ are~$m$ piecewise linear basis functions, defined on a FEM grid of dimension~$d$, and~$\bm{c}$ are coefficients as described by Roberts et al.~\cite{roberts2003approximation}. An example of a two-dimensional FEM grid with a data set is shown in Figure~\ref{fig:grid_data}, which consists of 25 nodes and 32 triangular elements. The mesh size is $h = 0.25$. The basis function~$b_p(x)$ is associated with the $p$-th node in the FEM grid. The value of $b_p(\bm{x}_i)$ is zero if the data point $\bm{x}_i$ falls outside the radius of support of $b_p$. We want $s(\bm{x}_i) = \sum_{p=1}^m c_p b_p(\bm{x}_i)$ to approximate the value of the data at location $\bm{x}_i$. 

    The smoothing term of the TPS in Equation~\eqref{eqn:tps_minimiser} requires a basis function with continuous high-order derivatives. The linear basis functions used to define $s$ do not, directly, satisfy that requirement. One solution is to define $s$ in terms of quadratic polynomials, as has been done in~\cite{ramsay2002spline}. Another approach is to use a mixed finite element approach, where an auxiliary variable is introduced to represent the gradient of $s$. The advantage of the mixed finite element approach is that we can construct a formulation in terms of linear basis functions, resulting in a sparse system of equations. In their mixed finite element formulation, Roberts et al. introduced auxiliary functions~$\bm{u}$ to represent the gradients of~$s$, where
	\begin{equation*}
	\bm{u}(\bm{x})=\begin{bmatrix} u_{1}(\bm{x}) \\ \vdots \\\ u_{d}(\bm{x}) \end{bmatrix}=\begin{bmatrix} \bm{b}(\bm{x})^{T}\bm{g}_{1} \\ \vdots \\\ \bm{b}(\bm{x})^{T}\bm{g}_{d} \end{bmatrix}.
	\end{equation*}
	Function~$u_{k}(\bm{x})$ approximates the gradient of~$s$ in dimension~$k$ for~$1 \le k \le d$ and~$\bm{g}_{1},\ldots,\bm{g}_{d}$ are the corresponding coefficients. The TPSFEM $s$ and~$\bm{u}$ needs to satisfy
    \begin{equation} \label{eqn:tpsfem_relationship}
        \int_{\Omega}\nabla s(\bm{x})\cdot\nabla b_{j}(\bm{x})\,d\bm{x}= \int_{\Omega}\bm{u}(\bm{x})\cdot\nabla b_{j}(\bm{x})\,d\bm{x}
    \end{equation}
    for every basis function~$b_{j}$, which ensures~$\nabla s$ and~$\bm{u}$ are equivalent in a weak sense. In other words, $\bm{u}$ approximates the gradient of $s$ up to a constant.
    
    We find~$s$ by minimising
	\begin{align} \label{eqn:tpsfem}
		J(\bm{c},\bm{g}_{1},\ldots,\bm{g}_{d}) 
		&= \frac{1}{n}\sum^{n}_{i=1}\left(\bm{b}(\bm{x}_{i})^{T}\bm{c}-y_{i}\right)^{2}+\alpha \int_{\Omega} \sum_{k=1}^{d} \nabla \left(\bm{b}\left(\bm{x}\right)^{T}\bm{g}_{k}\right)\nabla \left(\bm{b}\left(\bm{x}\right)^{T}\bm{g}_{k}\right)\,d\bm{x} \nonumber \\
		&= \bm{c}^{T}A\bm{c}-2\bm{d}^{T}\bm{c}+\bm{y}^{T}\bm{y}/n+\alpha \sum_{k=1}^{d} \bm{g_{k}}^{T}L\bm{g_{k}},
	\end{align}
	subject to the constraint~$L \bm{c} = \sum_{k=1}^{d}G_{k}\bm{g}_{k}$, where~$\bm{y}=[y_{1},\ldots,y_{n}]^{T}$, $\bm{d} = \frac{1}{n} \sum_{i=1}^{n} \bm{b}(\bm{x}_{i})y_{i}$ and $A = \frac{1}{n} \sum_{i=1}^{n} \bm{b}(\bm{x}_{i})\bm{b}(\bm{x}_{i})^{T}$.  Matrix~$L$ is a discrete approximation to the negative Laplacian and $L_{p,q}=\int_{\Omega}\nabla\bm{b}_{p}\nabla\bm{b}_{q}\,d\bm{x}$. Matrix~$G_{k}$ is a discrete approximation to the gradient operator in dimension~$k$ for~$1 \le k \le d$ and~$(G_{k})_{pq}=\int_{\Omega}\bm{b}_{p}\partial_{k}\bm{b}_{q}\,d\bm{x}$. Note that constraint~$L \bm{c} = \sum_{k=1}^{d}G_{k}\bm{g}_{k}$ is equivalent to Equation~\eqref{eqn:tpsfem_relationship}.
 
    The matrix~$A$ and vector~$\bm{d}$ are constructed by scanning the observed data, where~$A=BB^{T}$ and~$B_{p,i}=b_{p}(\bm{x}_{i})\slash\sqrt{n}$. Intuitively, the role of matrix~$B$ is to project the data down onto the FEM grid. This is done by averaging the data that falls within the radius of support of a given basis function.  Note that the $p$-th row of $B$ will be zero if none of the data points falls within the radius of support of $b_p$. Matrices~$L$, $G_{k}$ and~$A$ are~$m \times m$ sparse matrices, where~$m$ is the number of basis functions.  This mixed finite element formulation is efficient for large data sets where $n >> m$. The smoothing parameter~$\alpha$ may be calculated automatically using techniques such as the generalised cross-validation (GCV) discussed in Section~\ref{sec:smooth}.

	By using Lagrange multipliers, Equation~\eqref{eqn:tpsfem} and the given constraint can be rewritten as a system of linear equations~\cite{stals2006smoothing}. For example,  in a two-dimensional domain with Dirichlet boundary condition, the system is 
	\begin{equation} \label{eqn:system}
		\begin{bmatrix} A & \bm{0} & \bm{0} & L \\ \bm{0} & \alpha L & \bm{0} & -G_{1}^{T} \\ \bm{0} & \bm{0} & \alpha L & -G_{2}^{T} \\ L & -G_{1} & -G_{2} & \bm{0} \end{bmatrix}\begin{bmatrix} \bm{c} \\ \bm{g}_{1} \\ \bm{g}_{2} \\ \bm{w} \end{bmatrix} = \begin{bmatrix} \bm{d} \\ \bm{0} \\ \bm{0} \\ \bm{0} \end{bmatrix}-\begin{bmatrix} \bm{h}_{1} \\ \bm{h}_{2} \\ \bm{h}_{3} \\ \bm{h}_{4} \end{bmatrix},
	\end{equation}
	where~$\bm{w}$ is a Lagrange multiplier and Dirichlet boundary information is stored in~$\bm{h}_{1}$,~$\bm{h}_{2}$,~$\bm{h}_{3}$ and~$\bm{h}_{4}$. The $\bm{h}_1, \cdots, \bm{h}_4$ vectors are not required in the case of Neumann boundary conditions. The original theory formulated by Roberts et al. used zero Neumann boundary conditions. We included both Dirichlet and Neumann boundary conditions in the numerical experiments as the Dirichlet boundary conditions allow us to set up and analyse problems with known solutions. Furthermore, the resulting system of equations with Dirichlet boundaries is better conditioned.  In the examples presented in this paper, a sparse direct solver was used to solve the above saddle point problem. However, Stals~\cite{stals2015efficient} has studied the use of a preconditioned conjugate gradient method.

	When using Dirichlet boundary conditions, the~$\bm{g}_{k}$ vectors give an approximation to the gradient of the data. When using Neumann boundary conditions, the ~$L$ and~$G_{k}$ matrices have non-trivial null spaces and consequently~$\bm{g}_{k}$ gives an approximation to the gradient of the data up to a constant. 

	Roberts et al. presented a convergence analysis for the TPSFEM with uniformly distributed data perturbed by white noise. They showed that the error convergence of the TPSFEM~$s$ depends on the smoothing parameter~$\alpha$, mesh size~$h$ and the maximum distance between data points~$d_{X}$ of its FEM grid. The mesh size~$h=\max_{\tau \in \mathcal{T}^{h}} \text{diameter}(\tau)$ of a two-dimensional triangular FEM grid is defined as the longest edge of all triangles and the distance~$d_{X}=\sup_{\bm{x}\in \Omega} \text{distance}\left(\bm{x}, \{\bm{x}_{i}\}\right)$ be the maximum distance of any point in $\Omega$ from a 
 data point. Suppose that the finite element spaces satisfy Assumptions 1 to 5 of~\cite{roberts2003approximation}, then there exist constants~$C_{1} > 1$ and~$\alpha_{0}>0$ such that for all smooth function~$f \in H^{2}(\Omega)$, which models response values~$\bm{y}$, the expected errors of the TPSFEM satisfy
	\begin{equation} \label{eqn:tpsfem_convergence_exact_1}
		E||s(\bm{y})-f||_{L^{2}(\Omega)}^{2} \preceq \left(\alpha+h^{4}+d_{X}^{4}\right)||f||_{H^{2}(\Omega)}^{2}+\frac{\sigma^{2}}{\alpha^{1\slash 2}}\frac{\left(h^{4}+d_{X}^{4}\right)}{h^{2}},
	\end{equation}
	\begin{equation} \label{eqn:tpsfem_convergence_exact_2}
		E|s(\bm{y})-f|_{H^{1}(\Omega)}^{2} \preceq \frac{1}{\alpha^{1/2}}\left(\alpha+h^{4}+d_{X}^{4}\right)||f||_{H^{2}(\Omega)}^{2}+\frac{\sigma^{2}}{\alpha}\frac{\left(h^{4}+d_{X}^{4}\right)}{h^{2}},
	\end{equation}	
	provided that~$h$ and~$\alpha$ satisfy~$h>C_{1}d_{X}$ and~$d_{X}^{4}+h^{4}<\alpha<\alpha_{0}$. The~$L^{2}(\Omega)$ norm and Sobolev semi-inner-products~$H^{1}(\Omega)$ and~$H^{2}(\Omega)$ are defined in Section 2.1 of~\cite{roberts2003approximation} and~$\sigma^{2}$ is the variance of measurement errors.  

Observe that when~$\alpha$ or~$d_{X}$ is large, the difference between the FEM approximation $s$ and the observed data $f$ may not necessarily be reduced by decreasing~$h$~\cite{fang2018error}. Since the distribution of data points is predetermined,~$d_{X}$ cannot be modified. As opposed to PDE-based error indicators that focus solely on $h$, we must estimate an optimal $\alpha$ and $h$ to minimise the expected regression errors. 

Roberts et al. assume the data points are uniformly spread in the sense that there exists a constant $C > 1$ such that $d < C\min_{i\ne j}\left|\bm{x}_i-\bm{x}_j\right|$. See Assumption 5 of ~\cite{roberts2003approximation}. We ensured this assumption held in our original study of a set of model problems used to validate the performance of the TPSFEM and our code~\cite{fang2021error}. For the real-world problems shown in this article, the assumption no longer holds. Nevertheless, the behaviour of the numerical results is similar to what was observed previously with our model problems.

\subsection{Comparison with radial basis functions}
\label{sec:cost}

	To validate the TPSFEM's competitiveness, we compared the performance of the TPSFEM to the TPS~\cite{wahba1990spline} and several radial basis functions (RBFs) with compact support (CSRBF) from Wendland~\cite{wendland1995piecewise} and Buhmann~\cite{buhmann1998radial} using the Crater Lake survey~\cite{bathymetric2000lake} described in Section~\ref{sec:experiment}. 

	To find the control points for the RBFs, we place a uniform rectangular grid of mesh size $\bar h$ over the domain and find those data points closest to the nodes of the grid. The grid is similar to what is shown in Figure~\ref{fig:grid_data}, except rectangles are used instead of triangles. Regions of the domain where the distance between the grid nodes and the nearest data point is greater than $\bar h/3$ are ignored. 
    With this approach, it is possible to bound the maximum distance between neighbouring control points. Given the error estimate depends on the maximum distance between neighbouring points, we can construct a corresponding bound on the error estimate of the RBF that will decrease as the number of rectangles, and in turn the number of control points, is increased. Furthermore, the minimum distance between control points is bounded by $\bar h/3$. Let $\bar n$ be the number of control points extracted from the Crate Lake data set.

	Consider a control point $\bar{\bm{x}}_i$ where $1 \le i \le \bar{n}$. The RBF kernel defined at point $\bar{\bm{x}}_i$ is~$\Psi_i(\bm{x}) = \Phi_i(r)$ where~$r=||\bar{\bm{x}}_i-\bm{x}||_{2}$ is the Euclidean distance between $\bar{\bm{x}}_i$ and $\bm{x}$. Examples of $\Phi_i$ used in the paper are listed in Table~\ref{fig:comp_rbf}. The function~$(1-r)_{+}$ is a truncated power function and~$(1-r)_{+}=0$ if~$(1-r)<0$. Thus, the Wendland RBFs have compact support. Buhmann's RBFs also have compact support as the kernel is defined by~$\Phi_i(r) =0$ if $r > 1$. The three CSRBFs are scaled using a radius of support~$\rho$ and the kernel values are calculated using $\Phi_i(r/\rho)$. 
    A fixed number of data points fall within the circle with centre $\bar{\bm{x}}_i$ and radius $\rho$, see~\cite{deparis2014rescaled}. In this report, the radius~$\rho$ is calculated so that the local support of each kernel contains about 100 or 200 sampled data points. For instance, the radius of support that covers 200 sampled points when the number of sample points is 946, 3668 and 14477 are 0.16, 0.13 and 0.066, respectively. 

	\begin{table}
		\centering
		\caption{Comparison of computational costs}
		\label{fig:comp_rbf}
		\begin{tabular}{lllllllll}
			\hline\noalign{\smallskip}
			Technique & Kernel $\Phi_i$ & \# basis & Radius $\rho$ & \# nonzero & Ratio & Time & RMSE \\
			\noalign{\smallskip}\hline\noalign{\smallskip}
			TPSFEM & & 4,225 & & 181,352 & 0.072\% & 0.09 & 9.29 \\
			& & 16,641 & & 608,642 & 0.015\% & 0.43 & 5.50 \\
            & & 66,049 & & 2,454,404 & 0.004\% & 2.12 & 2.34 \\
			\hline\noalign{\smallskip}
			TPS & $r^{2}\log(r)$ & 946 & & 894,916 & 100\% & 0.22 & 9.92 \\
			& & 3,668 & & 13,454,224 & 100\% & 1.67 & 5.32 \\
			& & 14,477 & & 209,583,529 & 100\% & 18.10 & 2.49 \\
			\hline\noalign{\smallskip}
			Wendland & $(1-r)^{2}_{+}$ & 946 & 0.13 & 106,458 & 11.90\% & 0.14 & 11.10 \\
            & & 946 & 0.16 & 198,312 & 22.16\% & 0.18 & 10.84 \\
            & & 3,668 & 0.093 & 455,572 & 3.39\% & 0.29 & 5.80 \\
			& & 3,668 & 0.13 & 901,208 & 6.70\% & 0.40 & 5.64 \\
			& & 14,477 & 0.047 & 1,985,133 & 0.95\% & 1.80 & 2.96 \\
            & & 14,477 & 0.066 & 3,921,795 & 1.87\% & 2.22 & 2.75 \\
			\hline\noalign{\smallskip}
			Wendland & $(1-r)^{4}_{+}(4r+1)$ & 946 & 0.13 & 106,458 & 11.90\% & 0.15 & 10.03 \\
            & & 946 & 0.16 & 198,308 & 22.16\% & 0.21 & 9.99 \\
			& & 3,668 & 0.093 & 455,572 & 3.39\% & 0.30 & 5.37 \\
            & & 3,668 & 0.13 & 901,208 & 6.70\% & 0.42 & 5.37 \\
            & & 14,477 & 0.047 & 1,985,133 & 0.95\% & 1.80 & 2.56 \\
            & & 14,477 & 0.066 & 3,921,795 & 1.87\% & 2.20 & 2.51 \\
			\hline\noalign{\smallskip}
			Buhmann & $1/3+r^{2}-4r^{3}/3$ & 946 & 0.13 & 106,458 & 11.90\% & 0.15 & 10.09 \\
            & $+2r^{2}\log(r)$ & 946 & 0.16 & 198,308 & 22.16\% & 0.17 & 9.98 \\
            & & 3,668 & 0.093 & 455,572 & 3.39\% & 0.35 & 5.42 \\
			& & 3,668 & 0.13 & 901,208 & 6.70\% & 0.46 & 5.33 \\
			& & 14,477 & 0.047 & 1,985,133 & 0.95\% & 1.94 & 2.62 \\
            & & 14,477 & 0.066 & 3,921,795 & 1.87\% & 2.79 & 2.52 \\
			\hline\noalign{\smallskip}
		\end{tabular}
	\end{table}
 
	The resulting RBF approximation is $s_{rbf}(\bm{x}) = \sum_{i=1}^{\bar{n}}w_i \Psi_i(\bm{x})$ for coefficients $w_i \in \mathbb{R}$. The TPSFEM is constructed using Dirichlet boundaries with uniform grids. The Dirichlet boundary conditions~$\bm{h}_i$ are defined using constants described in Section~\ref{sec:experiment}.
 
    The smoothing parameter~$\alpha$ is calculated using the GCV as discussed in Section~\ref{sec:smooth}.

	We compared the accuracy and efficiency of the different kernels using the five attributes listed in Table~\ref{fig:comp_rbf}. The number of basis is $\bar n$ for the RBFs and $m$ for the TPSFEM. The accuracy of the approximation associated with each kernel is measured using the~$\text{RMSE} = \big(\frac{1}{n}\sum_{i=1}^{n}\big(\hat{y}_{i}-y_{i}\big)^{2}\big)^{\frac{1}{2}}$, where~$\hat{y}_{i} = s_{rbf}(\bm{x}_i)$ for the RBF kernels and~$\hat{y}_i = s(\bm{x}_i)$ for the TPSFEM. We also evaluated their efficiency using the number and ratio of nonzero entries in the system of equations and the time to solve the system. The number and ratio of nonzero entries is a measure of the sparsity of the system, which may be exploited by sparse solvers~\cite{beatson1997fast}. Since sparse solvers may be affected by other factors like the structure of the systems, we also included the time for solving the system of equations, measured in seconds. 
 
    The TPSFEM program\footnote{TPSFEM program, https://bitbucket.org/fanglishan/tpsfem-program/} was implemented using Python 3.8, and we used a direct solver for the TPS and a sparse solver (PyPardiso\footnote{PyPardiso,https://pypi.org/project/pypardiso/}) for both the TPSFEM and CSRBFs. 

    As shown in Table~\ref{fig:comp_rbf}, the CSRBFs achieved slightly higher RMSE compared to the TPS but required markedly less time to solve the system of equations. The matrix associated with the TPS is a full matrix and therefore computationally expensive. The CSRBFs using a radius that covers 200 points have approximately twice the number of nonzero entries compared to the ones using a radius that covers 100 points. Thus their systems took more time to solve, but they achieved relatively lower RMSE. The time required to solve the TPSFEM was bounded below and above by the time required to solve the CSRBFs, depending on the radius $\rho$. However, irrespective of $\rho$, the TPSFEM achieved lower RMSE.
 
	We are aware the performance of the RBFs may be improved by adopting additional modifications such as a localised radius of support or more efficient parallel solvers~\cite{deparis2014rescaled}. Nevertheless, we believe this example demonstrates the TPSFEM's competitiveness. As discussed above, the system associated with the TPSFEM is built using low-order basis functions with local support, which is beneficial when smoothing large data sets.

\subsection{B-splines, an alternative spline function}\label{B-spline}
Although we have focused on thin plate splines, other spline functions, such as the tensor B-splines 
are available. Tensor B-splines do not require a solution of a global matrix, unlike the TPS method. Instead, they solve a local dense matrix whose size depends on the order of the polynomial spline. The solution of the local matrix problems is efficient for low-order polynomials. The cost of the dense matrix solve is listed as one of the disadvantages of the TPS  method compared to B-splines. Our TPSFEM method is specifically designed to reduce that cost.

 Tensor B-splines are defined on a rectangular domain. The splines can be extended to more general domains by using harmonic transformations \cite{Xu_2014, 10.1145/1364901.1364938, VUILLOD2024116913}. However, the FEM grid used in the TPSFEM implementation is more flexible and only makes minor assumptions on the structure of the domain boundary. 
 
Adaptive refinement techniques for tensor B-splines have also been developed to adjust the shape of the grid according to the data distribution \cite{BRACCO2018239, SKYTT2022102894}. However, Skytt et. al. \cite{SKYTT2022102894} observe that using statistical metrics as an error indicator is unreliable when working with noisy data. As discussed in Section \ref{sec:regression}, we observed similar behaviour and therefore concentrated on error indicators that work directly on Equation (1).

An advantage of the TPS splines is that they are optimal in the appropriate Sobolev space. Furthermore, well-established techniques exist for choosing the smoothing parameter $\alpha$. Observe that only one smoothing parameter is required, which is advantageous compared to the tensor splines that require as many smoothing parameters as the number of predictors \cite{KALOGRIDIS2023}. 

\section{Adaptive refinement}
\label{sec:adaptive}

	Adaptive refinement is a process that iteratively refines elements in sensitive regions, such as peaks, to improve the accuracy of a solution. Consequently, compared to uniform grids, adaptive grids require fewer nodes to reach a satisfactory solution. An error indicator is critical to adaptive refinement as it identifies these sensitive regions~\cite{gratsch2005posteriori}.

	In this article, we examine adaptive refinement and error indicators over two-dimensional grids with triangular elements. The triangles are refined using the newest node bisection method. Each triangle will have one, and only one, vertex labelled as the {\it newest node}. A triangle is divided into two children triangles by adding an edge from the newest node to the middle point of the opposite edge~\cite{zienkiewicz2005finite}. That midpoint becomes the newest node for the next refinement iteration. This helps to prevent long thin triangles and limits the growth of the interpolation error~\cite{babuka1976angle}.

Our data structure is designed to be flexible enough to store one-, two- or three-dimensional grids. It does not explicitly store the triangles. The triangles are constructed when needed using the geometric and topological information contained in the structure. Within such a setting, it is easier to work with the edges (topological information) of the triangles rather than the nodes. Given a triangle, the edge opposite the newest nodes is labelled a base-edge. There is a one-to-one correspondence between the base-edges and the newest nodes.

 
 Any two triangles that share a base-edge are called a triangle pair. During uniform refinement, all the triangles are bisected along the base-edge. The remaining old edges in each triangle become the new set of base-edges. This controls the order of refinement to, again, avoid long thin triangles. During adaptive refinement, only some of the triangle pairs are bisected. This may result in edges that are shared by coarse (old) and fine (new) triangles. Such edges are labelled as {\it interface base-edges}. Triangles can only be bisected along base-edges. If an interface base-edge is to be bisected, the coarse triangle sharing that edge is bisected first. That interface base-edge is then marked as a base-edge, and bisection can proceed for the triangle pair. This avoids hanging nodes. It is possible the neighbouring coarse triangle mentioned previously also contains an interface base-edge. A recursive algorithm is used to move down the neighbouring coarse triangles until a base-edge is found. Mitchell~\cite{mitchell1989comparison} and Stals~\cite{stals1995parallel} showed this recursive algorithm will terminate. These error indicators and the newest node bisection may also be extended to other dimensions~\cite{gratsch2005posteriori}. We show more details of our iterative adaptive refinement process in Section~\ref{sec:iterative} and how the smoothing parameter~$\alpha$ is updated iteratively in this process in Section~\ref{sec:smooth}.

\subsection{Iterative adaptive refinement process}
\label{sec:iterative}

	The PDE-based iterative adaptive refinement process starts with a coarse initial grid and refines it iteratively using an error indicator until the solution reaches a given error tolerance.
    In each iteration, elements with an error indicator value higher than a given error tolerance are marked and refined. A new solution is built using the refined grid and the error indicator values will be compared to the error tolerance to determine whether the iterative process should continue or terminate.

	Error indicators developed for PDE formulations are generally based on the discretisation error of the FEM approximation. Consequently, their calculations are purely based on~$h$. However, in our TPSFEM formulation, the approximation error is based on several factors, including~$\alpha$,~$h$ and~$d_{X}$, as mentioned in Section~\ref{sec:discrete}. Therefore, standard PDE-based approaches may not work, especially in the case with large~$\alpha$ or~$d_{X}$ as demonstrated in Section 4 of~\cite{fang2018error}.

	Some PDE-based error indicators measure the rate of change and a smooth approximation will give small error indicator values, thus terminating the iterative process. However, this approach in TPSFEM may lead to premature termination. A smooth function does not necessarily fit the data so the RMSE may still be high~\cite{fang2021error}. Therefore, the standard procedure described in Mitchell~\cite{mitchell1989comparison} does not work for the TPSFEM. A modified iterative process for the TPSFEM is shown in Algorithm~\ref{alg:tpsfem_adaptive} and it contains an additional inner loop from steps 9 to 14. We iteratively mark and refine a small subset of elements with relatively high error indicator values in each iteration of the inner loop. This inner loop is guaranteed to add new nodes, and it terminates as soon as the FEM grid contains a preset number of nodes. In our experiments, we terminated the inner loop when the number of nodes has been doubled. The reason is that we intend to explore the use of the multigrid method as an iterative solver in future studies. However, this is a somewhat arbitrary choice and specific problems may benefit from a different choice. This inner loop avoids the limitations of a pre-set error tolerance for large~$\alpha$ and~$d_{X}$. 

	Referring to Algorithm~\ref{alg:tpsfem_adaptive}, when elements are refined, new nodes added in triangular grid~$\mathcal{T}_{k,j+1}$ take average values of their neighbouring nodes. We do not solve Equation~\eqref{eqn:system} until the inner loop terminates. We use~$s_{k}$ and~$s_{k,j+1}$ to represent the smoother~$s$ defined on grid~$\mathcal{T}_{k}$ and~$s$ that is defined on grid~$\mathcal{T}_{k,j+1}$. Since new nodes take average values of neighbouring nodes, $s_{k,j+1}$ is only an approximation. We work with the approximation because solving Equation~\eqref{eqn:system} each time a new node is added would be exorbitantly expensive. However, only using an approximation may affect the performance of the error indicators. This is a trade-off between cost and accuracy. After this procedure terminates, the full system is solved to improve the accuracy of~$s$. This allows us to estimate the RMSE of updated~$s$ stored in~$s_{k+1}$ for the stopping criteria. 

	\begin{algorithm}
		\caption{Solve TPSFEM with adaptive refinement}
		\label{alg:tpsfem_adaptive}
		\begin{algorithmic}[1]
            \Statex \textbf{Input:} Data and initial FEM grid~$\mathcal{T}_{0}$
            \Statex \textbf{Output:} TPSFEM smoother $s$
		\State $k \gets 0$
		\State calculate optimal~$\alpha_{k}$ on~$\mathcal{T}_{k}$
		\State solve Equation~\eqref{eqn:system} to build~$s_{k}$ on~$\mathcal{T}_{k}$
		\State evaluate RMSE of~$s_{k}$
		\While{RMSE is higher than error tolerance}
			\State calculate error indicators for edges in~$\mathcal{T}_{k}$
			\State $j \gets 0$
			\State $\mathcal{T}_{k,j} \gets \mathcal{T}_{k}$
			\While{size of~$\mathcal{T}_{k,j}$ is smaller than twice of size of~$\mathcal{T}_{k}$}
				\State mark a subset of elements with high error indicator values
				\State refine marked elements to construct~$\mathcal{T}_{k,j+1}$ and obtain~$s_{k,j+1}$ by averaging
				\State calculate error indicators for new edges in~$\mathcal{T}_{k,j+1}$
				\State $j \gets j+1$
			\EndWhile
			\State $\mathcal{T}_{k+1} \gets \mathcal{T}_{k,j}$
			\State calculate optimal~$\alpha_{k}$ on refined grid~$\mathcal{T}_{k+1}$
			\State solve Equation~\eqref{eqn:system} to build~$s_{k+1}$ on~$\mathcal{T}_{k+1}$
			\State evaluate RMSE of~$s_{k+1}$
			\If{RMSE is reduced less than 10\% for two iterations}
			\State break
			\EndIf
			\State $k \gets k+1$
		\EndWhile
        \State $s \gets s_{k}$
		\end{algorithmic}
	\end{algorithm}

	Traditional stopping criteria of the iterative process are based on fixed error tolerance for errors approximated by error indicators~\cite{gratsch2005posteriori}. In contrast, our new stopping criteria are a combination of fixed error tolerance and rates of change of the RMSE~\cite{fang2021adaptive}. The iterative process stops when the RMSE of~$s$ is lower than the error tolerance. Note the RMSE may not converge to zero if the data is not smooth enough or perturbed by noise. 
 We terminate the iterative process when the addition of new nodes ceases to significantly reduce the RMSE. Specifically, we terminate the iterative process when the RMSE is reduced by less than 10\% by one iteration of refinement for two consecutive iterations. This criterion performed well in previous experiments given various model problems, data sizes and noise levels shown in Chapter 8 of~\cite{fang2021error}.

\subsection{Smoothing parameter}
\label{sec:smooth}

	Generalised cross-validation (GCV) is a popular statistical technique for model validation and variable selection. Golub~\cite{golub1979generalized} developed the GCV method to estimate the ridge parameter for use in ridge regression. The GCV must calculate the trace of a matrix inverse to estimate~$\alpha$, which is expensive for medium to large-scale problems. Hutchinson~\cite{hutchinson1990stochastic} showed the trace term can be approximated using an unbiased stochastic estimator. This approach was used with the TPSFEM in ~\cite{roberts2003approximation, stals2006smoothing, stals2015efficient}.

	The error convergence rates of~$s$ depend on the smoothing parameter~$\alpha$ as discussed in Section~\ref{sec:discrete}. The value of an optimal~$\alpha$ is affected by a range of factors, including the size, noise level and smoothness of the observed data; and~$h$. Stals~\cite{stals2015efficient} investigated a wide range of~$\alpha$ values from~$10^{-10}$ to~$1$ to test the robustness of solvers. In practice, the optimal~$\alpha$ of the TPSFEM often ranges from $10^{-10}$ to~$10^{-4}$.

	\begin{figure}
		\centering
		\includegraphics[width=0.7\textwidth]{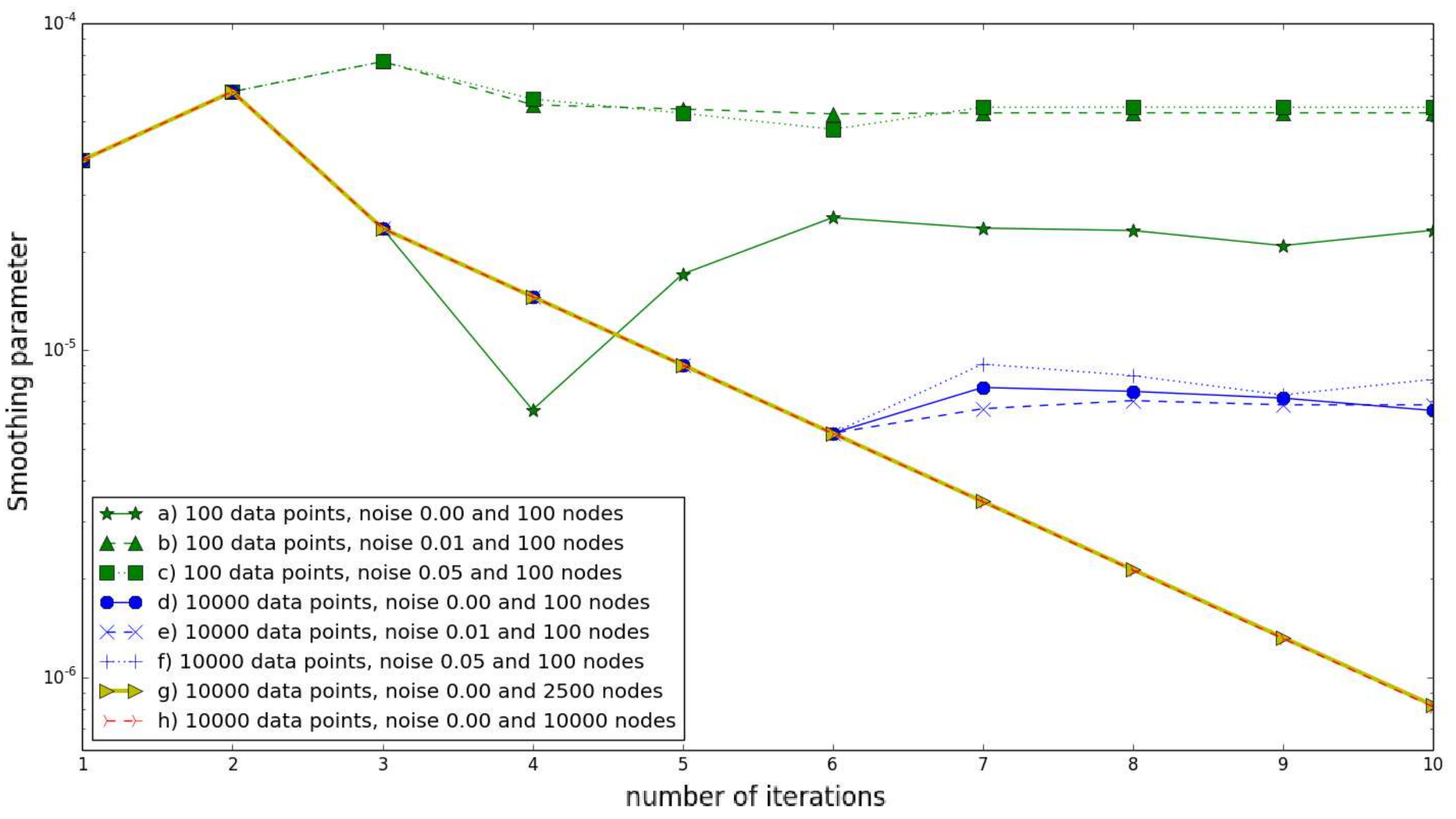}
		\caption{$\alpha$ values at each iteration of bounded minimisation with various numbers of data, nodes and noise levels.}
		\label{fig:alpha_changes}
	\end{figure}

	The iterative process in Algorithm~\ref{alg:tpsfem_adaptive} calculates a new~$\alpha$ in each iteration. Note that the value of~$\alpha$ found in iteration~$k$ is labelled as~$\alpha_{k}$. We speed up the computation of~$\alpha_{k}$ in iteration~$k$ using~$\alpha_{k-1}$ from the previous iteration~\cite{fang2020adaptive}. Fang~\cite{fang2021error} argued the optimal~$\alpha$ generally decreases as~$h$ decreases with each refinement iteration until stabilising inside a fixed range. An example from Fang is shown in Figure~\ref{fig:alpha_changes}, where~$\alpha$ values were obtained for different~$h$ using the Brent method~\cite{brent2013algorithms}. Therefore,~$\alpha_{k}$ is calculated as~$\alpha_{k} = \arg\min\{V(\alpha_{k-1}),V(r_{1}\alpha_{k-1}),V(r_{2}\alpha_{k-1})\}$ where~$0<r_{1},\,r_{2}<1$. Function~$V(\alpha)$ is the GCV function evaluated at $\alpha$ using Hutchinson's method to approximate the trace. In our experiments, we chose~$r_{1}=0.1$ and $r_{2}=0.3$, which are consistent with the convergence of~$\alpha$. The~$\alpha_{0}$ value for the initial grid in Algorithm~\ref{alg:tpsfem_adaptive} is still computed using bounded minimisation within interval~$[10^{-10},~10^{-4}]$. This approach is simple but robust, especially given that small changes in~$\alpha$ do not have a marked impact on~$s$. Fang showed this approach produces results similar to Brent's method which required ten evaluations of~$V$. The experiments were carried out on both uniform and adaptive refinement for various data sets and noise levels.

\section{Error indicators}
\label{sec:indicator}

Given our focus is on data regression, we first explore the use of regression metrics as an error indicator in Section~\ref{sec:regression}. As an alternative approach, four PDE-based error indicators were adapted for the TPSFEM. Namely, the auxiliary problem error indicator, residual-based error indicator, recovery-based error indicator and norm-based error indicator. More details are given in Chapter 5 of~\cite{fang2021error}. Two of the error indicators use global information like the grid and smoother~$s$, and the other two also use local data points directly.  A brief description of the four PDE-based error indicators is given in Sections~\ref{sec:auxiliary} to~\ref{sec:norm}.

The error indicator values are assigned to the base and interface base edges introduced in Section \ref{sec:adaptive}. In the approaches described below, the calculations are performed on all triangles sharing the specified edge. In the interior of the domain this will be a pair of triangles, along the boundary it may be a single triangle. 

\subsection{Regression metric as an error indicator}
\label{sec:regression}

	Since the observed data is available for the TPSFEM, we explored several approaches that use data directly to indicate errors, including regression metrics~\cite{shcherbakov2013survey}. While many regression metrics like the RMSE were developed to assess surface fitting techniques, they have not been used for adaptive refinement. 
 
 The RMSE is measured locally by computing the difference between smoother~$s$ and the data points lying inside a triangle. Specifically, consider an edge~$e$. Let $\tau_e$ be the set of triangles that have $e$ as an edge. Furthermore, let $\{(\hat{\bm{x}}_{i}, \hat y_{i}):i=1,2,\ldots,n_{e}\}$ be the set of data points lying inside $\tau_e$. If a data point sits on a vertex or an edge, it is only included in the calculation for one triangle. Thus, it is not counted multiple times. The error indicator value~$\eta_{e}$ is calculated as
    \begin{equation*}
        \boxed{\eta_{e} = \left(\frac{1}{n_{e}}\sum_{i=1}^{n_{e}}\left(s\left(\hat{\bm{x}}_{i}\right)-\hat{y}_{i}\right)^{2}\right)^{\frac{1}{2}}.}
    \end{equation*}
    
	\begin{figure}
		\centering
		\begin{subfigure}[b]{0.38\textwidth}
		\centering
		\includegraphics[width=\textwidth]{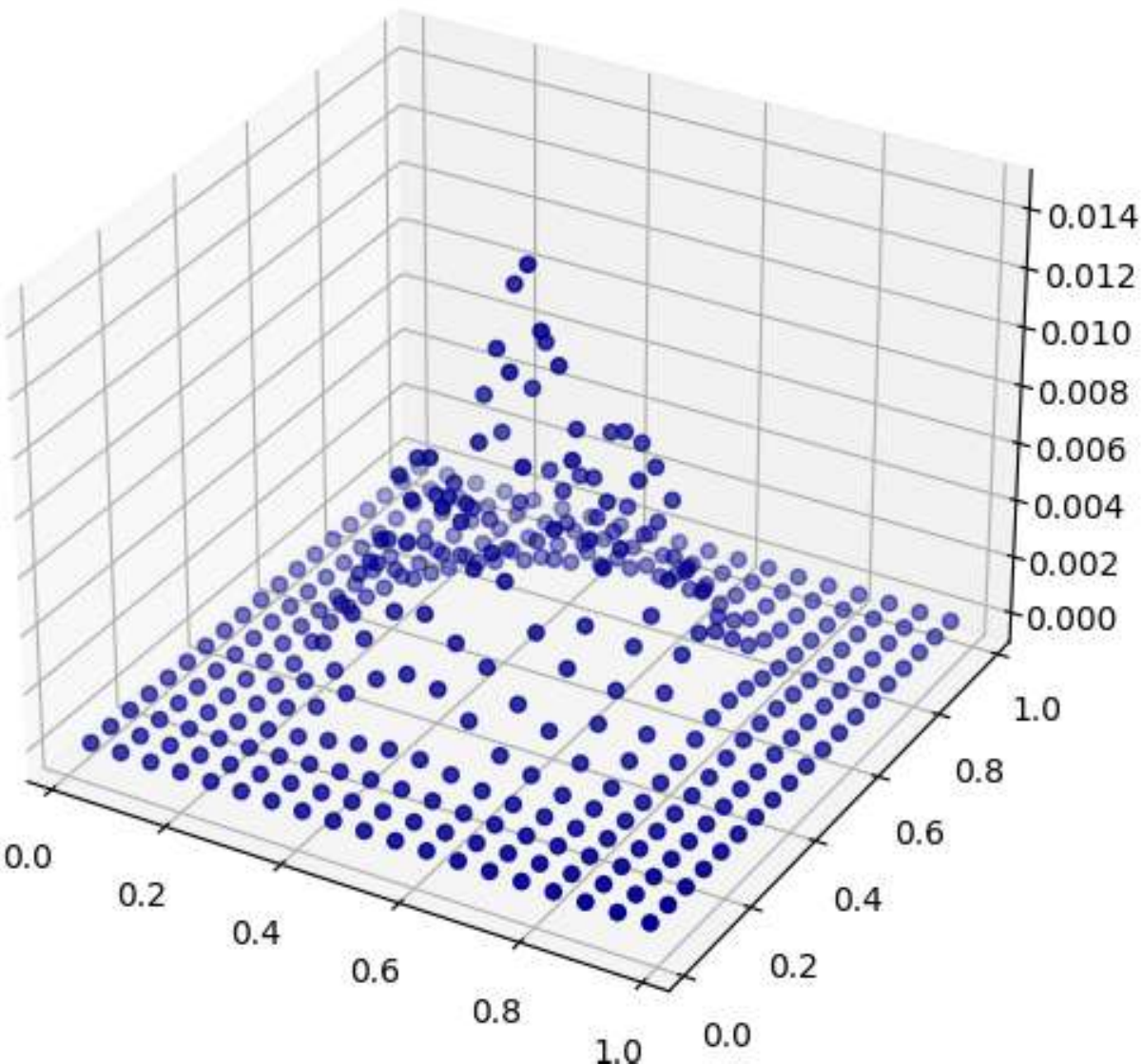}
		\caption{}
		\end{subfigure}
		\hspace{0.5cm}
		\begin{subfigure}[b]{0.38\textwidth}
		\centering
		\includegraphics[width=\textwidth]{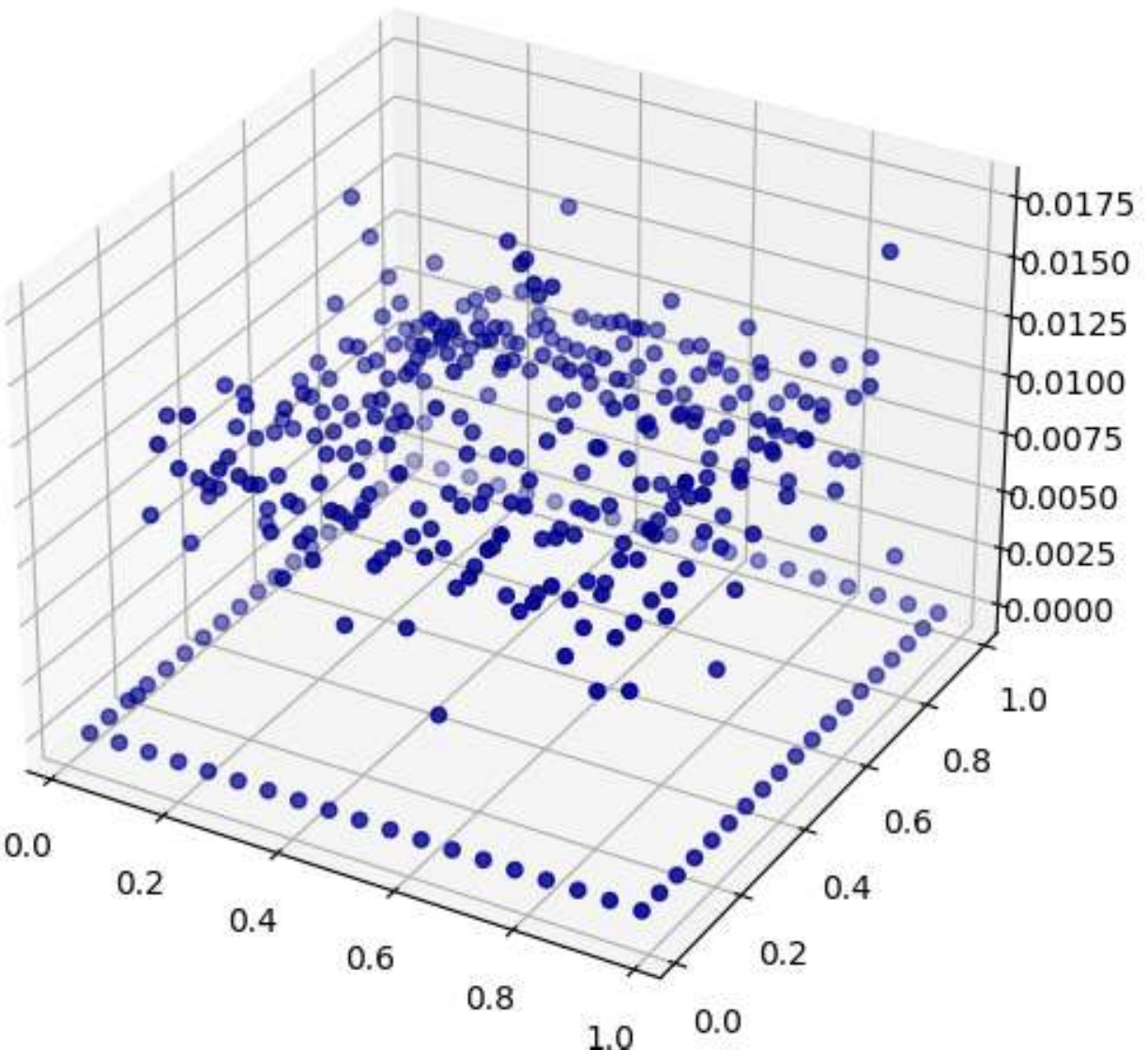}
		\caption{}
		\end{subfigure}
		\caption{RMSE of~$s$ against data sets modelled by~$y=e^{-50(x_{1}-0.5)^{2}}e^{-50(x_{2}-0.5)^{2}}$ (a) without noise; and (b) with Gaussian noise. The points represent RMSE in each triangle pair.}
		\label{fig:regression}
	\end{figure}

	\begin{figure}
		\centering
		\begin{subfigure}[b]{0.38\textwidth}
		\centering
		\includegraphics[width=\textwidth]{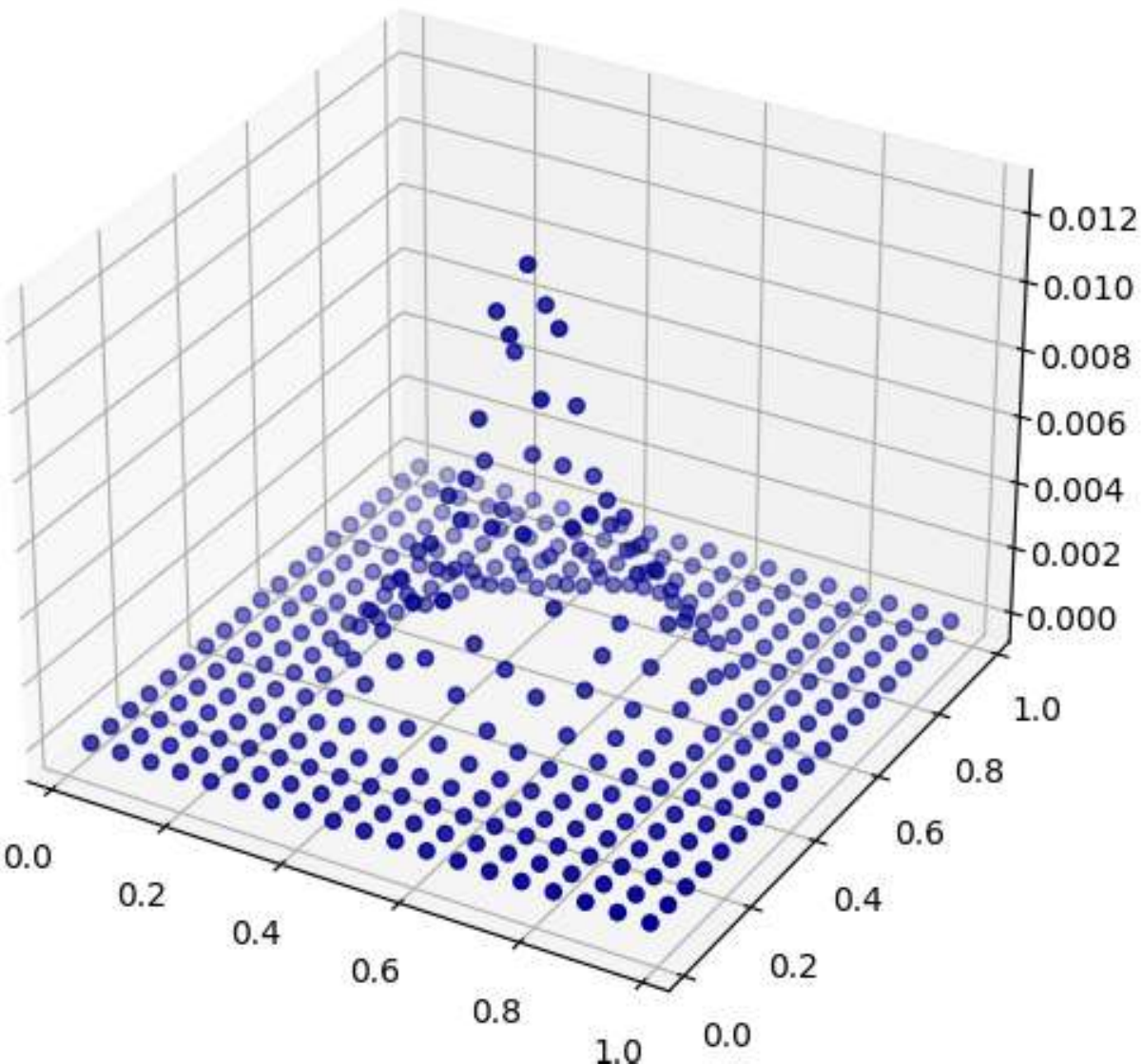}
		\caption{}
		\end{subfigure}
		\hspace{0.5cm}
		\begin{subfigure}[b]{0.38\textwidth}
		\centering
		\includegraphics[width=\textwidth]{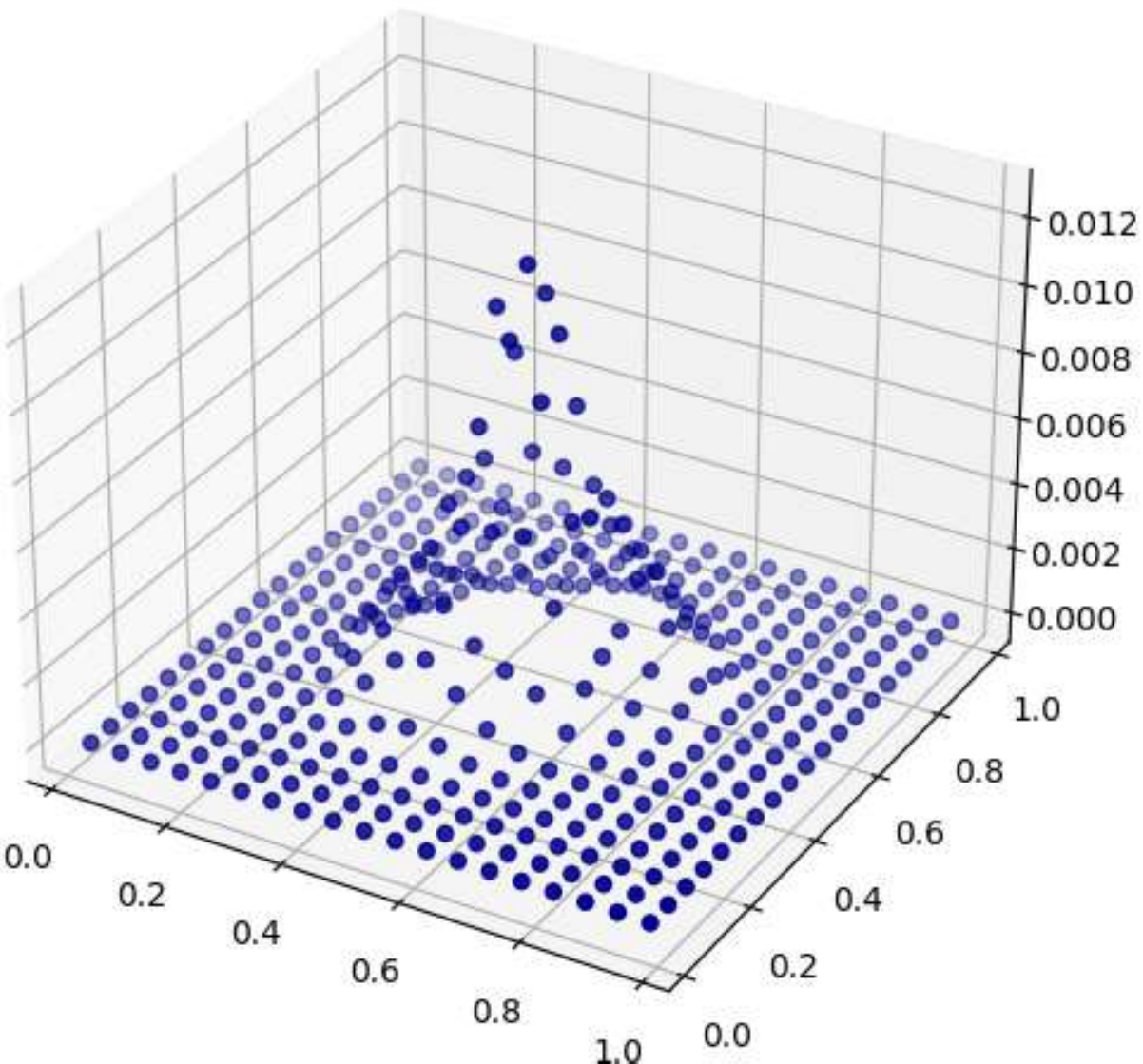}
		\caption{}
		\end{subfigure}
		\caption{Recovery-based error indicator values for the data set modelled by~$y=e^{-50(x_{1}-0.5)^{2}}e^{-50(x_{2}-0.5)^{2}}$ (a) without noise; and (b) with Gaussian noise. The points represent error indicator values in each triangle pair.}
		\label{fig:indicator}
	\end{figure}
 
We tested the RMSE as an error indicator using a two-dimensional uniformly distributed data set of size~$1 \times 10^4$ in~$[0,1]^{2}$ square domain using a uniform FEM grid with 400 nodes. The tests were carried out without noise and in the presence of Gaussian noise with mean zero and standard deviation 0.01. The range of~$\hat y$ values is 0 to 1. The resulting error indicator values are shown in Figure~\ref{fig:regression}. When the data is not perturbed by noise, the RMSE is high at the centre as shown in Figure~\ref{fig:regression}(a). The model problem changes more rapidly in this region and indicated high RMSE suggests more refinement is required at the centre. However, the RMSE no longer provides useful information when the data is perturbed by noise as shown in Figure~\ref{fig:regression}(b), and is thus not an effective error indicator in that scenario. Furthermore, the RMSE does not work if there are elements without data points, which has consequences for applications like surface reconstructions~\cite{carr2001reconstruction}. Several other regression metrics were tested, and they performed poorly in the presence of noise. Moreover, they are sensitive to changes in data distribution patterns and may lead to over-refinement. 
 
For comparison, we reran the experiment using the PDE-based error indicator defined in Section~\ref{sec:recovery}. The results are provided in Figure~\ref{fig:indicator}. The recovery-based error indicator accesses the data indirectly and is less sensitive to noise. Its error indicator values in Figures~\ref{fig:indicator}(a) and~\ref{fig:indicator}(b) are almost identical. We tested the regression as an error indicator and provided more statistics and adaptive grids in Section~\ref{sec:result_model}.

\subsection{Auxiliary problem error indicator}
\label{sec:auxiliary}

The auxiliary problem error indicator solves an auxiliary problem on smaller domains to obtain locally more accurate approximations~\cite{ainsworth2011posteriori,babuvvska1978error}. The error of the solution is approximated locally by the difference between the global approximation~$s$ and a local approximation with improved accuracy. Similar ideas have also been applied in surface fitting techniques to assess approximation quality using a cluster of neighbouring data points~\cite{zhang2017adaptive}. 



\begin{figure}
		\centering
		\includegraphics[width=0.4\textwidth]{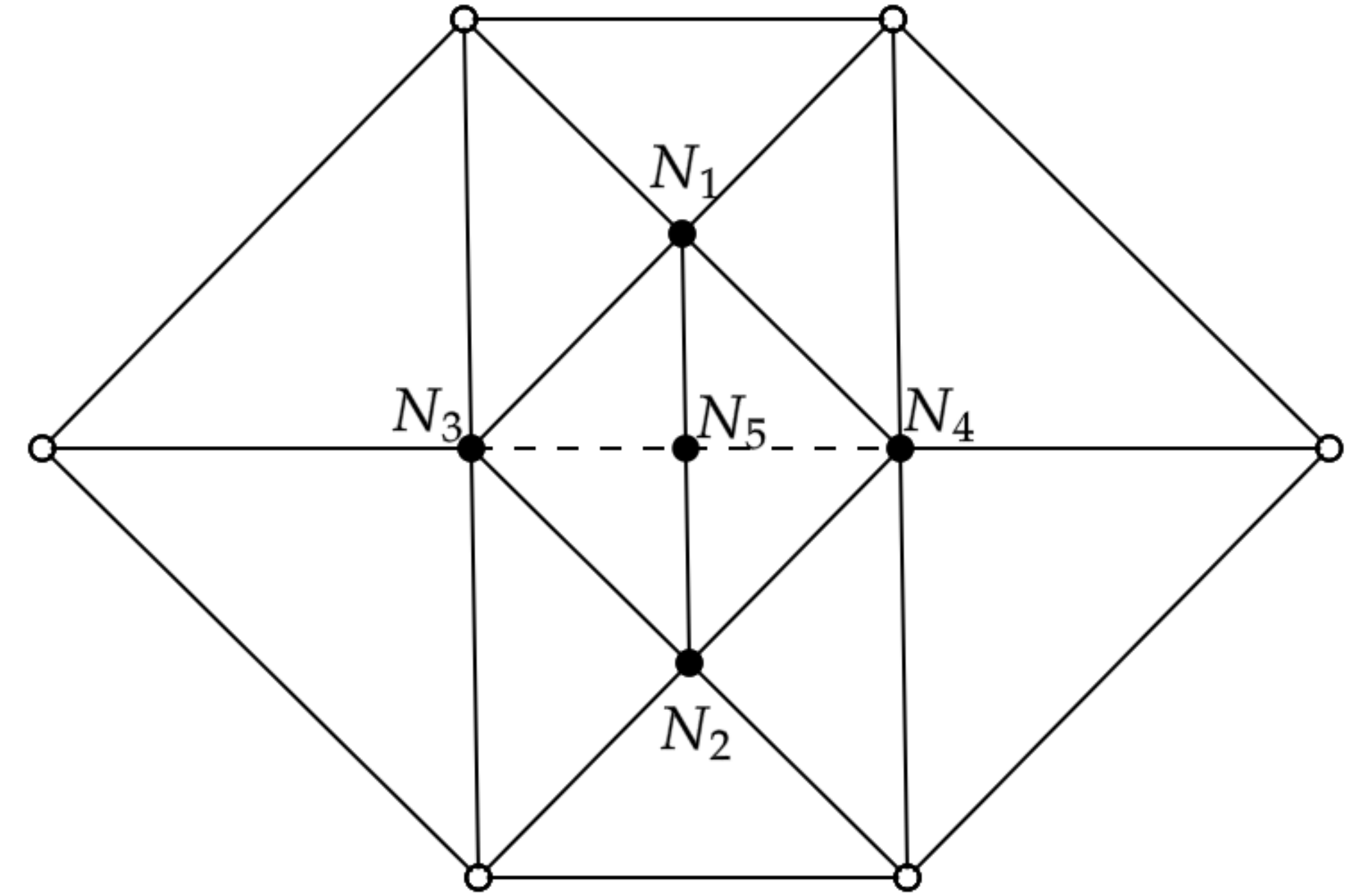}
		\caption{An example local domain. Interior nodes and boundary nodes are represented as filled circles and open circles, respectively. Accuracy is improved by refining edge~$N_{1}\text{-}N_{2}$ to introduce a new vertex~$N_{5}$ and new edges represented as dashed lines.}
		\label{fig:local_grid}
	\end{figure}

As in Section~\ref{sec:regression}, let $e$ be a base or interface base edge, and $\tau_e$ be the set of triangles sharing that edge. We then define a local grid $\mathcal{T}_e$ to be all the triangles that share a node with the triangles in $\tau_e$. For example, referring to Figure~\ref{fig:local_grid}, $e$ is the edge between nodes $N_1$ and $N_2$, $\tau_e$ contains the triangles $\{N_1, N_2, N_3\}$ and $\{N_1, N_2, N_4\}$, and the local grid is all the triangles that share the nodes $N_1$, $N_2$, $N_3$ and $N_4$. Finally, let $\Omega_e$ be the local domain whose boundary coincides with the boundary of $\mathcal{T}_e$. We take $\Omega_e$ to have Dirichlet boundary conditions.
 
To define the improved local approximation, we bisect the triangles along the edge $e$ to obtain the grid $\hat{\mathcal{T}}$. We treat interface-base and base edges the same as this is a temporary change to the grid. The grid in Figure~\ref{fig:local_grid}, with node $N_5$, is an example of $\hat{\mathcal{T}}$. Let ~$\{(\hat{\bm{x}}_{i}, \hat{y}_{i}):i=1,2,\ldots,n_e\}$ be the subset of the observed data that sits inside $\Omega_e$. The local approximation~$\hat{s}=\hat{\bm{b}}(\bm{x})^{T}\hat{\bm{c}}$ is built by minimising the functional 
    \begin{equation} \label{eqn:auxi_minimiser}
		J_{\alpha}\left(\hat{\bm{c}},\hat{\bm{g}}_{1},\ldots,\hat{\bm{g}}_{d}\right) = \frac{1}{{n_e}}\sum^{{n_e}}_{i=1}\left(\hat{\bm{b}}(\hat{\bm{x}}_{i})^{T}\hat{\bm{c}}-\hat{y}_{i}\right)^{2}+\alpha \int_{{\Omega_e}} \sum_{k=1}^{d} \nabla \left(\hat{\bm{b}}(\hat{\bm{x}})^{T}\hat{\bm{g}}_{k}\right)\nabla \left(\hat{\bm{b}}(\hat{\bm{x}})^{T}\hat{\bm{g}}_{k}\right)\,d\bm{x}
    \end{equation}
subject to constraint~$\hat{L}\hat{\bm{c}} = \sum_{k=1}^{d}\hat{G}_{k}\hat{\bm{g}}_{k}$. Dirichlet boundary conditions applied to the local domain~${\Omega_e}$ are taken from the current global solution~$s$. So, if node $p$ sits on the boundary of $\Omega_e$,
	\begin{equation*}
		\hat{c}_{p} = c_{p}, \quad \hat{g}_{1,p} = g_{1,p},\quad\hat{g}_{2,p} = g_{2,p},\, \quad
        \cdots \quad
        \hat{g}_{d,p} = g_{d,p},\, \quad
        \hat{w}_{p} = w_{p}
    \end{equation*}
    where~${g}_{k,p}$ is the $p$-th entry of~$\bm{g}_k$. 
    
    The error indicator value~$\eta_{e}$ of triangular elements~$\tau_e$ is based on the energy norm as suggested by Babu{\v{s}}ka and Rheinboldt~\cite{babuvvska1978error}, where
    \begin{equation*} \label{eqn:auxiliary_indicator}
        \boxed{\eta_{e}^{2}= \int_{\tau_e}\left[
        \sum_{k=1}^{2}
        \left( \frac{\partial}{\partial x_{k}}\left(s-\hat{s}\right)\right)^{2}
        \right]\,d\bm{x}.}
    \end{equation*}

\begin{figure}
		\centering
		\includegraphics[width=0.35\textwidth]{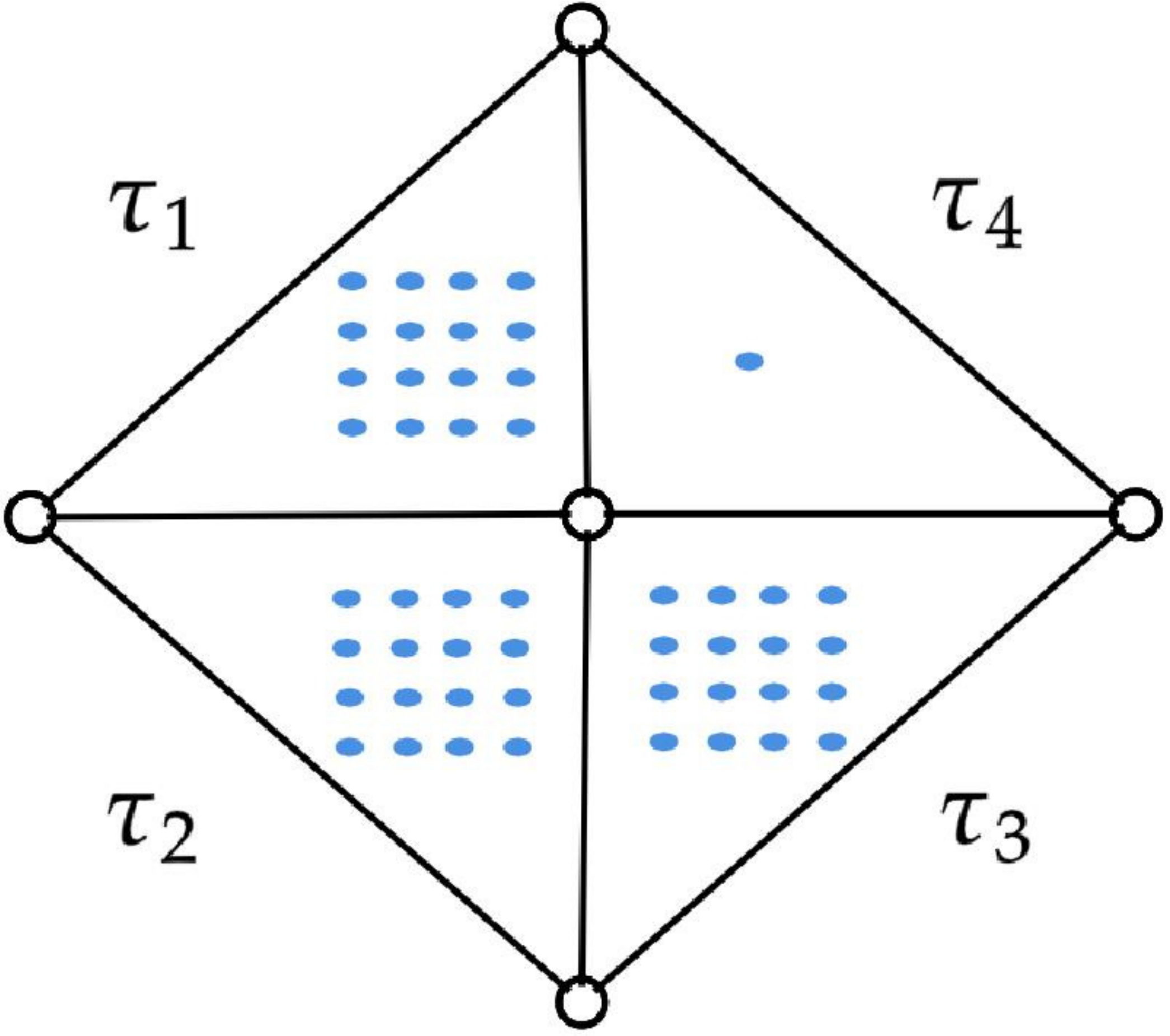}
		\caption{An example local domain with uneven data densities. Nodes and data points are represented as open circles and filled blue circles, respectively. The local domain contains four triangular elements, where elements~$\tau_{1}$,~$\tau_{2}$,~$\tau_{3}$ have markedly more data points than~$\tau_{4}$.}
		\label{fig:tri_data}
	\end{figure}

	Recall from Equation~\eqref{eqn:tpsfem_convergence_exact_1}, the expected error of the TPSFEM is written in terms of~$\alpha+h^{4}+d_{X}^{4}$. Consider the data distribution illustrated by data points in an example domain with triangles~$\tau_{1}$,~$\tau_{2}$,~$\tau_{3}$ and~$\tau_{4}$ in Figure~\ref{fig:tri_data}. We can calculate local values for~$d_{X}$ using the data given in triangles~$\tau_{1}$,~$\tau_{2}$ or~$\tau_{3}$ directly using the definition given in Section~\ref{sec:discrete}. However, the subdomain in~$\tau_{4}$ only has one data point and~$d_{X}$ is no longer locally defined. Consequently, the solution~$s$ in $\tau_{4}$, defined on~$\Omega$, can be very different to the solution~$\hat{s}$ in~$\tau_{4}$, defined on~${\Omega_e}$. Furthermore, the second term~$\frac{\sigma^{2}}{\alpha^{1\slash 2}}\frac{(h^{4}+d_{X}^{4})}{h^{2}}$ of equation~\eqref{eqn:tpsfem_convergence_exact_1} shows the effects of the noise may be more prominent in subdomains with few data points and thus, potentially, larger values of~$d_{X}$.
	More specifically, due to the potential difference in data distribution patterns, the behaviour of the solution in a local domain may vary dramatically compared to that in a global domain. In examples where the data is distributed evenly across the whole domain, such as image processing, the auxiliary problem error indicator is a viable option. However, if the data distribution is uneven, the behaviour of the solution can vary rapidly in the local subdomains. The consequences are readily seen in examples in Figure~\ref{fig:adaptive_grid} from the numerical experiments in Section~\ref{sec:result_crater} below.



\subsection{Residual-based error indicator}
\label{sec:residual}
    
	The residual-based error indicator estimates the energy norm of FEM approximation errors, which is represented as a combination of integrals of interior element residuals~$r_e$ and jumps of gradients~$j_e$ across the element boundaries~\cite{gratsch2005posteriori}. The TPSFEM has a different formulation compared to PDEs and the residual cannot be calculated directly. Our alternative approach is to work with locally refined grids, similar to the one shown in Figure \ref{fig:local_grid}.
 
Let $e$ be a base or interface base edge, and $\tau_e$ be the set of triangles sharing that edge. Using the procedure described in Section~\ref{sec:auxiliary}, construct a grid $\mathcal{T}_e$ defined on a local domain~${\Omega}_e$. Finally, take $\hat{\mathcal{T}}$ to be the result of bisecting edge $e$ in $\mathcal{T}_e$. 
 
 Let~$S\bm{p}=\bm{q}$ be the system defined in Equation~\eqref{eqn:system} where~$S$, $\bm{p}$ and $\bm{q}$ are the matrix, solution vector and right-hand side vector, respectively.  And let~$\hat S \hat{\bm{p}}=\hat{\bm{q}}$ be the corresponding system for the local grid $\hat{\mathcal{T}}$. Set $\mathbf{p}_e$ be the values of $\mathbf{p}$ restricted to grid $\mathcal{T}_e$. The residual~$\bm{r}_e$ is then calculated as
    \begin{equation*} \label{eqn:residual_R}
		\mathbf{r}_e = \hat{S}(I\bm{p}_e)-\hat{\bm{q}} 
    \end{equation*}  
    where $I$ is a linear interpolation operator that maps from functions defined on $\mathcal{T}_e$ to functions defined on $\hat{\mathcal{T}}$. In terms of Figure~\ref{fig:local_grid}, the values assigned to node $N_5$ are the average of the values assigned to nodes $N_1$ and $N_2$. When calculating the error indicator, only the first block of $\mathbf{r}_e$ that corresponds to the $\mathbf{c}$ coefficients is used. We use $\left.\mathbf{r}_e\right|_c$ to represent that block. 

 The jump~$\mathbf{j}_e$ for the edge $e$ shared between triangle~$\tau$ and triangle~$\tau'$ is 
    \begin{equation*} \label{eqn:residual_J}
        \mathbf{j}_e = \bm{n}\cdot \nabla s_{\tau} + \bm{n}' \cdot \nabla s'_{\tau},
    \end{equation*}    
    where~$s_{\tau}$ an~$s'_{\tau}$ are FEM approximations of $s$ on~$\tau$ and~$\tau'$, respectively; and~$\bm{n}$ and~$\bm{n}'$ are unit outward normal vectors to the edge $e$ in ~$\tau$ and~$\tau'$, respectively. If $e$ lies on Dirichlet boundaries, $\mathbf{j}_e=0$. If $e$ lies on Neumann boundaries, $\mathbf{j}_e=-\bm{n}\cdot \nabla s_{\tau}$. The error indicator value~$\eta_{e}$ of element~$\tau_e$ is computed as
    \begin{equation*} \label{eqn:residual_indicator}
        \boxed{\eta_{e}^{2} = c_{1}h_{e}^{2}\left\|\left.\mathbf{r}_e\right|_c\right\|_{L^{2}(\hat{\mathcal{T}})}^{2}+c_{2}h_{e}\left\|j_e\right\|_{L^{2}(\partial \tau_e)}^{2},}
    \end{equation*}     
    where~$h_{e}$ is the mesh size of $\tau_e$,~$c_{1}$ and~$c_{2}$ are constants that balance the contributions from the residual and jump, respectively~\citep{gratsch2005posteriori}. 


\subsection{Recovery-based error indicator}
\label{sec:recovery}
	
\begin{figure}
		\centering
		\includegraphics[width=0.5\textwidth]{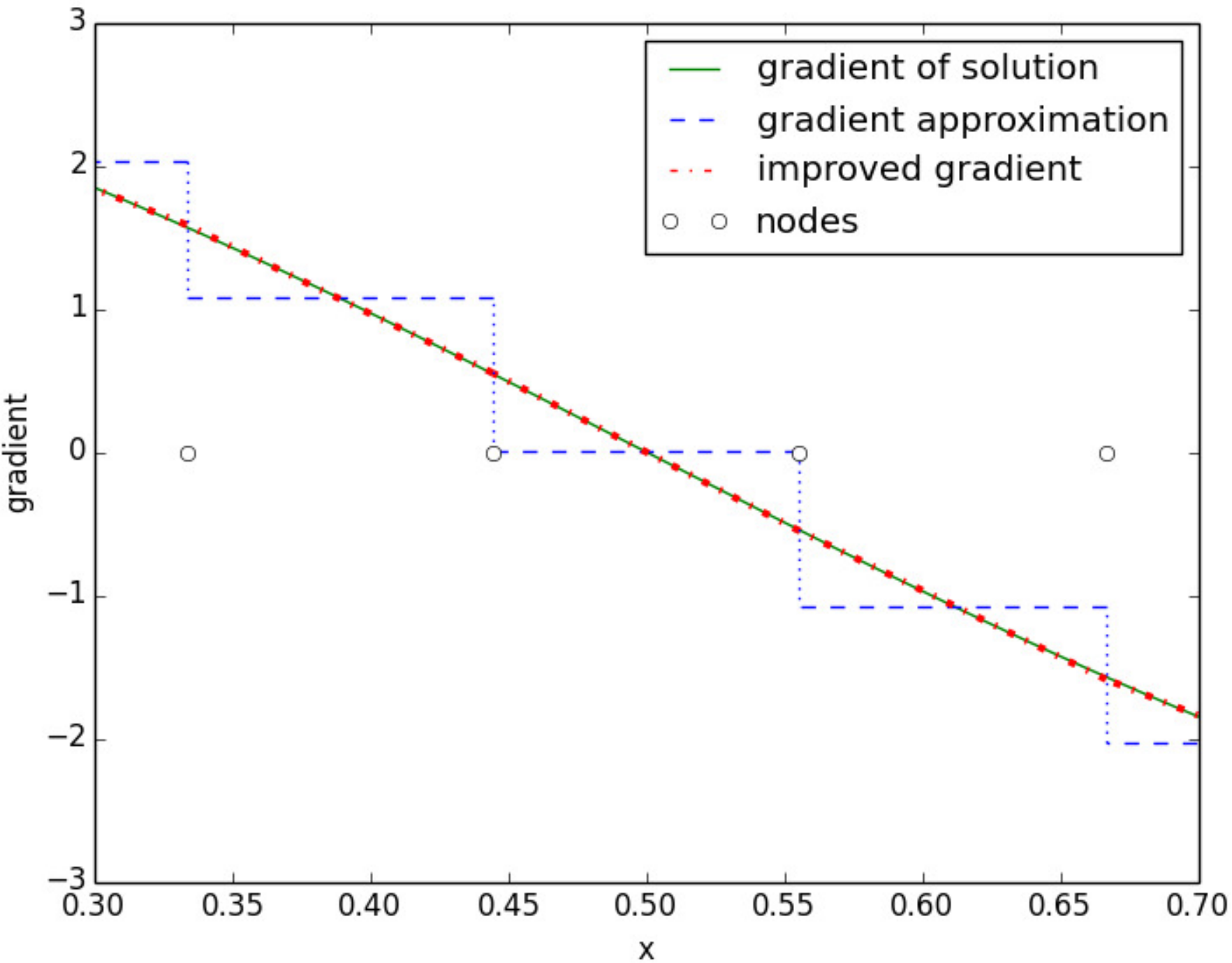}
		\caption{A one-dimensional example of improved gradients of FEM approximations.}
		\label{fig:improved}
	\end{figure}

	The recovery-based error indicator calculates the error norm by post-processing discontinuous gradients of the TPSFEM smoother~$s$ across inter-element boundaries~\cite{zienkiewicz1987simple}. The TPSFEM~$s$ is represented by piecewise linear basis functions with zero-order continuity~$C^{0}$ and the corresponding gradient approximation~$D^{1}s$ are piecewise constants as illustrated in Figure~\ref{fig:improved}. The gradient approximation~$D^{1}s$ is post-processed to obtain an improved discrete gradient approximation~$\hat{D}^{1}s$, which is closer to true gradients if the solution is smooth~\citep{gratsch2005posteriori}. 


	The improved gradient~$\hat{D}^{1}s$ is determined by the~$L^{2}$-projection
    \begin{equation*} \label{eqn:recovery_improved_project}
        \int_{\Omega}b_{p}\left(\hat{D}^{1}s-D^{1}s\right)\,d\mathbf{x} = 0, \qquad p =1,\ldots,m.
    \end{equation*}
	This leads to a system of equations
     \begin{equation*}
        \sum_{q=1}^{m}\int_{\Omega}b_{p}b_{q} \hat{D}^{1}s_{q} \,d\mathbf{x} = \int_{\Omega}b_{p}D^{1}s \,d\mathbf{x}, \qquad p =1,\ldots,m,
    \end{equation*}
where~$\hat{D}^{1}s_{q}$ is the improved gradient~$\hat{D}^{1}s$ evaluated at~$q$-th node. The error is then estimated using the difference between current and improved gradients and the error indicator value~$\eta_{e}$ of element~$\tau_e$ is thus set as
    \begin{equation*} \label{eqn:recovery_indicator}
        \boxed{\eta_{\tau_e}^{2}= \int_{\tau_e}\left(\hat{D}^{1}s-D^{1}s\right)^{2}\,d\bm{x}.}
    \end{equation*}	
	For PDEs, theoretical and experimental results show this error bound is asymptotically exact up to higher order terms, and it is shown to be effective for smooth functions~$f\in H^{1}(\Omega)$ with surface load~$g \in H^{1}(\partial\Omega_{N})$ on Neumann boundaries in Theorem 2.1 of~\cite{carstensen2001averaging}. The theory developed by Roberts et al.~\cite{roberts2003approximation} assumes functions~$f\in H^{2}(\Omega)$ with zero Neumann boundary conditions. Consequently, since the functions and boundary conditions for the TPSFEM are assumed to be smooth, we deduce this error indicator is worth considering for the TPSFEM.


	The process above only uses~$\bm{c}$ values to calculate improved gradient values~$\hat{D}^{1}s$. We calculate both~$\bm{c}$ and~$\bm{g}$ values when solving Equation~\eqref{eqn:system}, and~$\bm{g}$ values approximate the gradients of~$s$. Suggesting we could use $\bm{u}$, defined in Section~\ref{sec:discrete}, when calculating the error indicator.  Equations~\eqref{eqn:tpsfem_convergence_exact_1} and~\eqref{eqn:tpsfem_convergence_exact_2} show that the $L_2$ and $H_1$ norms of the TPSFEM are bounded by~$\frac{\sigma^{2}}{\alpha^{1\slash 2}}\frac{(h^{4}+d_{X}^{4})}{h^{2}}$ and~$\frac{\sigma^{2}}{\alpha}\frac{(h^{4}+d_{X}^{4})}{h^{2}}$, respectively, when the variance~$\sigma^{2}$ of noise is nonzero. Given that generally $\alpha \le 1$, this suggests the gradient is more sensitive to noise. This is illustrated by a numerical experiment in Section 5.5 of~\cite{fang2021error}. Therefore, we chose to compute improved gradients using~$\bm{c}$ values for the recovery-based error indicators.



\subsection{Norm-based error indicator}
\label{sec:norm}

	The norm-based error indicator uses an error bound on the~$L_{\infty}$ norm of the interpolation error~\cite{sewell2012analysis}. Assuming the solution~$f$ has bounded second-order derivatives, Sewell demonstrated the error of a two-dimensional FEM grid with first-order basis functions is bounded in terms of second-order derivatives~$D^{2}_{\text{max}}f$, where
    \begin{equation*} \label{eqn:indicator_norm_Dn}
        D^{2}_{\text{max}}f\left(x_{1},x_{2}\right) = \max_{i+j=2}\left|\partial^{2}f\left(x_{1},x_{2}\right)/\partial x_{1}^{i}x_{2}^{j}\right|.
    \end{equation*}
	He showed that a near-optimal grid with piecewise linear approximation will have~$\int_{\tau_{i}}D^{2}_{\text{max}}f\,d\bm{x}$ equally distributed over all elements~$\tau_{i}$. This error indicator identifies regions where the solution changes rapidly (high integral values) and resolves them using finer elements (smaller~$h$) to achieve the required accuracy. 

	In contrast to the PDEs, the accuracy of the TPSFEM smoother~$s$ is also affected by~$\alpha$ and~$d_{X}$ as shown in Equation~\eqref{eqn:tpsfem_convergence_exact_1}. If either~$\alpha$ or~$d_{X}$ is large, the TPSFEM will interpolate a flat surface and this error indicator will refine evenly to distribute~$\int_{\tau_{i}}D^{2}_{\text{max}}s\,d\bm{x}$ across the domain. If~$\alpha$ and~$d_{X}$ is sufficiently small, this error indicator will concentrate on reducing~$h$ on regions with high expected errors to produce an efficient grid. 


	For two-dimensional grids, we approximate second-order derivatives with respect to~$x_{1}$ and~$x_{2}$ on each node and~$D^{2}_{\text{max}}s$ is calculated as
    \begin{equation*}
        D^{2}_{\text{max}}s \approx \max\left\{\left|D_{x_{1}x_{1}}s\right|,\left|D_{x_{1}x_{2}}s\right|,\left|D_{x_{2}x_{2}}s\right|\right\}.
    \end{equation*}
	The error indicator value assigned to edge $e$ is defined as
    \begin{equation*} \label{eqn:indicator_norm_indicator_2d}
        \boxed{\eta_{e} = \int_{\tau_e}D^{2}_{\text{max}}s\,d\bm{x},}
    \end{equation*}
    where $\tau_e$ is the set of triangles sharing the edge $e$.
	The error bound was derived under the assumption that the model problem has bounded second-order derivatives. Since the smoother~$s$ of the TPSFEM in Equation~\eqref{eqn:tpsfem} is minimised with respect to its second-order derivatives, the error bound naturally holds for the TPSFEM.
	

\section{Numerical experiment}
\label{sec:experiment}

	We evaluated adaptive refinement for the TPSFEM using three data sets: data generated using the peaks function from MATLAB, the 2000 Multibeam sonar survey of Crater Lake~\cite{bathymetric2000lake} and a bathymetric sounding survey of a coastal region~\cite{bathymetric2021coastal}. The first data set consists of 62,500 data points uniformly distributed inside the~$[-2.4,2.4]^{2}$ region of a~$[-3,3]^{2}$ finite element domain and is modelled by
	\begin{equation} \label{eqn:peak}
		y =  3\left(1-x_1\right)^2\exp\left(-\left(x_1^2\right) - \left(x_2+1\right)^2\right) - 10\left(\frac{x_1}{5} - x_1^3 - x_2^5\right)\exp\left(-x_1^2-x_2^2\right)- \frac{1}{3}\exp\left(-\left(x_1+1\right)^2 - x_2^2\right)+\epsilon.
	\end{equation}
	The peaks function has three local maxima and three local minima around the centre of the domain, and small flat surfaces near four corners. The data is perturbed by Gaussian noise~$\epsilon$ with mean~$0$ and standard deviation~$0.02$. The performance of the TPSFEM on this model problem is discussed in Section~\ref{sec:result_model}. The Crater Lake survey comprises 12,936,068 data points consisting of latitude, longitude and elevations above sea level from a multibeam bathymetric survey of Crater Lake, Oregon in the United States.
 
    The survey collected over 16 million soundings and the data portrays the bottom of Crater Lake at a spatial resolution of 2 meters. The Crater Lake data is used to compare the performance of the TPSFEM with Dirichlet and Neumann boundaries in Section~\ref{sec:result_crater}. The Coastal Region survey comprises 48,905 data points consisting of North American Datum of 1983\footnote{North American Datum, https://en.wikipedia.org/wiki/North\_American\_Datum} and modern readings of mean low water in meters. The mean low water is the average of all the low water heights observed over the National Tidal Datum Epoch. These measurements were taken to provide an estimate of historical bathymetry for the Mississippi-Alabama coastal region to aid geologic and coastal hazard studies. In Section~\ref{sec:result_coast}, the Coastal Region data forms the basis of an investigation into the performance of the TPSFEM in a non-rectangular domain. 
	The performance of uniform and adaptive grids is measured using the root-mean-square error (RMSE) versus the number of nodes, where~$\text{RMSE} = \big(\frac{1}{n}\sum_{i=1}^{n}\left(s(\bm{x}_{i})-y_{i}\right)^{2}\big)^{\frac{1}{2}}$.


	\begin{figure}
		\centering
		\begin{subfigure}[b]{0.45\textwidth}
		\centering
		\includegraphics[width=\textwidth]{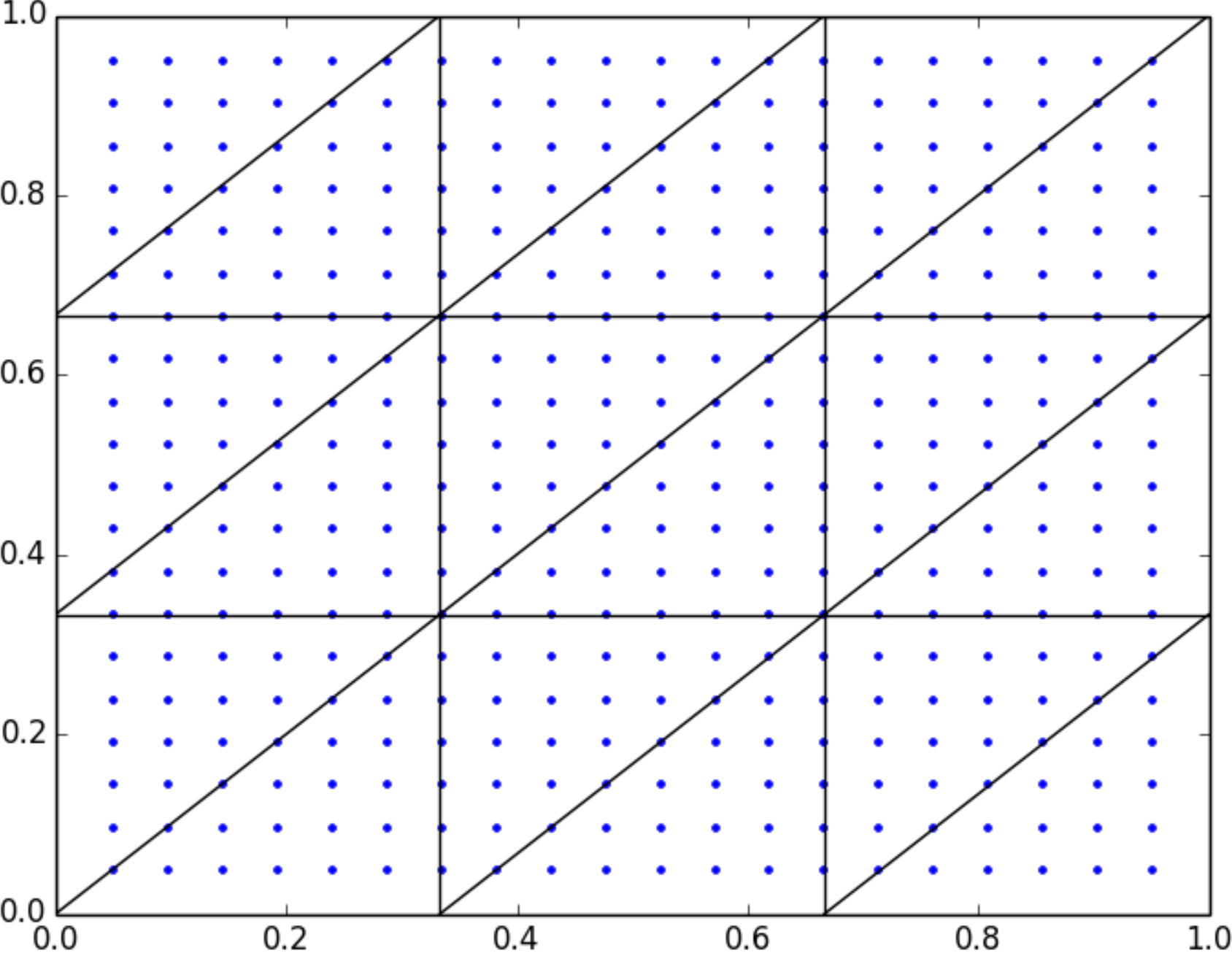}
		\caption{}
		\end{subfigure}
		\hspace{0.5cm}
		\begin{subfigure}[b]{0.45\textwidth}
		\centering
		\includegraphics[width=\textwidth]{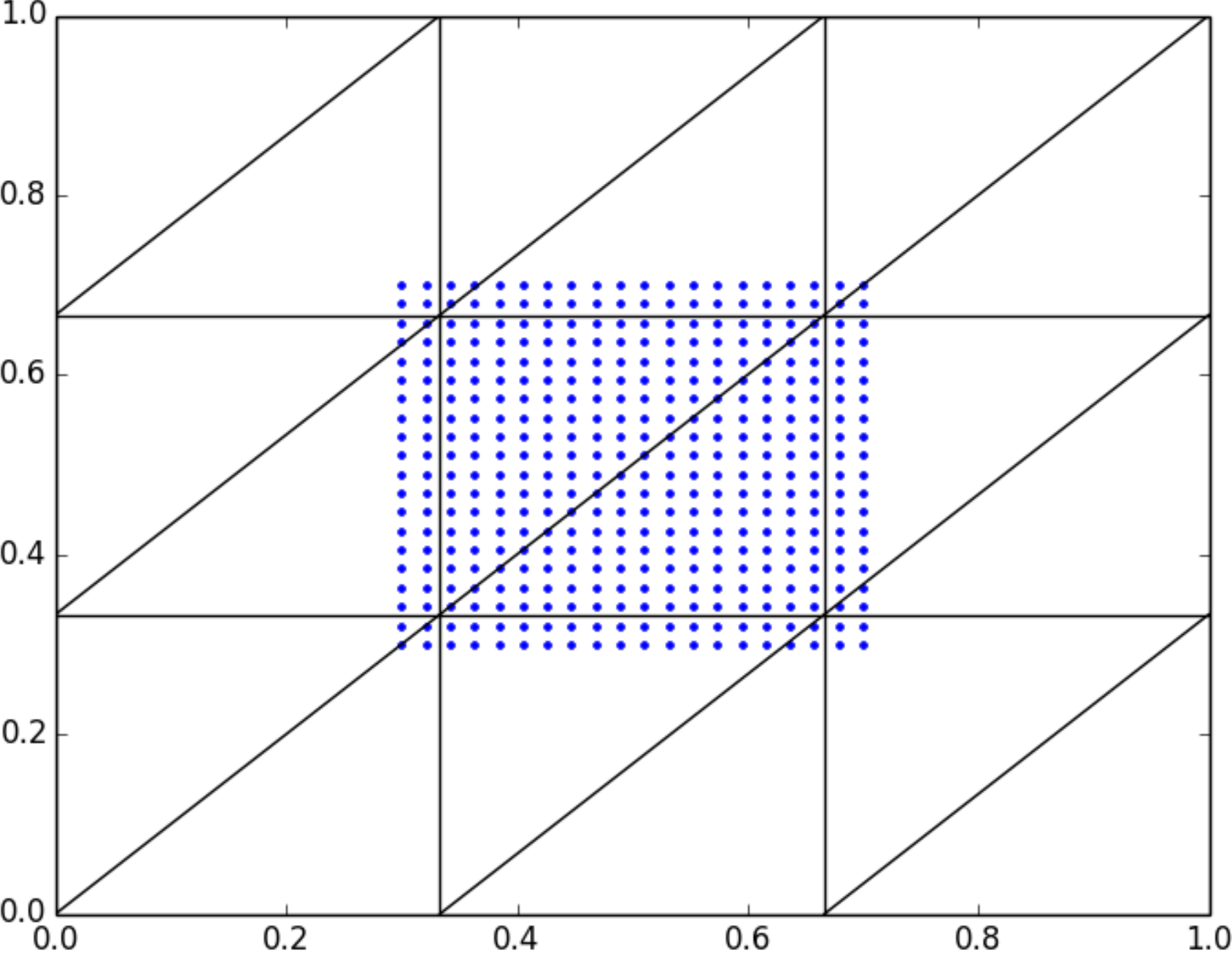}
		\caption{}
		\end{subfigure}
		\caption{Example FEM grids with data points (a) close to boundaries; and (b) far away from boundaries. Data points are represented as blue points.}
		\label{fig:fem_grid_data}
	\end{figure}

	Dirichlet boundary conditions are set using the $\bm{h}_1, \cdots, \bm{h}_4$ vectors in Equation~\eqref{eqn:system}. If the exact boundary conditions are known, the data can be placed near the boundaries of the domain as illustrated in Figure~\ref{fig:fem_grid_data}(a). In general, however, such an approach is not possible. The theoretical analysis given in~\cite{roberts2003approximation} assumes zero Neumann boundary conditions. In other words, it assumes the solution is flat near the boundary, which is not necessarily the case. Consequently, if the boundary conditions are not known, the data points are placed well inside the domain as shown in Figure~\ref{fig:fem_grid_data}(b). In summary, when it is possible to approximate the boundary conditions, it is advantageous to use Dirichlet boundary conditions. Otherwise, the FEM grid is extended well past the boundary of the data set and Neumann boundary conditions are applied.



	We observe in Figure~\ref{fig:fem_grid_data}(b) that the distribution of the data points can vary greatly from element to element. Some of the elements near the boundary only have 1 or 2 data points, whereas those in the interior have many data points. As discussed in Section~\ref{sec:experiment}, this may affect the accuracy of the local approximation of the auxiliary problem error indicator.


\subsection{Results: peaks function}
\label{sec:result_model}

\begin{figure}
		\centering
		\begin{subfigure}[b]{0.45\textwidth}
		\centering
		\includegraphics[width=\textwidth]{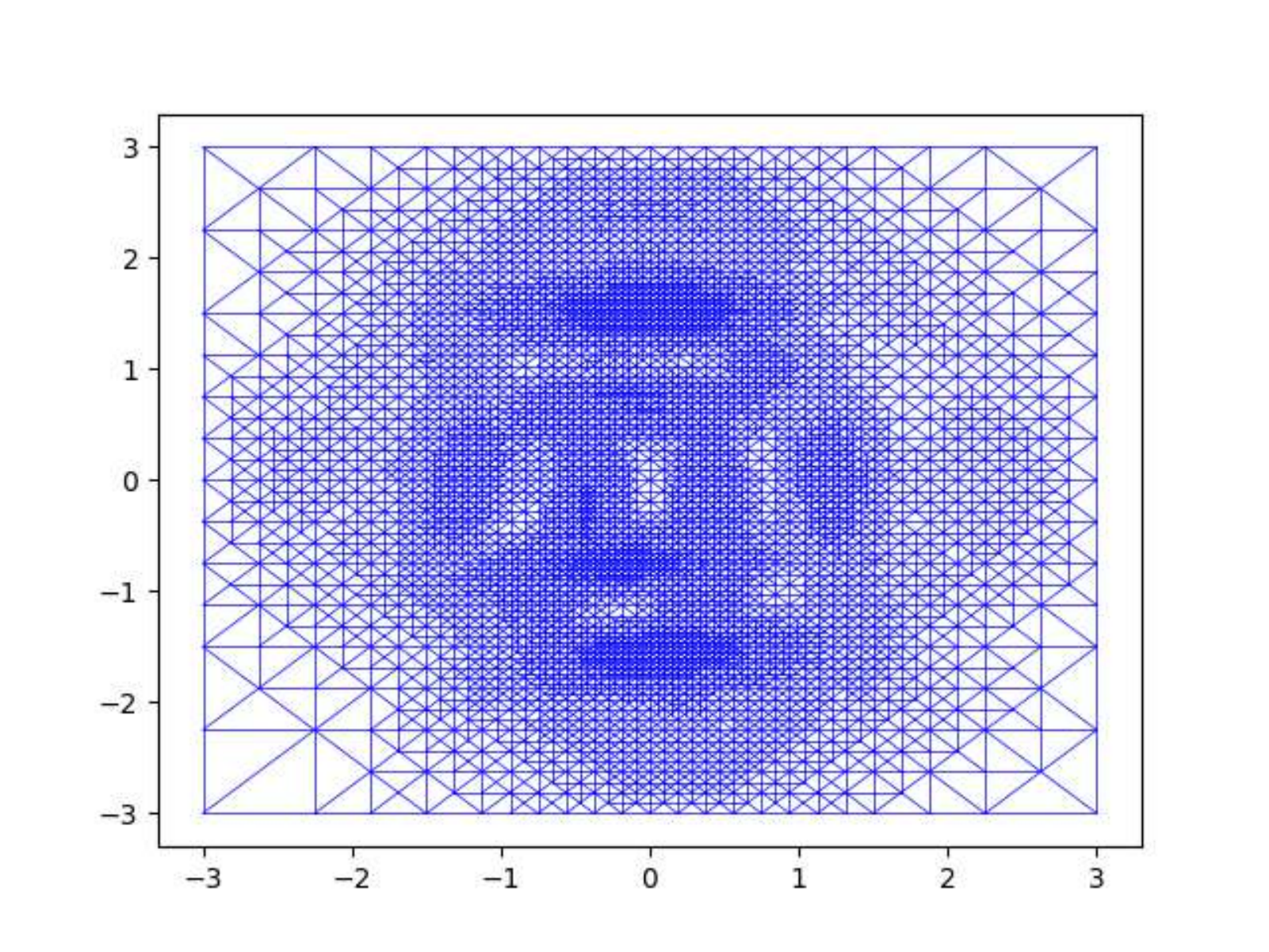}
		\caption{}
		\end{subfigure}
		\hspace{0.5cm}
		\begin{subfigure}[b]{0.45\textwidth}
		\centering
		\includegraphics[width=\textwidth]{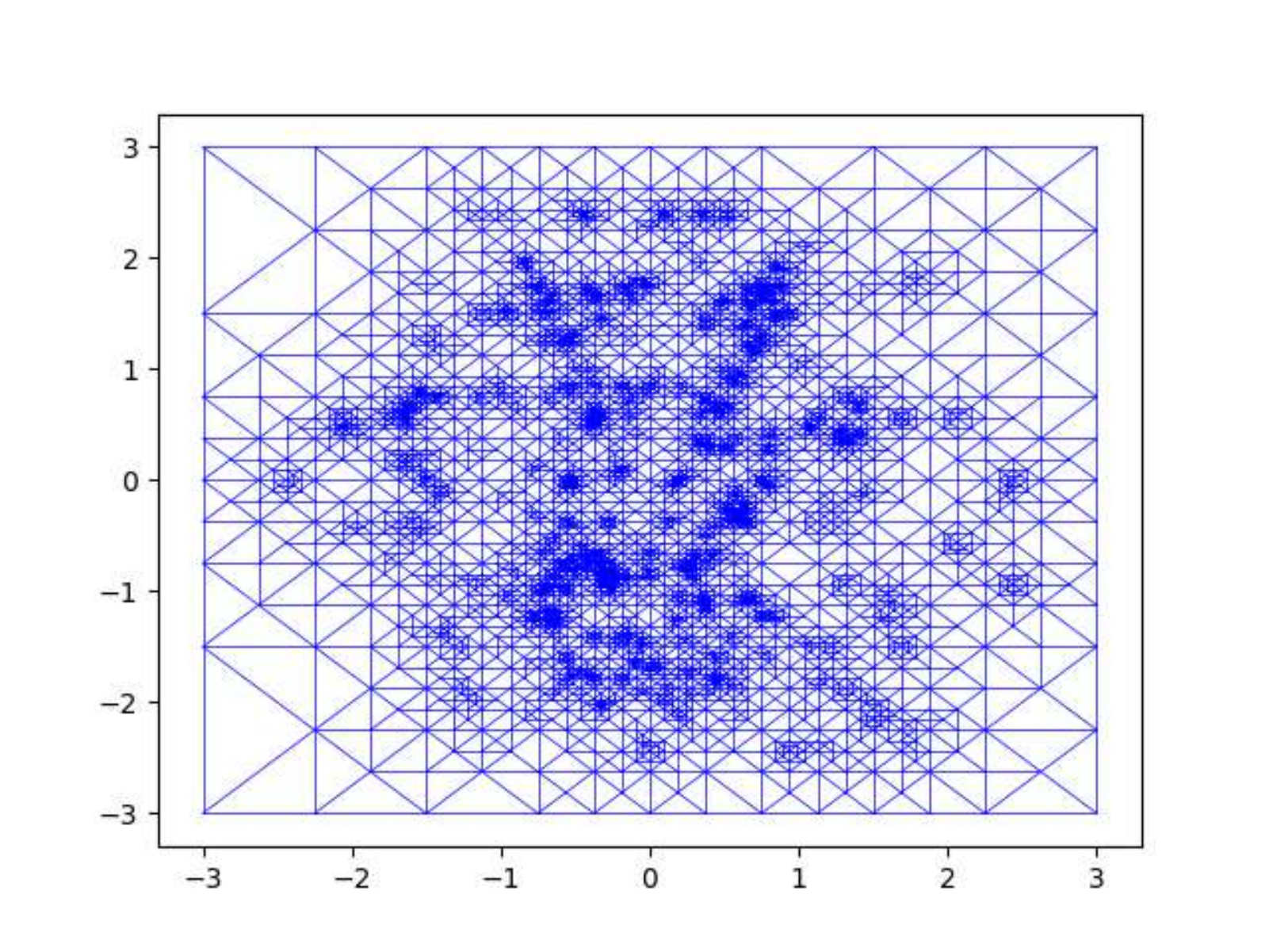}
		\caption{}
		\end{subfigure}
		\caption{Adaptive grids using (a) norm-based error indicator; and (b) regression metric.}
		\label{fig:grid_rmse}
	\end{figure}

	Two adaptive grids for the peaks function data in Equation~\eqref{eqn:peak} generated using norm-based error indicator and the regression metric are shown in Figures~\ref{fig:grid_rmse}(a) and~\ref{fig:grid_rmse}(b), respectively. The description of these two error indicators is given in sections~\ref{sec:regression} and~\ref{sec:norm}. 
 Refinement in the grid shown in Figure~\ref{fig:grid_rmse}(a) is concentrated at the peaks in the centre of the domain. There is some refinement away from the centre of the domain. This is partly due to the use of interface base-edges to avoid hanging nodes as described in Section~\ref{sec:adaptive}. It is also due to the fact that once the grid is fine enough in the centre, the error indicator will pick up some of the regions away from the centre of the domain.
 The norm-based error indicator has identified regions that change more rapidly. In contrast, the adaptive grid generated using the regression metric as an error indicator shown in Figures~\ref{fig:grid_rmse}(b) contains over-refinement in many small regions. These results further confirm that the regression metric is sensitive to noise as discussed in Section~\ref{sec:regression}. Consequently, they are not considered for further testing.


\begin{figure}
		\centering
		\begin{subfigure}[b]{0.45\textwidth}
		\centering
		\includegraphics[width=\textwidth]{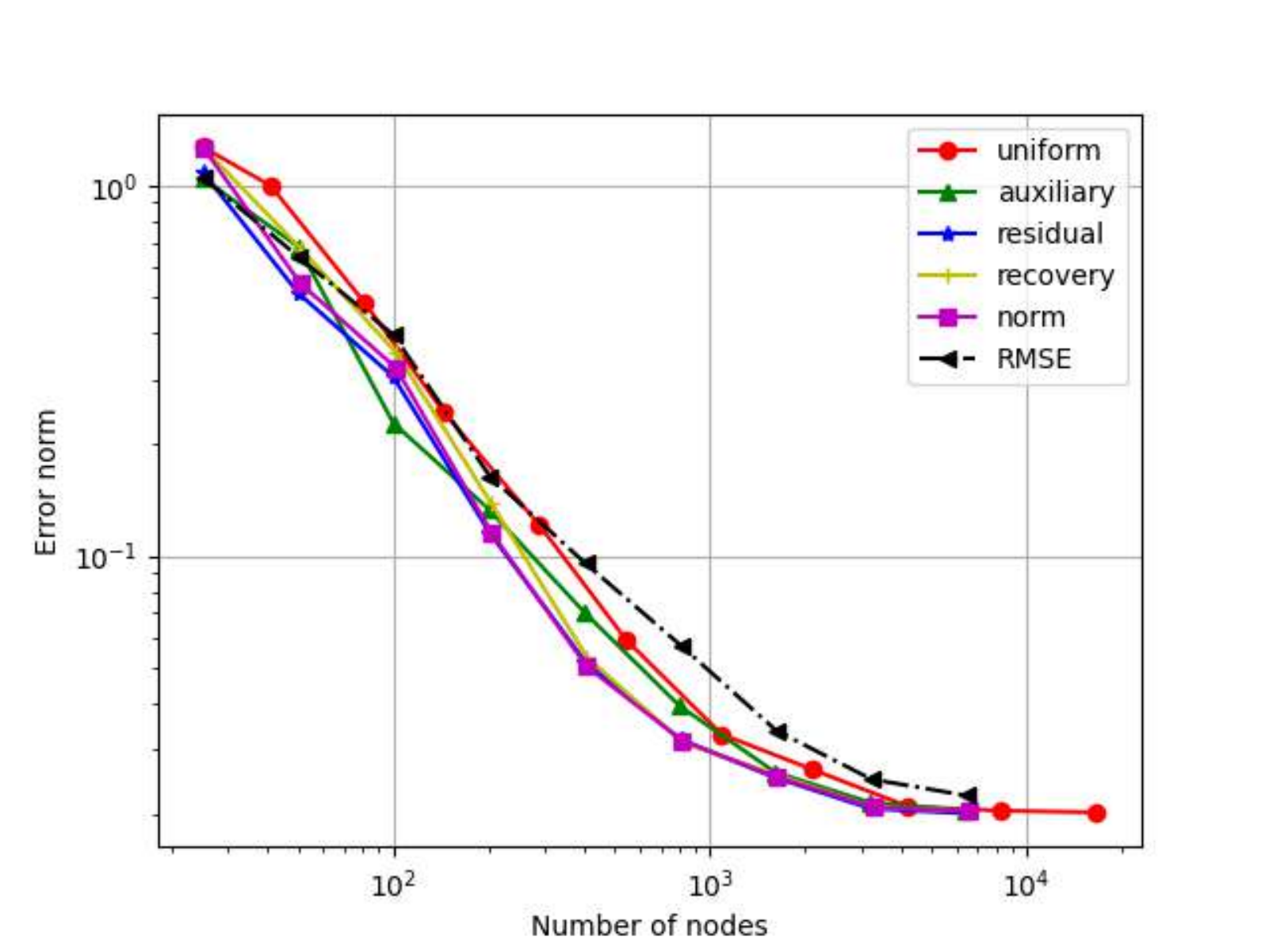}
		\caption{}
		\end{subfigure}
		\hspace{0.5cm}
		\begin{subfigure}[b]{0.45\textwidth}
		\centering
		\includegraphics[width=\textwidth]{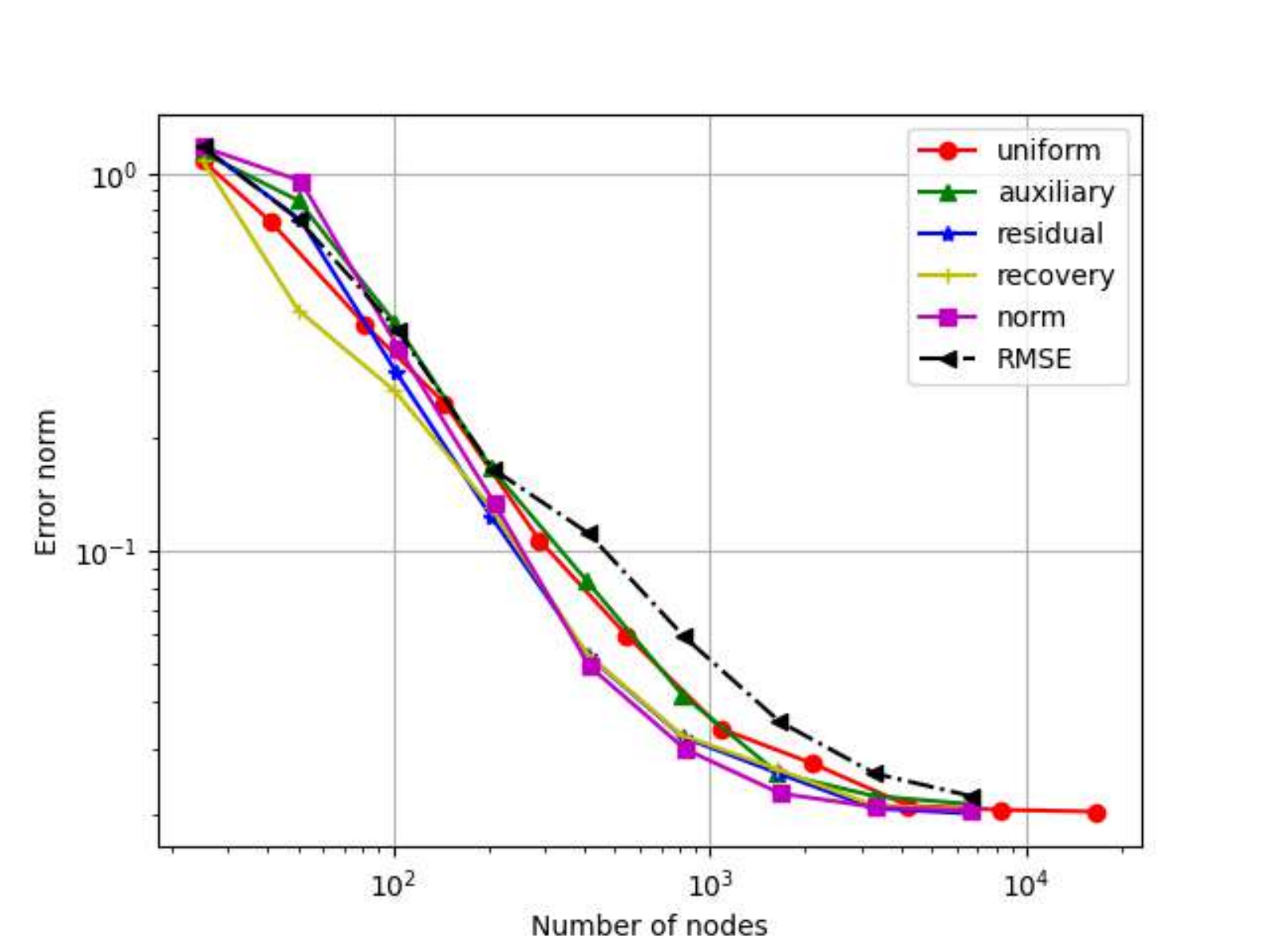}
		\caption{}
		\end{subfigure}
		\caption{RMSE of TPSFEM using uniform and adaptive grids with (a) Dirichlet boundaries; and (b) Neumann boundaries. The six convergence curves correspond to TPSFEM generated using uniform refinement, adaptive refinement with auxiliary problem, residual-based, recovery-based and norm-based error indicators and RMSE as an error indicator, respectively.}
		\label{fig:rmse_model}
	\end{figure}

	The convergence of the RMSE of the TPSFEM using uniform and adaptive grids for the peaks function data with Dirichlet and Neumann boundaries is shown in Figures~\ref{fig:rmse_model}(a) and~\ref{fig:rmse_model}(b), respectively. The initial square grids contain 25 nodes and are refined using uniform and adaptive refinement for at most 10 and 8 iterations, respectively. The adaptive grids produced by the four error indicators have similar error convergence rates for both Dirichlet and Neumann boundaries, which are significantly higher than that of the uniform grid. While the adaptive grids produced by the regression metric achieve similar error convergence for the first three iterations, the error convergence slows down, and the final grids have the highest RMSE among the evaluated indicators.

We list the RMSE, root-mean-square percentage errors (RMSPE), maximum residual errors (MAX), and number of nodes of the final grids in Table~\ref{tab:metric_model}. The RMSPE is used to compare the performance of the three data sets, which have various $\bm{y}$ value ranges. It is calculated as~$\text{RMSPE} = \big(\frac{1}{n}\sum_{i=1}^{n}\left((s(\bm{x}_{i})-y_{i})/{y_{\max}}\right)^{2}\big)^{\frac{1}{2}}$, where~$y_{\max}=\max(\bm{y})$. We also include the runtime, which is measured in seconds, for system-solving, system-building and the error indicators for the final grids of the TPSFEM in Table~\ref{tab:metric_model}. The construction of Equation~\eqref{eqn:system} needs to scan the data, which is affected by the number of nodes and data sizes. The error indicator is also critical to the Algorithm~\ref{alg:tpsfem_adaptive}. They are displayed and compared to illustrate the effectiveness of adaptive refinement.

The adaptive grids using the four error indicators achieve similar RMSE, RMSPE and MAX using about $39.22\%$ of the number of nodes, $46.39\%$ of the solve time and $44.80\%$ of the build time compared to the uniform grids. While the adaptive grids produced by the regression metric have similar MAX to some other adaptive grids, their RMSE and RMSPE are significantly higher and are therefore less efficient. Note that the runtime of the auxiliary problem and residual-based error indicators is markedly higher than that of the recovery-based and norm-based error indicators since they both access local data. While regression metric also uses data, it has a simpler formulation and took the least time in the last iteration compared to the others.

	\begin{table}
		\centering
		\caption{Regression metrics of the TPSFEM for peak function}
		\label{tab:metric_model}
		\begin{tabular}{llllllll}
		\hline\noalign{\smallskip}
		Boundary & Metric & Uniform & Auxiliary & Residual & Recovery & Norm & Regression \\
		\noalign{\smallskip}\hline\noalign{\smallskip}
		Dirichlet & RMSE & 0.021 & 0.021 & 0.021 & 0.021 & 0.021 & 0.023 \\
		 & RMSPE & $2.50\times 10^{-3}$ & $2.52\times 10^{-3}$ & $2.47\times 10^{-3}$ & $2.51\times 10^{-3}$ & $2.51\times 10^{-3}$ & $2.77\times 10^{-3}$\\
		 & MAX & 0.093 & 0.093 & 0.11 & 0.094 & 0.093 & 0.11 \\
		 & \# nodes & 16641  & 6404 & 6456 & 6561 & 6536 & 6468 \\
		 & Solve & 0.41 & 0.21 & 0.21 & 0.20 & 0.19 & 0.21 \\
          & Build & 24.22 & 11.28 & 11.61 & 10.60 & 10.65 & 10.83 \\
		 & Indicator &  & 286.19 & 270.33 & 134.57 & 86.45 & 54.97 \\
		\noalign{\smallskip}\hline\noalign{\smallskip}
		Neumann & RMSE & 0.021 & 0.022 & 0.020 & 0.021 & 0.021 & 0.022\\
		 & RMSPE & $2.51\times 10^{-3}$ & $2.63\times 10^{-3}$ & $2.48\times 10^{-3}$ & $2.53\times 10^{-3}$ & $2.51\times 10^{-3}$ & $2.76\times 10^{-3}$\\
		 & MAX & 0.095 & 0.10 & 0.10 & 0.11 & 0.10 & 0.11 \\
		 & \# nodes & 16641 & 6538 & 6564 & 6496 & 6660 & 6627 \\
		 & Solve & 0.49 & 0.22 & 0.22 & 0.22 & 0.20 & 0.22 \\
          & Build & 24.18 & 11.15 & 10.62 & 10.42 & 10.42 & 10.83 \\
		 & Indicator &  & 287.03 & 270.28 & 132.06 & 88.23 & 52.04  \\
		\noalign{\smallskip}\hline
		\end{tabular}
	\end{table}


\subsection{Results: Crater lake}\label{sec:result_crater}
 
\begin{figure}
		\centering
		\includegraphics[width=0.45\textwidth]{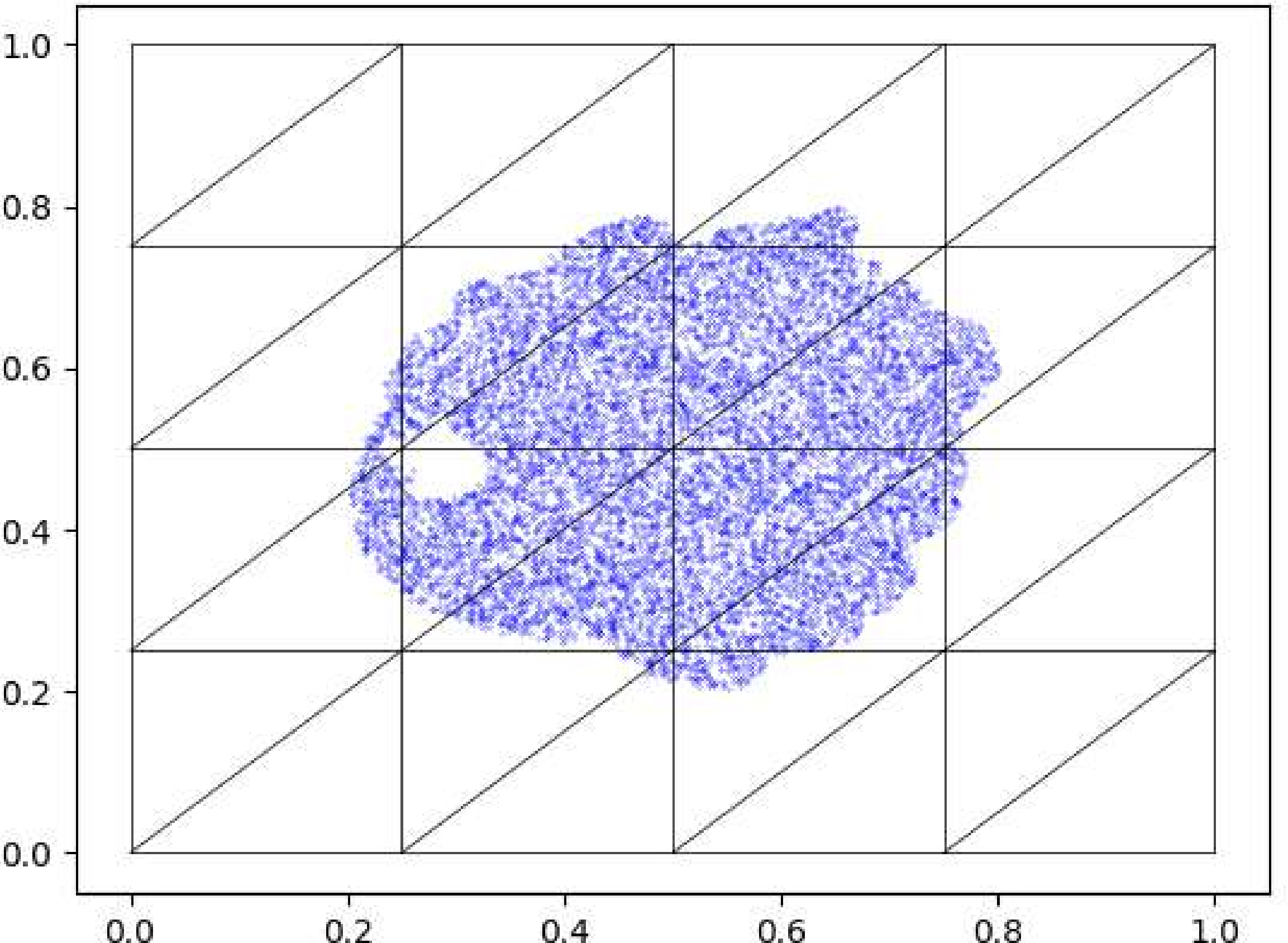}
		\caption{Randomly sampled data points in initial FEM grids for Crater Lake data. Data points are represented as blue dots.}
		\label{fig:grid_data_lake}
	\end{figure}

\begin{figure}
		\centering
		\begin{subfigure}[b]{0.45\textwidth}
		\centering
 		\includegraphics[width=\textwidth]{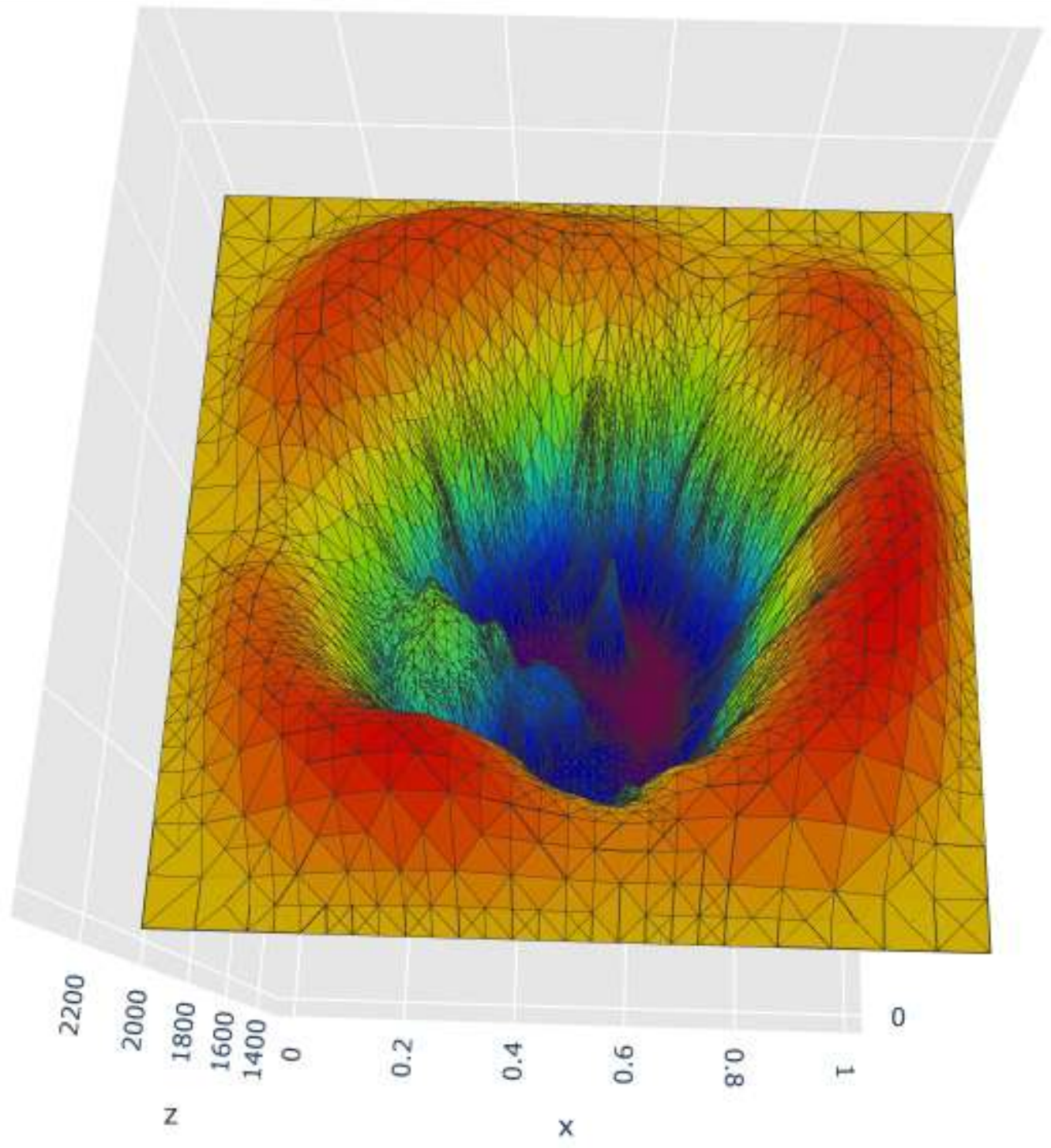}
		\caption{}
		\end{subfigure}
		\hspace{0.5cm}
		\begin{subfigure}[b]{0.45\textwidth}
		\centering
		\includegraphics[width=\textwidth]{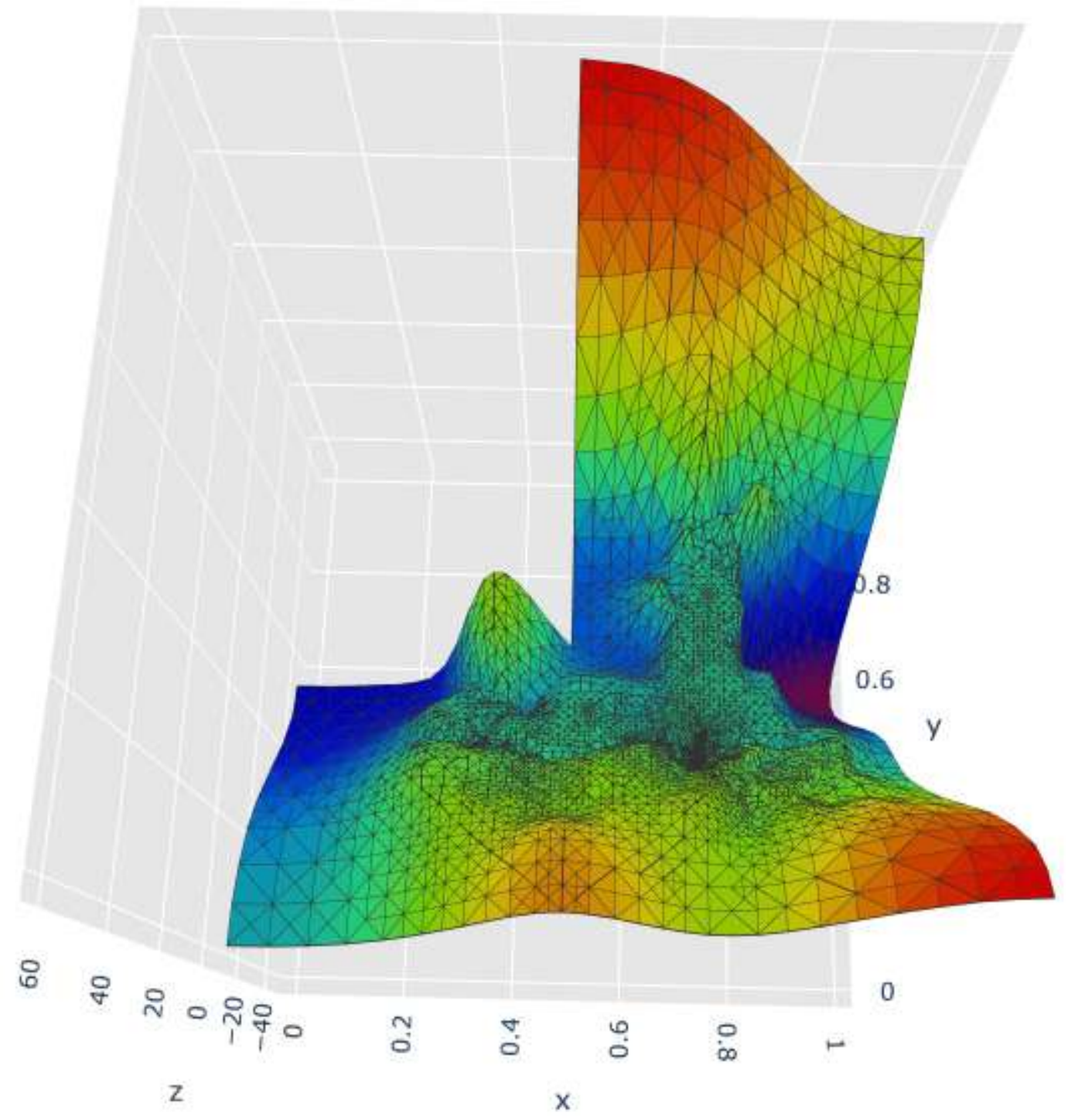}
		\caption{}
		\end{subfigure}
		\caption{Triangulated surfaces of TPSFEM for (a) Crater Lake data; and (b) Coastal Region data.}
		\label{fig:surface_3D}
    \end{figure}

We scaled and fitted the Crater Lake data inside the~$[0.2,0.8]^{2}$ region and used~$[0,1]^{2}$ as the finite element domain as illustrated in Figure~\ref{fig:grid_data_lake}. Note that while the FEM domain is square, the data has a complicated boundary. We tested the performance for the Crater Lake data using both Dirichlet and zero Neumann boundary conditions. The data values increase from about 1289.1 in the interior of the lake to about 1880.9 in the shoreline. As such, we would expect $s$ to smoothly increase from the interior of the lake to the boundary of the FEM domain. Furthermore, for this particular data set, the data values near the boundary are roughly constant. So, when applying Dirichlet boundary conditions, we set~$s=2100$, $u_{1}=0$ and $u_{2}=0$ along the boundary of the domain. The example surface plots of the TPSFEM for the Crater Lake data and Coastal Region data are shown in Figures~\ref{fig:surface_3D}(a) and~\ref{fig:surface_3D}(b), respectively. The corresponding interactive 3D plots are also available in the TPSFEM program repository.

\begin{figure}
		\centering
		\begin{subfigure}[b]{0.45\textwidth}
		\centering
		\includegraphics[width=\textwidth]{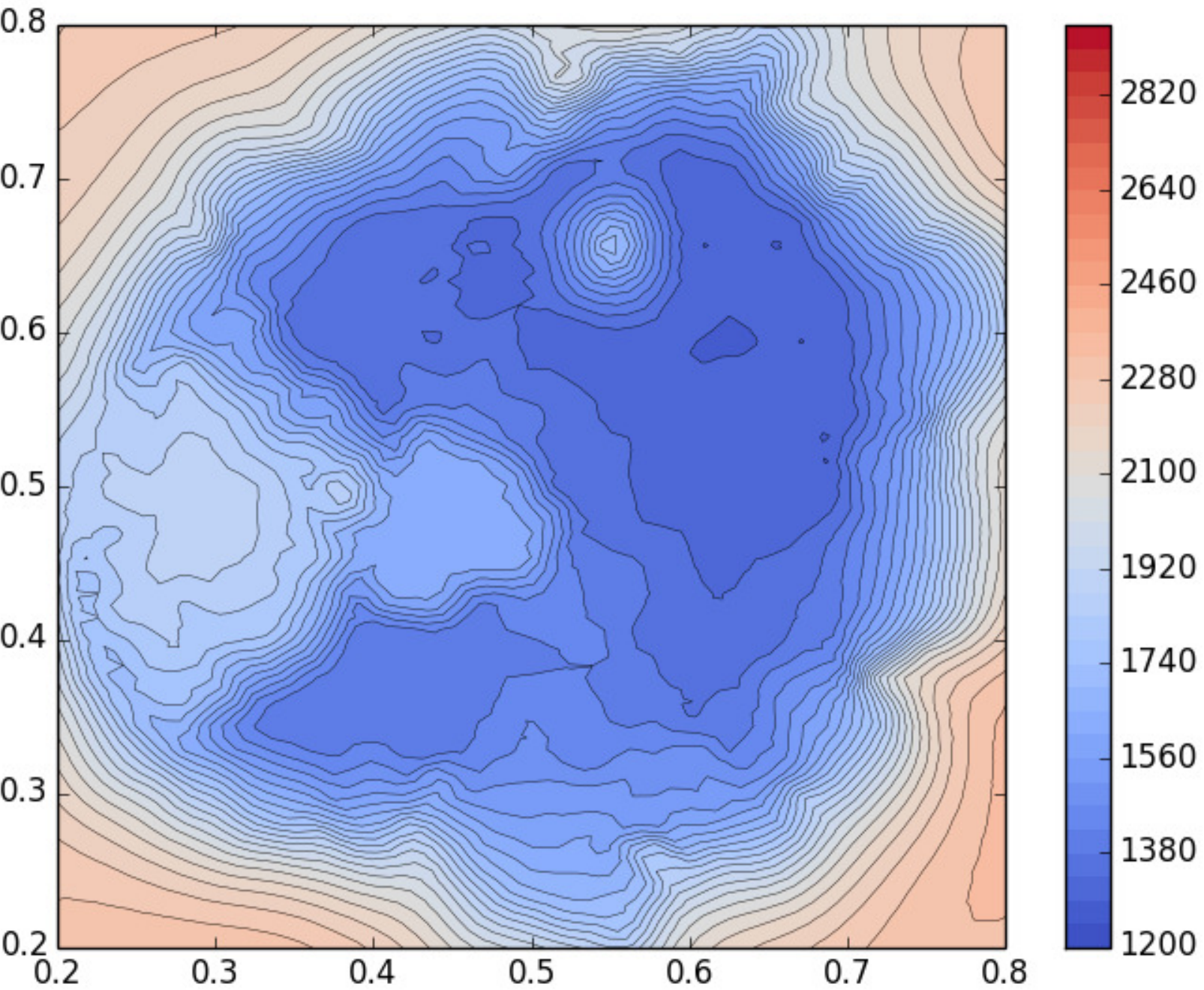}
		\caption{}
		\end{subfigure}
		\hspace{0.5cm}
		\begin{subfigure}[b]{0.45\textwidth}
		\centering
		\includegraphics[width=\textwidth]{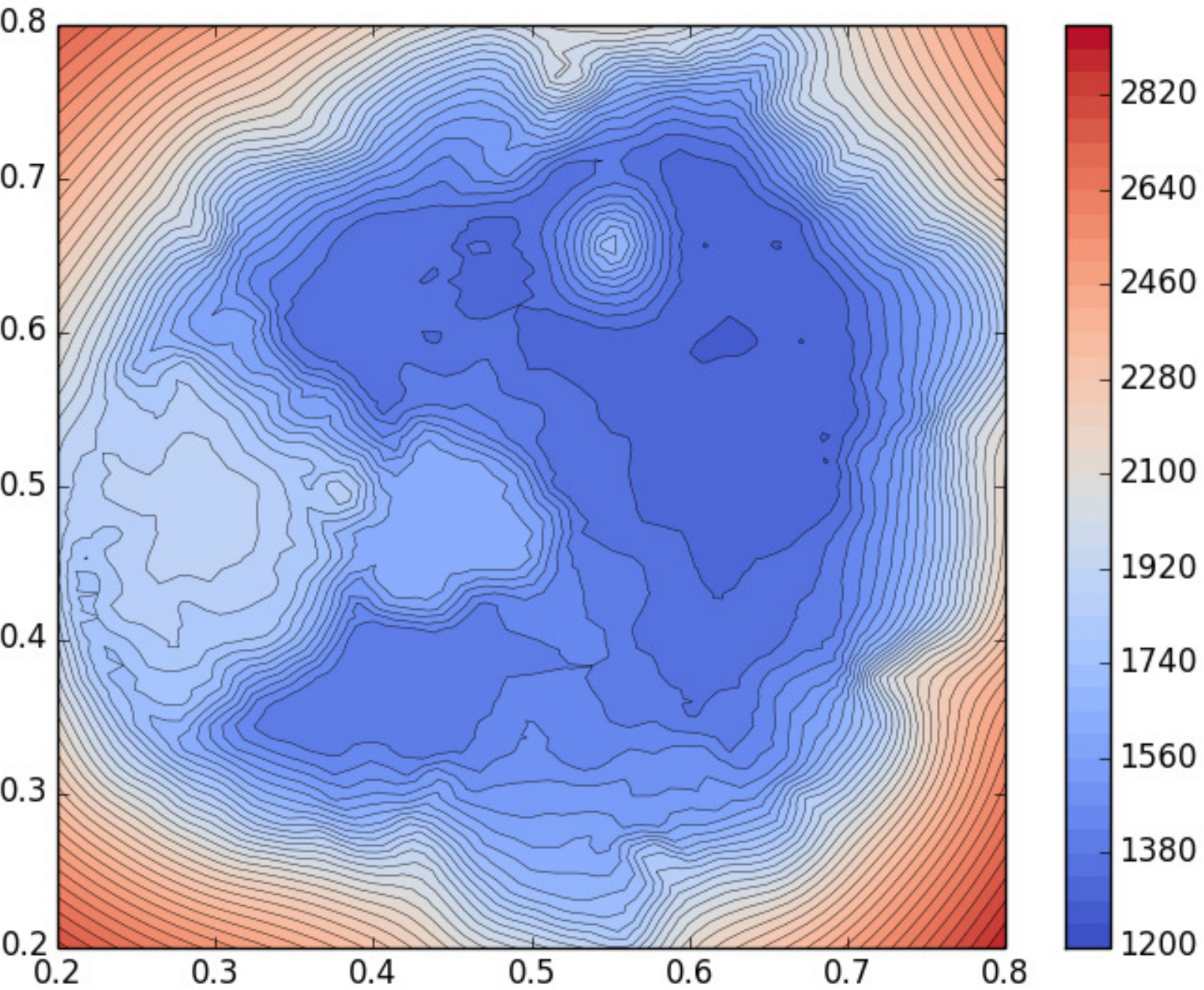}
		\caption{}
		\end{subfigure}
		\caption{Contour plots of $s$ using a uniform grid of 8,321 nodes with (a) Dirichlet boundaries; and (b) Neumann boundaries.}
		\label{fig:tpsfem_contour}
	\end{figure}

	Two contour maps of $s$ with Dirichlet and Neumann boundaries for the Crater Lake data are shown in Figures~\ref{fig:tpsfem_contour}(a) and~\ref{fig:tpsfem_contour}(b), respectively.  The scaling for Figure~\ref{fig:tpsfem_contour} was determined by the minimum/maximum values of the $\bm{c}$ coefficient. Both results are similar and capture bumps and ridges of the lake in regions with densely populated data. In contrast, they behave differently near corners of the domain where there is no data and the boundary conditions are different. The TPSFEM smoother built using Neumann boundaries has a steeper ascent near the boundary compared to the one with Dirichlet boundaries.

\begin{figure}
		\centering
		\begin{subfigure}[b]{0.45\textwidth}
		\centering
		\includegraphics[width=\textwidth]{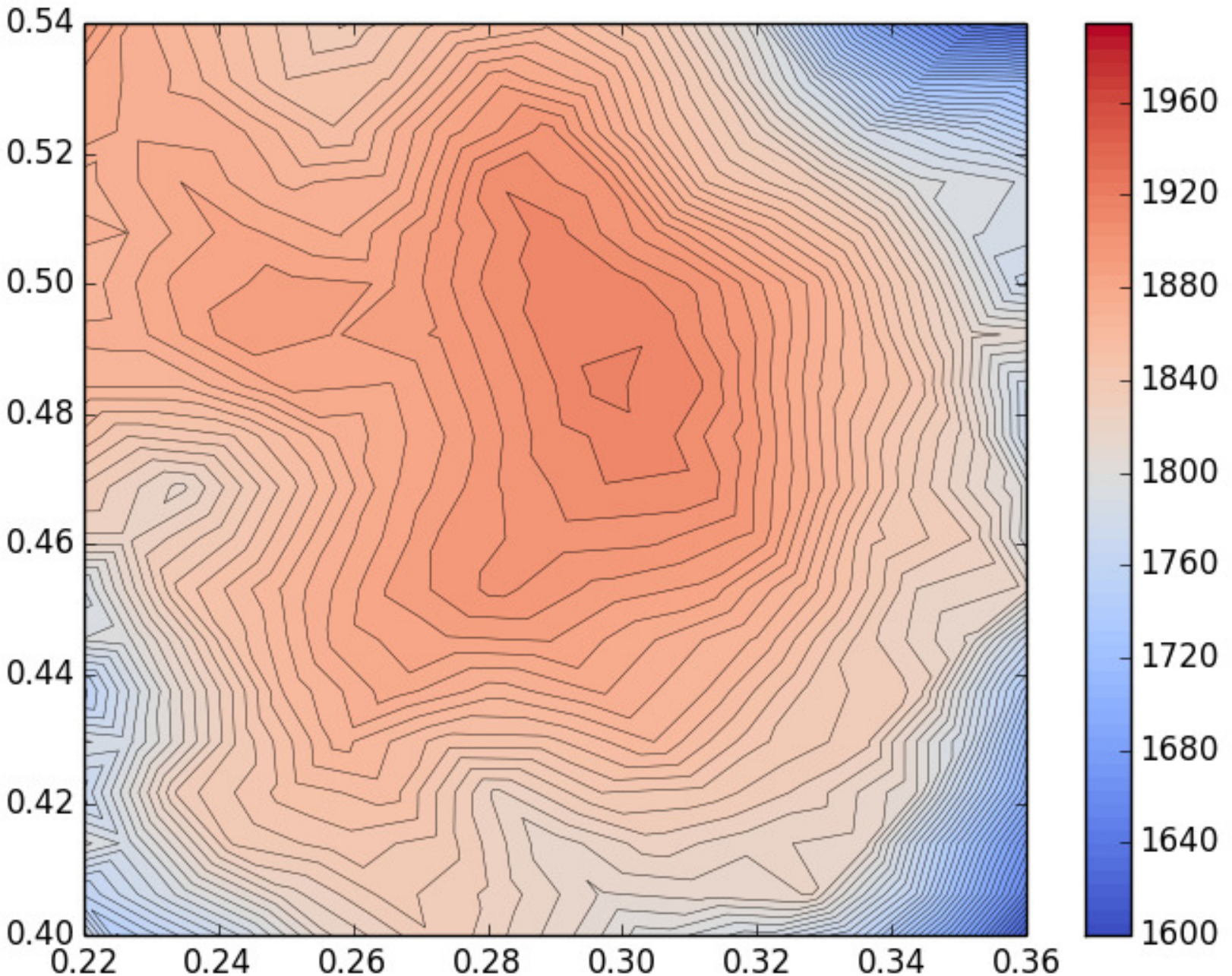}
		\caption{}
		\end{subfigure}
		\hspace{0.5cm}
		\begin{subfigure}[b]{0.45\textwidth}
		\centering
		\includegraphics[width=\textwidth]{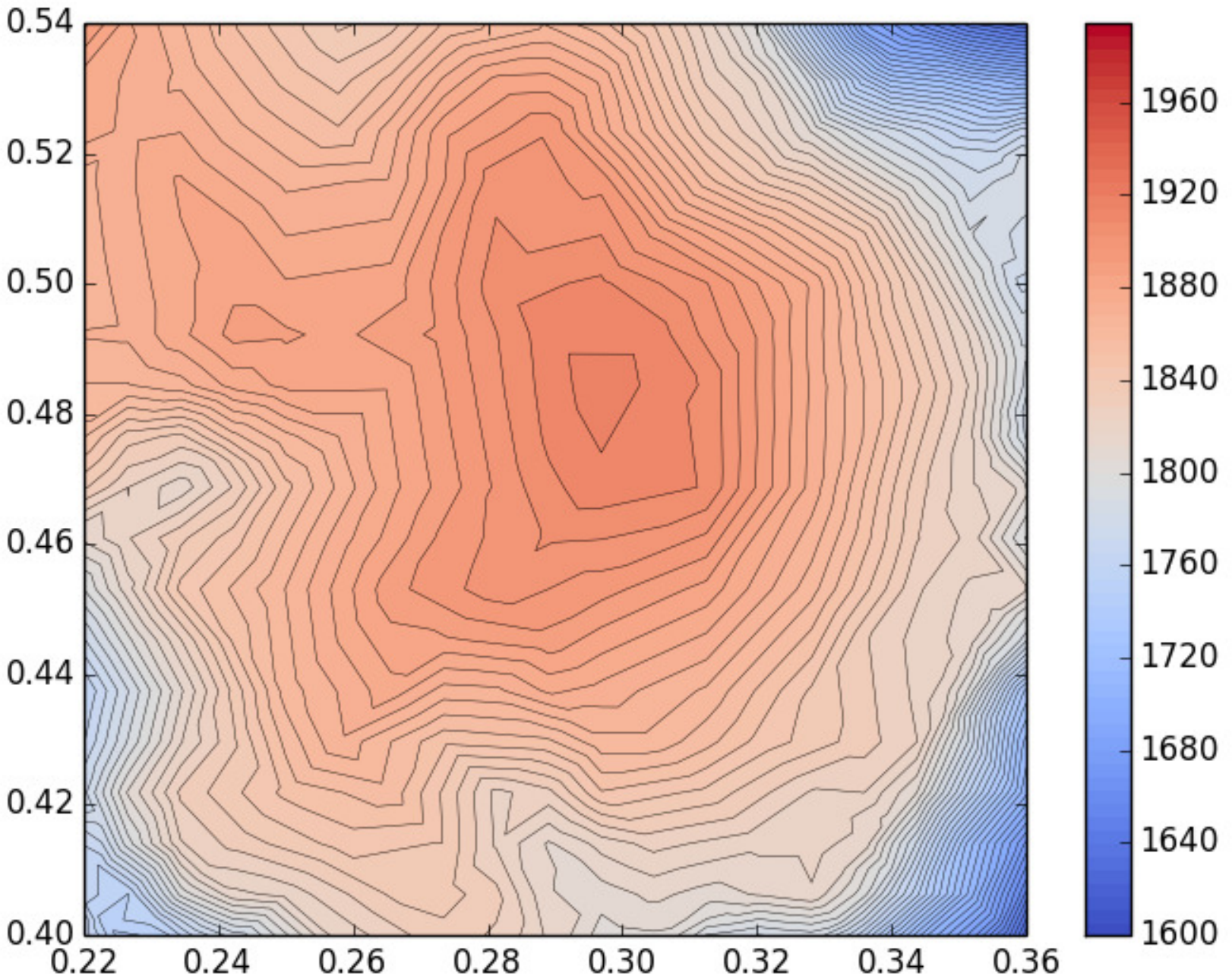}
		\caption{}
		\end{subfigure}
		\caption{Contour plots of $s$ in a small region without data points using (a) uniform refinement; and (b) adaptive refinement.}
		\label{fig:tpsfem_contour_hole}
	\end{figure}

	There is an island on the Crater Lake near $\bm{x}=[0.3,0.5]$ as shown in Figure~\ref{fig:tpsfem_contour}, and this region does not have any data points. Figure~\ref{fig:tpsfem_contour_hole} compares the approximation in this region using uniform refinement (a) and adaptive refinement (b). Note the plots in Figure~\ref{fig:tpsfem_contour_hole} are scaled differently compared to the plots in Figure~\ref{fig:tpsfem_contour} to better show the details.  The plots were obtained using the auxiliary problem error indicator. Although, as discussed in Section~\ref{sec:auxiliary}, this error indicator depends on local data distributions, it is applicable to regions without data points. The resulting smoother shown in Figure~\ref{fig:tpsfem_contour_hole}(a) is similar to that obtained using a uniform grid shown in Figure~\ref{fig:tpsfem_contour_hole}(b).

\begin{figure}
		\centering
		\begin{subfigure}[b]{0.45\textwidth}
		\centering
		\includegraphics[width=\textwidth]{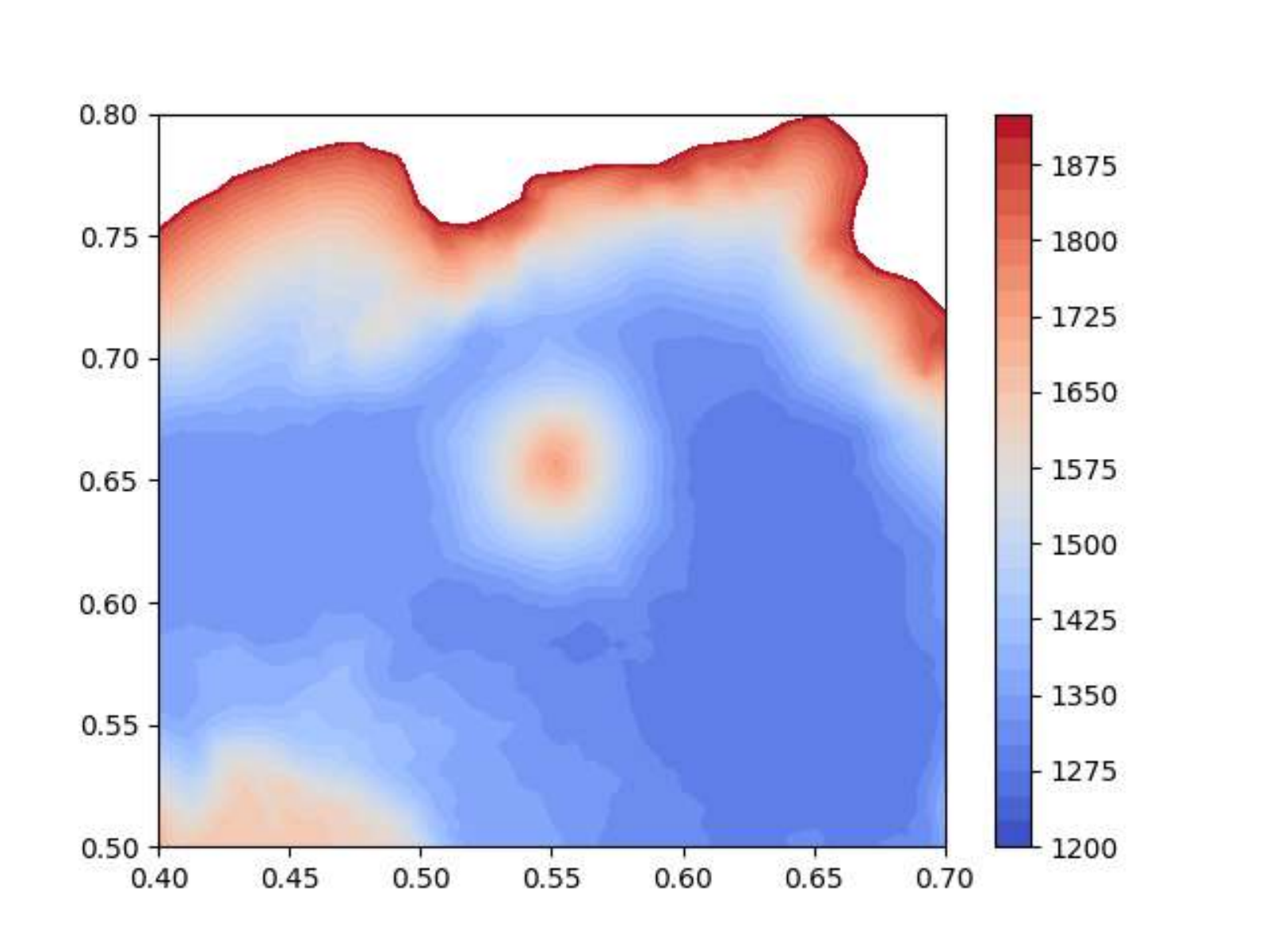}
		\caption{}
		\end{subfigure}
		\hspace{0.5cm}
		\begin{subfigure}[b]{0.45\textwidth}
		\centering
		\includegraphics[width=\textwidth]{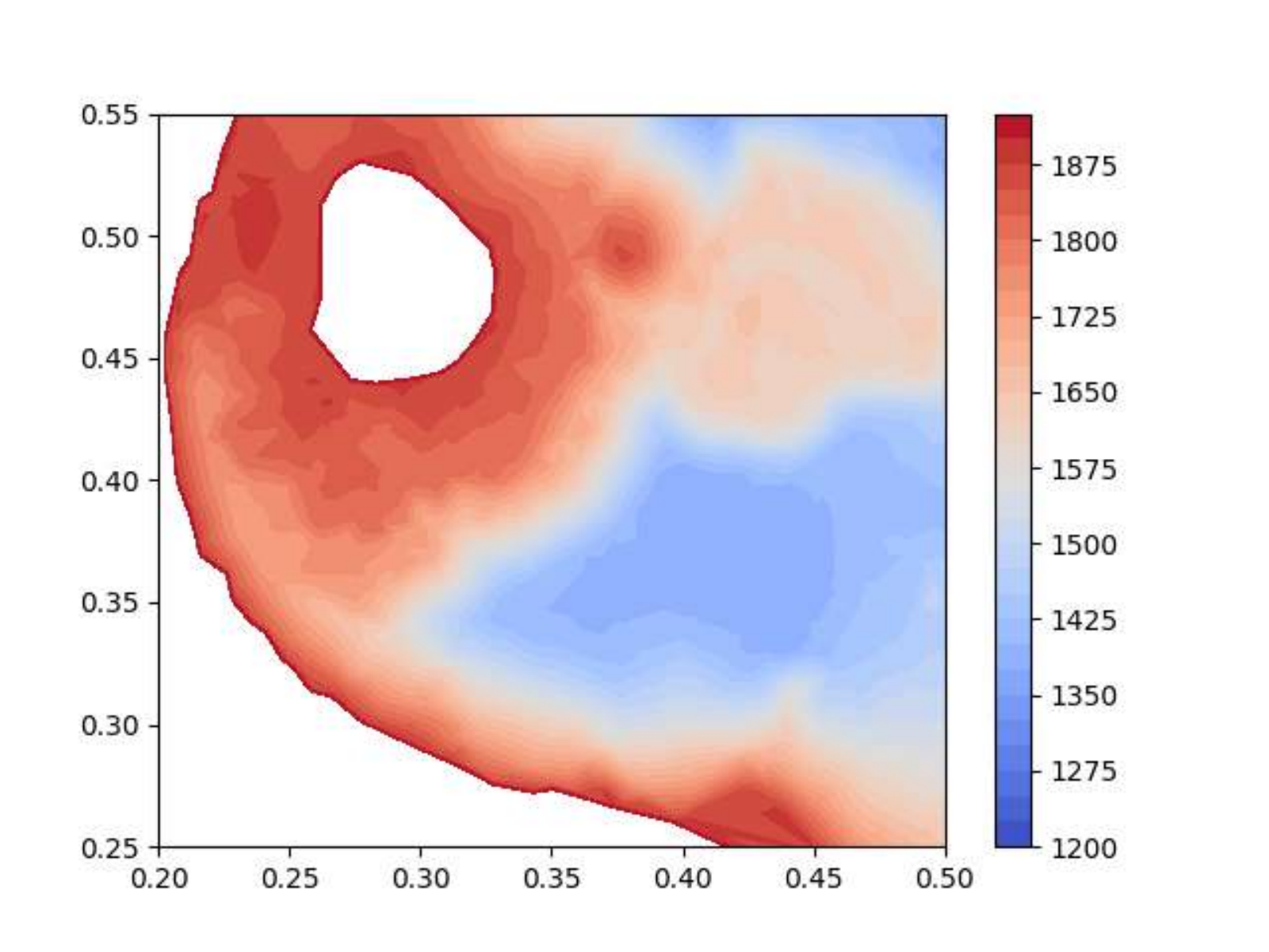}
		\caption{}
		\end{subfigure}
		\caption{Contour maps for Crater Lake data at (a) underwater hill; and (b) shoreline. White spaces represent regions without data points.}
		\label{fig:data_error_hill_shoreline}
	\end{figure}

	The Crater Lake data set consists of oscillatory and smooth regions. The contour maps of the Crater Lake data of two example oscillatory regions on an underwater hill and near shorelines are provided in Figures~\ref{fig:data_error_hill_shoreline}(a) and~\ref{fig:data_error_hill_shoreline}(b), respectively. Areas with white space represent regions without data points. The peak of the underwater hill sits near the point~[0.55,0.65] in Figure~\ref{fig:tpsfem_contour}. Figure~\ref{fig:data_error_hill_shoreline}(b) shows rapid changes in heights near a shoreline.  In contrast, the lower right corner in Figure~\ref{fig:data_error_hill_shoreline}(b) is a smooth region.

\begin{figure}
		\centering
		\begin{subfigure}[b]{0.45\textwidth}
		\centering
		\includegraphics[width=\textwidth]{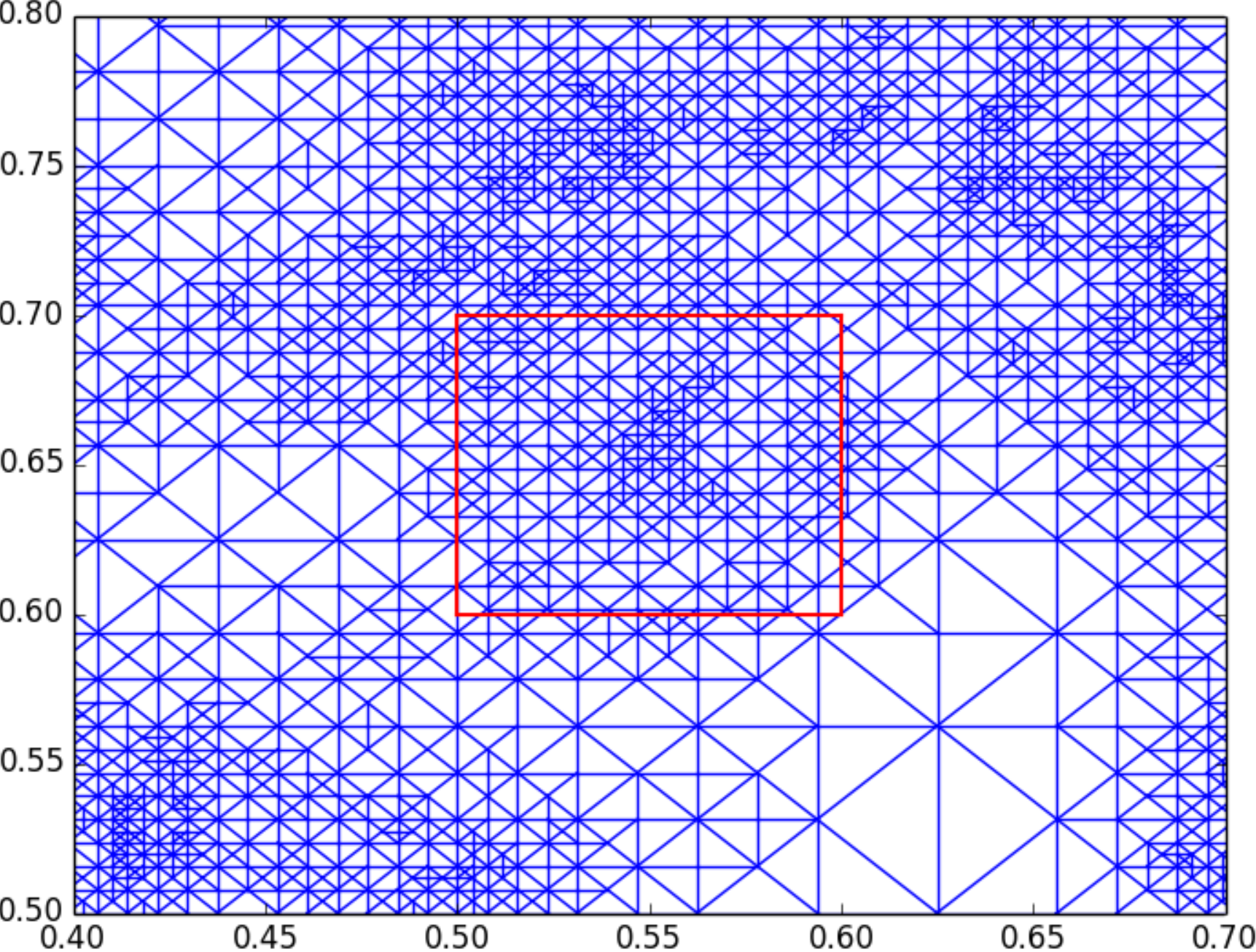}
		\caption{}
		\end{subfigure}
		\hspace{0.5cm}
		\begin{subfigure}[b]{0.45\textwidth}
		\centering
		\includegraphics[width=\textwidth]{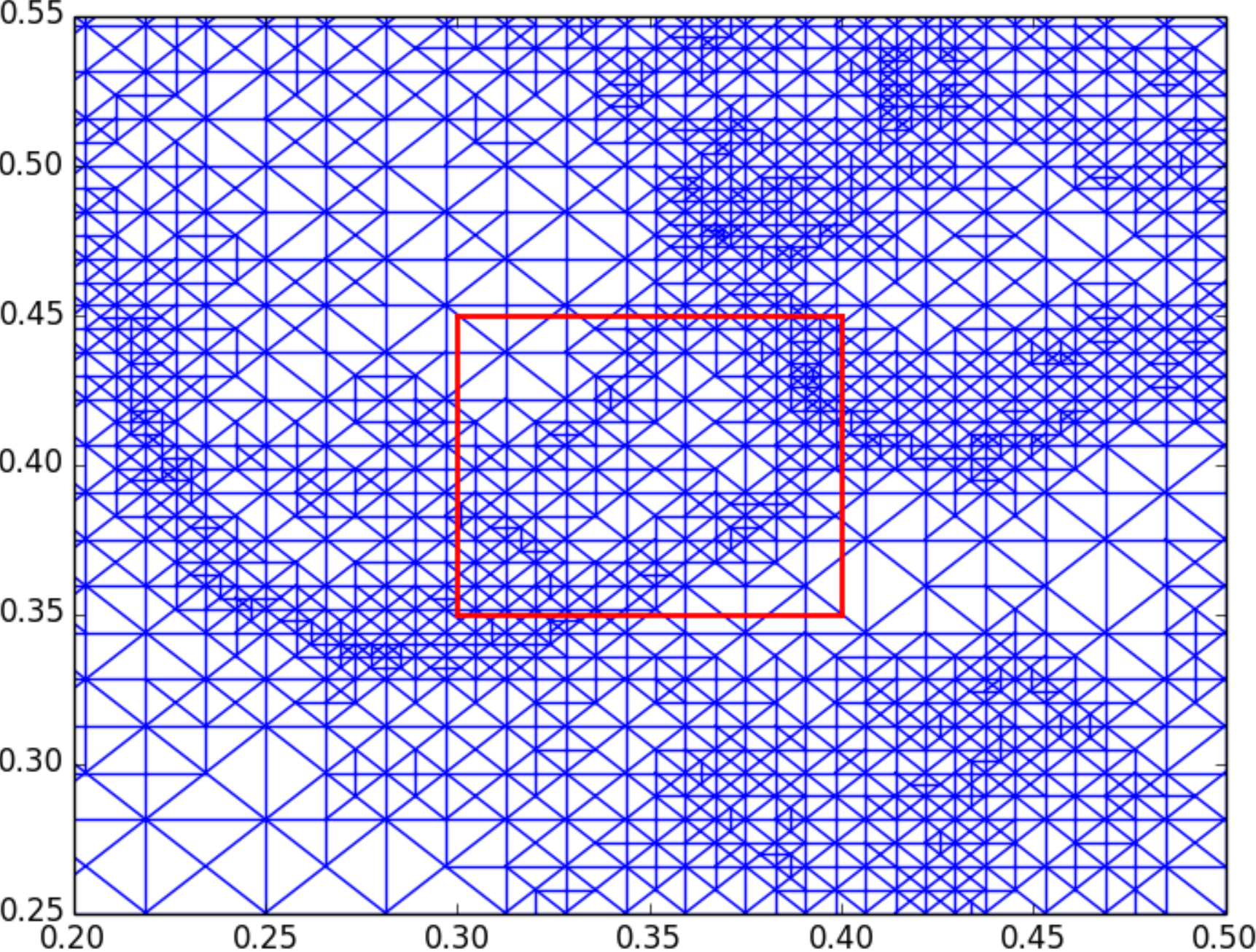}
		\caption{}
		\end{subfigure}
		\caption{Adaptive grids for Crater Lake at (a) underwater hill; and (b) shoreline. The underwater hill and shoreline are highlighted by red squares in subfigures (a) and (b), respectively.}
		\label{fig:grid_hill_shoreline}
	\end{figure}

\begin{figure}
		\centering
		\begin{subfigure}[b]{0.45\textwidth}
		\includegraphics[width=\textwidth]{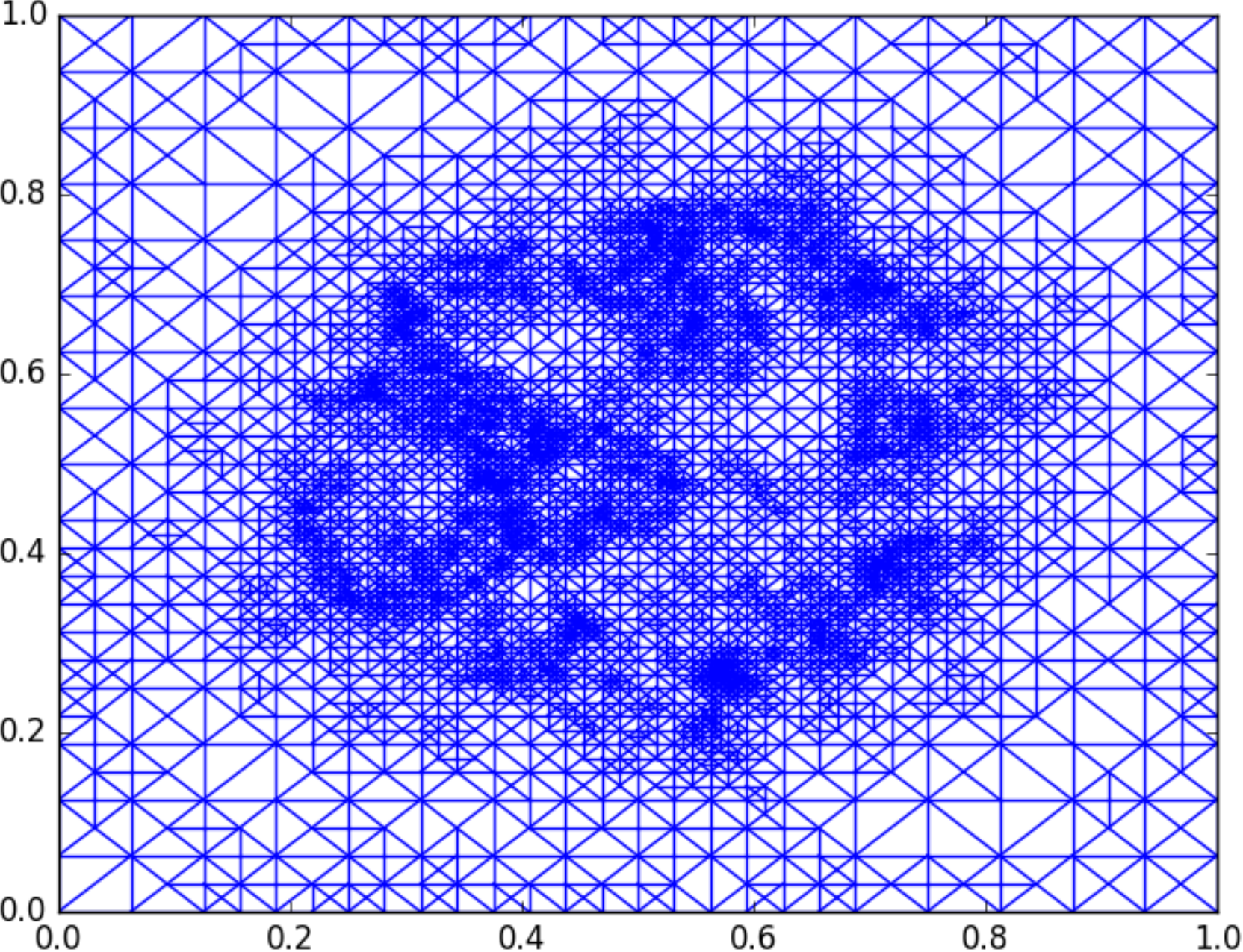}
		\caption{}
		\end{subfigure}
		\hspace{0.5cm}
		\centering
		\begin{subfigure}[b]{0.45\textwidth}
		\includegraphics[width=\textwidth]{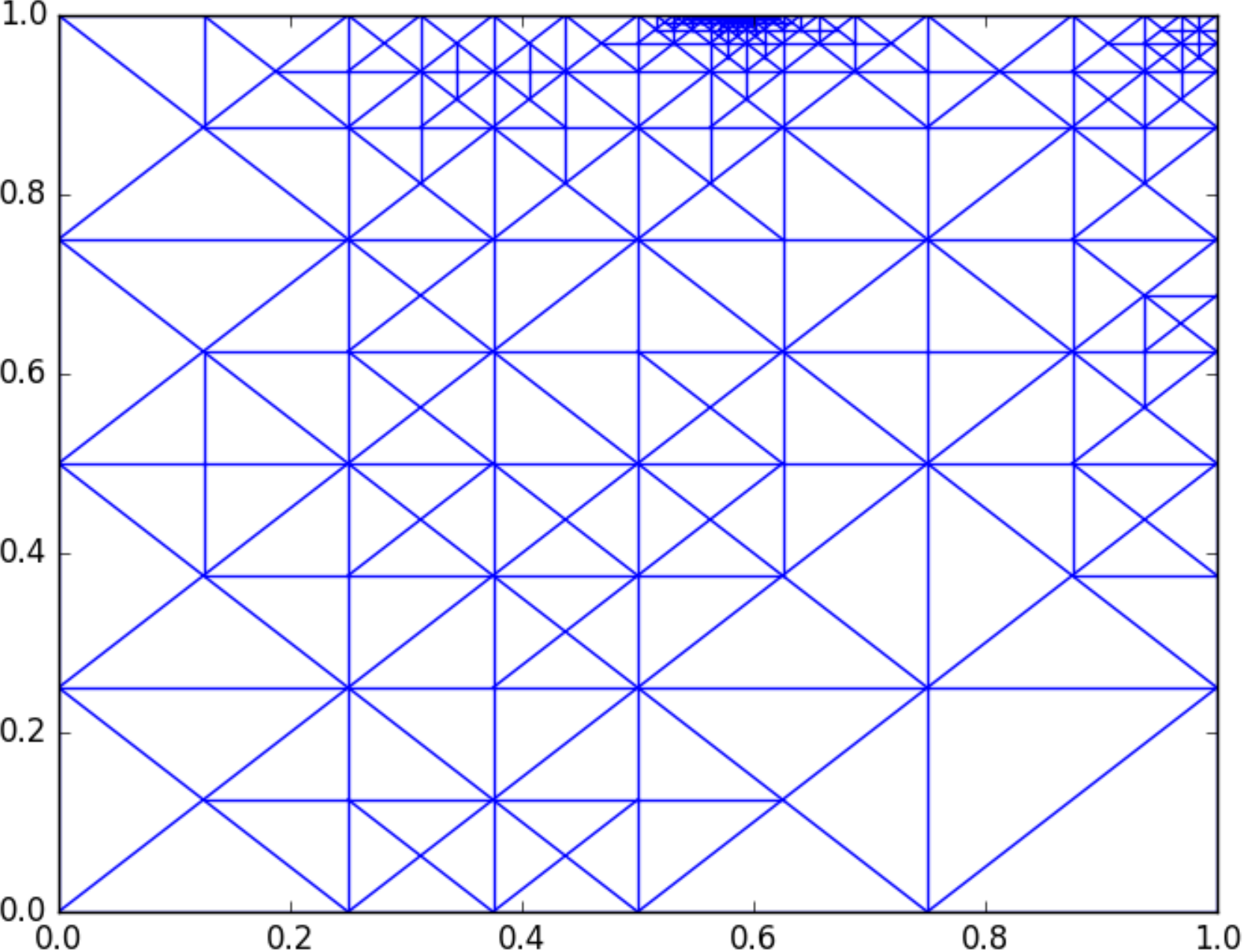}
		\caption{}
		\end{subfigure}
		\caption{Adaptive grid using auxiliary problem error indicator with (a) Dirichlet boundaries; and (b) Neumann boundaries.}
		\label{fig:adaptive_grid}
	\end{figure}

Figure~\ref{fig:grid_hill_shoreline} shows the parts of the adaptive grid corresponding to the data plots in Figure~\ref{fig:data_error_hill_shoreline}. This is the adaptive grid constructed using the norm-based error indicator. Observe the oscillatory regions are refined to improve the accuracy of the FEM approximation. In Figure~\ref{fig:grid_hill_shoreline}(a), finer elements are evident near the peak of the underwater hill that sits near the point~[0.55,0.65]. In contrast, the lower right corner has coarser elements. Similarly, shorelines in Figure~\ref{fig:grid_hill_shoreline}(b) have finer elements compared to the bottom of the lake at the bottom right corner. An example adaptive grid of the whole domain is shown in Figure~\ref{fig:adaptive_grid}(a).
    
\begin{figure}
		\centering
		\begin{subfigure}[b]{0.45\textwidth}
		\centering
		\includegraphics[width=\textwidth]{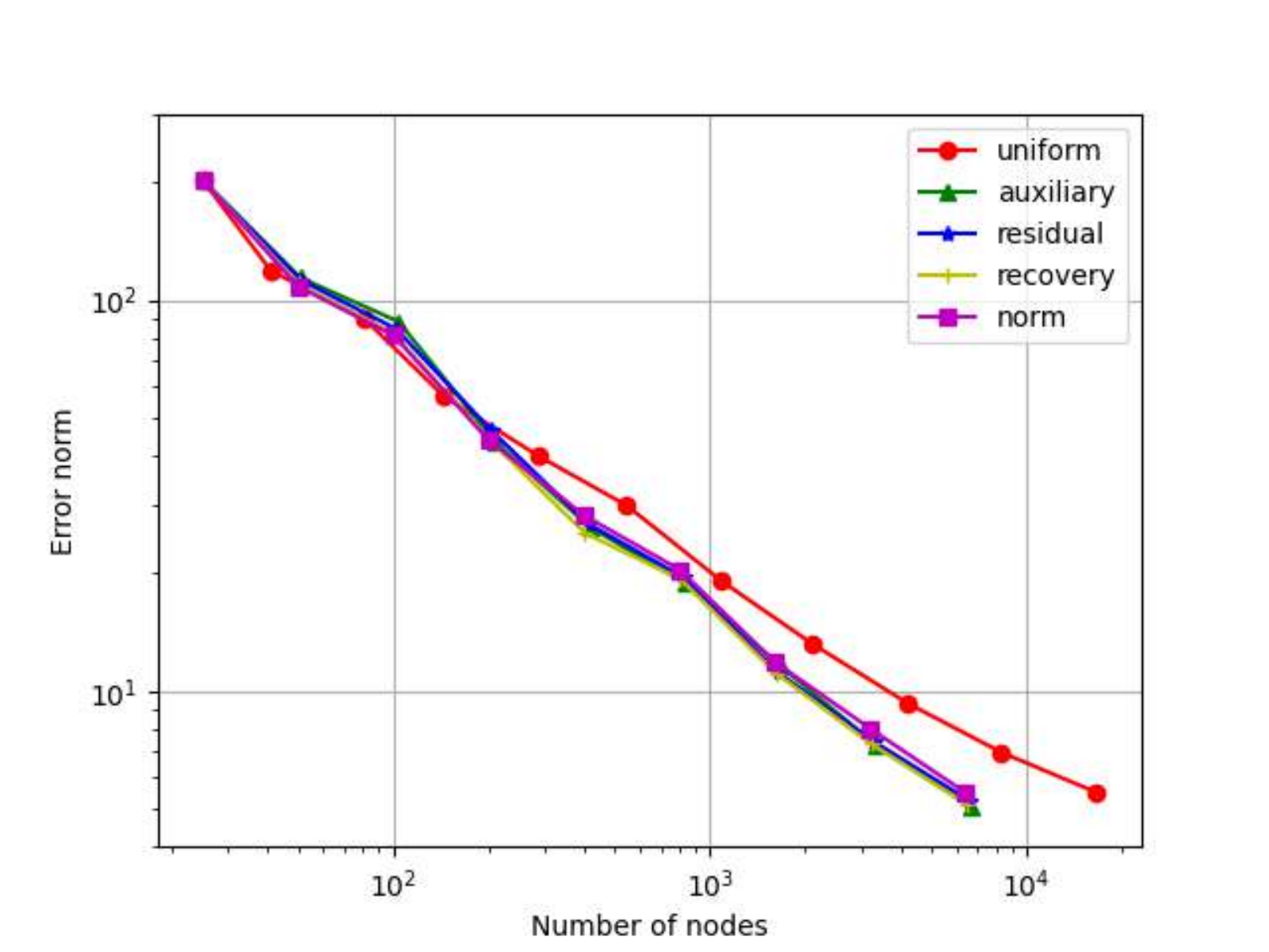}
		\caption{}
		\end{subfigure}
		\hspace{0.5cm}
		\begin{subfigure}[b]{0.45\textwidth}
		\centering
		\includegraphics[width=\textwidth]{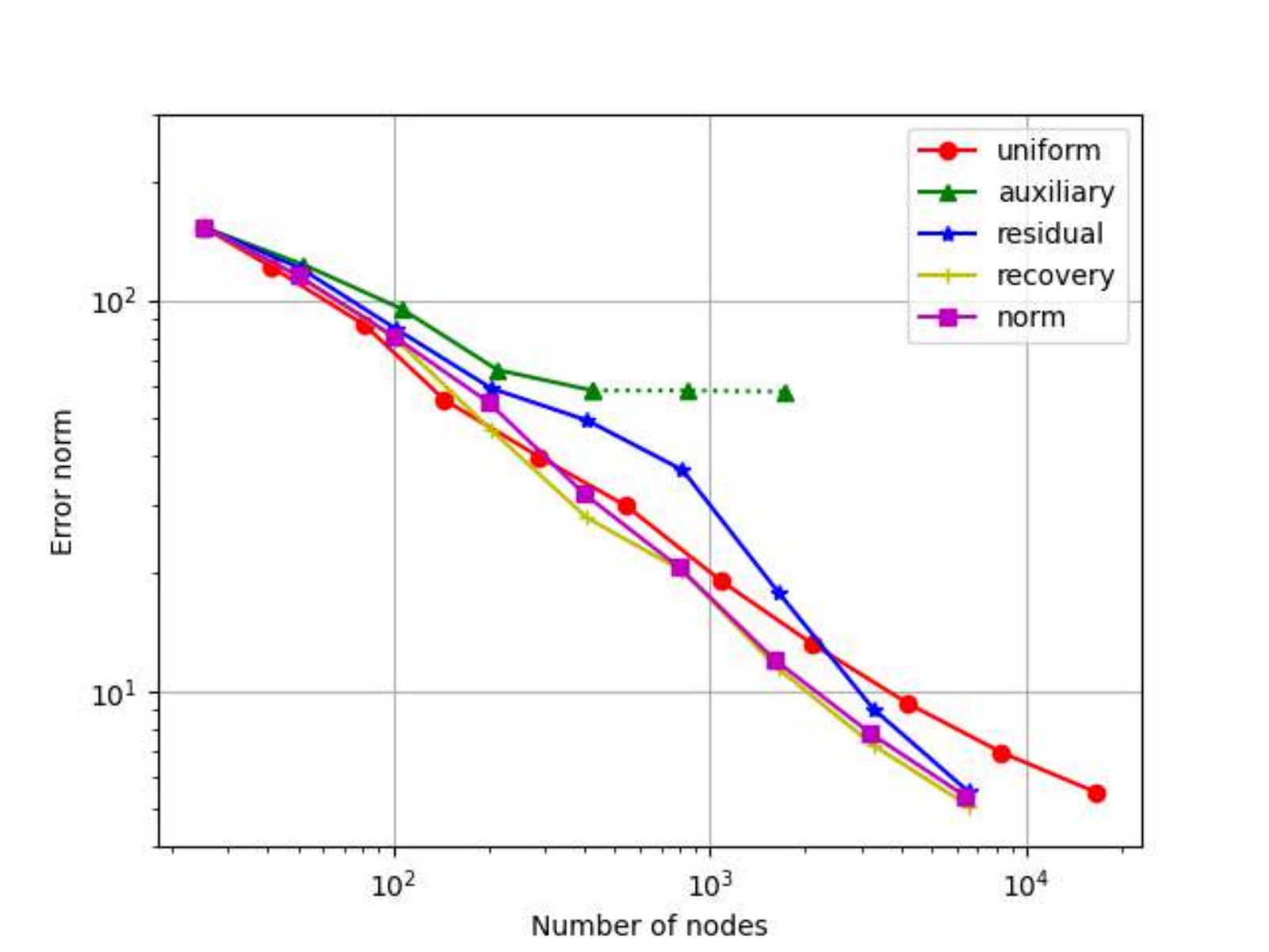}
		\caption{}
		\end{subfigure}
		\caption{RMSE of TPSFEM using uniform and adaptive grids with (a) Dirichlet boundaries; and (b) Neumann boundaries.}
		\label{fig:rmse}
	\end{figure}

	The convergence of the RMSE of the TPSFEM using uniform and adaptive grids for the Crater Lake survey with Dirichlet and Neumann boundaries is shown in Figures~\ref{fig:rmse}(a) and~\ref{fig:rmse}(b), respectively. We applied the same initial grid and refinement setting used in Section~\ref{sec:result_model}. Statistics of refined grids for the Crater Lake data are provided in Table~\ref{tab:metric}. The adaptive grids produced by the four error indicators with Dirichlet boundaries have similar error convergence rates, which are significantly higher than that of the uniform grid as shown in Figure~\ref{fig:rmse}(a). Apart from the auxiliary problem error indicator with Neumann boundaries, the adaptive grids achieve lower RMSE compared to the uniform grids using about $39.21\%$ of the number of nodes, $47.62\%$ of solve time and $85.02\%$ of build time for Dirichlet boundaries; and about $39.04\%$ of the number of nodes, $44.68\%$ of solve time and $80.34\%$ of build time for Neumann boundaries as shown in Table~\ref{tab:metric}. The runtime of the auxiliary problem and residual-based error indicators is markedly higher than that in Table~\ref{tab:metric_model}. In contrast, the runtime of the recovery-based and norm-based error indicators remains relatively stable. The Crater Lake data consists of significantly more data points than the peaks function and the efficiency of the two data-dependent error indicators deteriorates.

	\begin{table}
		\centering
		\caption{Regression metrics of the TPSFEM for Crater Lake data}
		\label{tab:metric}
		\begin{tabular}{lllllll}
		\hline\noalign{\smallskip}
		Boundary & Metric & Uniform & Auxiliary & Residual & Recovery & Norm \\
		\noalign{\smallskip}\hline\noalign{\smallskip}
		Dirichlet & RMSE & 5.50 & 5.06 & 5.25 & 5.14 & 5.50 \\
		 & RMSPE & $2.92\times 10^{-3}$ & $2.72\times 10^{-3}$ & $2.82\times 10^{-3}$ & $2.76\times 10^{-3}$ & $2.95\times 10^{-3}$ \\
		 & MAX & 106.12 & 107.19 & 104.16 & 102.63 & 104.78 \\
		 & \# nodes & 16641 & 6676 & 6563 & 6450 & 6408 \\
          & Solve & 0.42 & 0.22 & 0.19 & 0.19 & 0.20 \\
          & Build & 128.11 & 110.22 & 107.53 & 109.40 & 109.82 \\
          & Indicator &  & 5780.78 & 4061.07 & 145.89 & 95.47 \\
		\noalign{\smallskip}\hline\noalign{\smallskip}
		Neumann & RMSE & 5.50 & 58.39 & 5.53 & 5.09 & 5.38 \\
		 & RMSPE & $2.92\times 10^{-3}$ & $0.03$ & $2.97\times 10^{-3}$ & $2.72\times 10^{-3}$ & $2.89\times 10^{-3}$ \\
		 & MAX & 106.86 & 366.04 & 114.21 & 112.07 & 105.12 \\
		 & \# nodes & 16641 & 1707 & 6544 & 6544 & 6401 \\
		 & Solve & 0.47 & 0.047 & 0.21 & 0.22 & 0.20 \\
          & Build & 137.24 & 105.61 & 107.53 & 113.55 & 114.04 \\
          & Indicator &  & 2894.11 & 4008.11 & 163.14 & 96.20  \\
		\noalign{\smallskip}\hline
		\end{tabular}
	\end{table}

	The recovery-based and norm-based error indicators perform similarly with both the Dirichlet and Neumann boundaries and the former has lower RMSE than the latter. In contrast, the residual-based error indicator did not work as well for Neumann boundaries and the auxiliary problem error indicator did not converge as shown in Figure~\ref{fig:rmse}(b). Thus, the recovery-based error indicator achieves the best performance.

	Example adaptive grids obtained using the auxiliary problem error indicator with Dirichlet and Neumann boundaries are shown in Figure~\ref{fig:adaptive_grid}. With Dirichlet boundary conditions, this error indicator focuses on ridges and shorelines of the lake as illustrated in Figure~\ref{fig:adaptive_grid}(a). Compare Figure~\ref{fig:adaptive_grid}(a) with Figure~\ref{fig:tpsfem_contour}. However, with Neumann boundary conditions all refinement is concentrated on the boundaries as shown in Figure~\ref{fig:adaptive_grid}(b). Given that the refinement concentrated on regions not containing any data points, the RMSE will remain high. Recall that the iterative process of the TPSFEM terminates when the RMSE is not reduced fast enough as shown in the stopping criteria on line 19 of Algorithm~\ref{alg:tpsfem_adaptive}. Consequently, the adaptive refinement process terminated with high errors as shown in Table~\ref{tab:metric}. We have argued that when using Neumann boundary conditions the data needs to be placed well within the domain. The data distribution shown in Figure~\ref{fig:grid_data_lake} demonstrates that the FEM elements near the boundary in the initial coarse grid interact with the data points. The grid in Figure~\ref{fig:grid_data_lake} contains $5 \times 5$ nodes. We reran the above experiment with a different initial coarse grid, one containing $20 \times 20$ nodes, and found that the over-refinement near the boundary disappeared and the resulting grid more closely resembled that in Figure \ref{fig:adaptive_grid}(a).

\subsection{Results: Coastal region}
\label{sec:result_coast}

\begin{figure}
		\centering
		\includegraphics[width=0.45\textwidth]{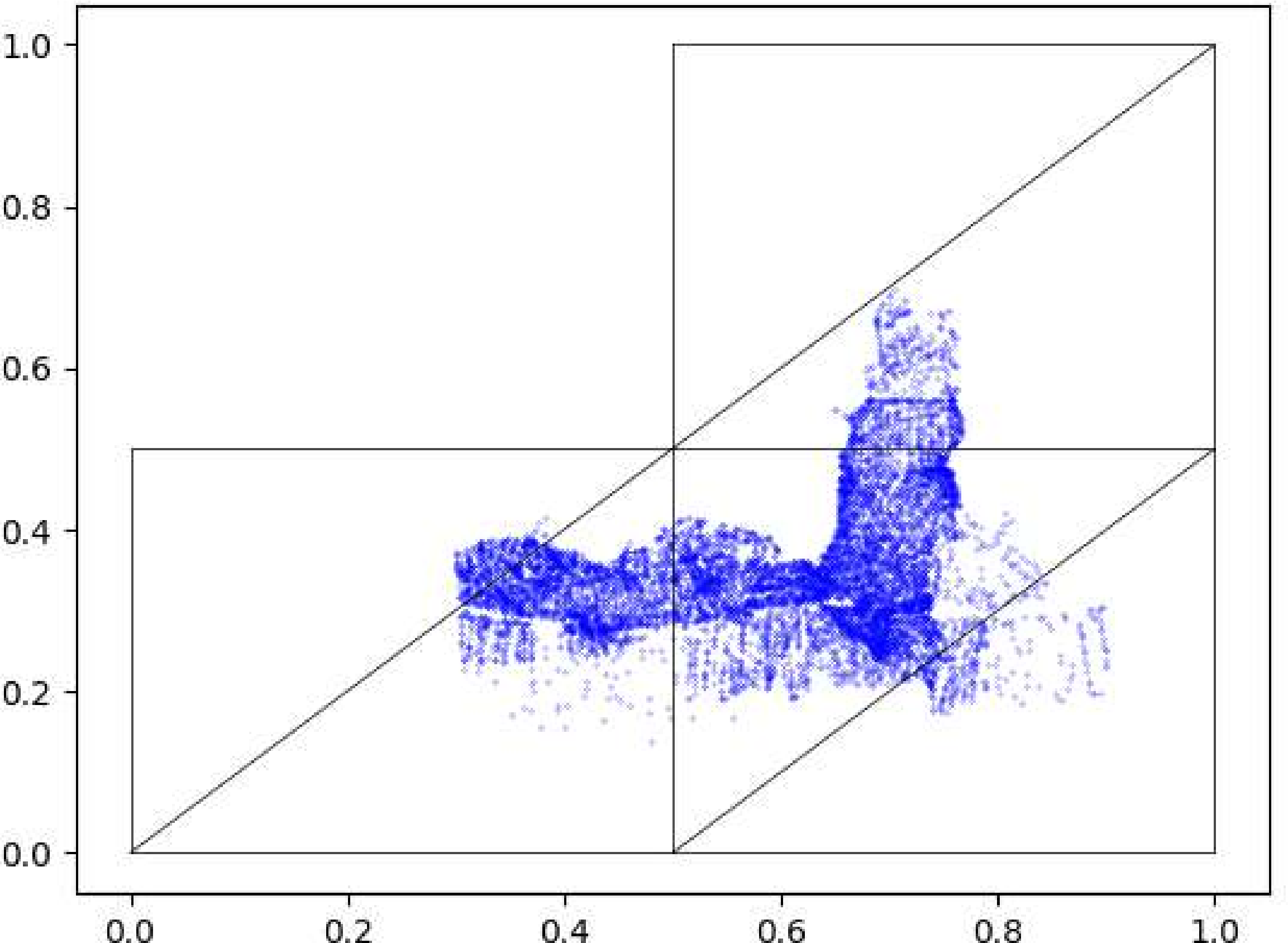}
		\caption{Randomly sampled data points in initial FEM grids for Coastal Region data. Data points are represented as blue dots.}
		\label{fig:grid_data_coast}
	\end{figure}

	With the Coastal Region data, we compared the performance of the TPSFEM in both the L-shaped domain and~$[0,1]^{2}$ square domains. Similar to the Crater Lake data, the Coastal Region data was scaled so the data points fall inside the domain $[0.2, 0.8]^2$. See Figure~\ref{fig:grid_data_coast} for an example coarse grid used in the L-shaped domain. Note that we chose this scale since Coastal Region data has relatively small oscillations near~$\bm{x}=[0.4,0.5]$. As the data values change drastically from the shore to nearby reefs, it is not possible to estimate appropriate values to accompany any Dirichlet boundary conditions. Therefore, we use zero Neumann boundaries.  

	\begin{figure}
		\centering
		\begin{subfigure}[b]{0.32\textwidth}
		\centering
		\includegraphics[width=\textwidth]{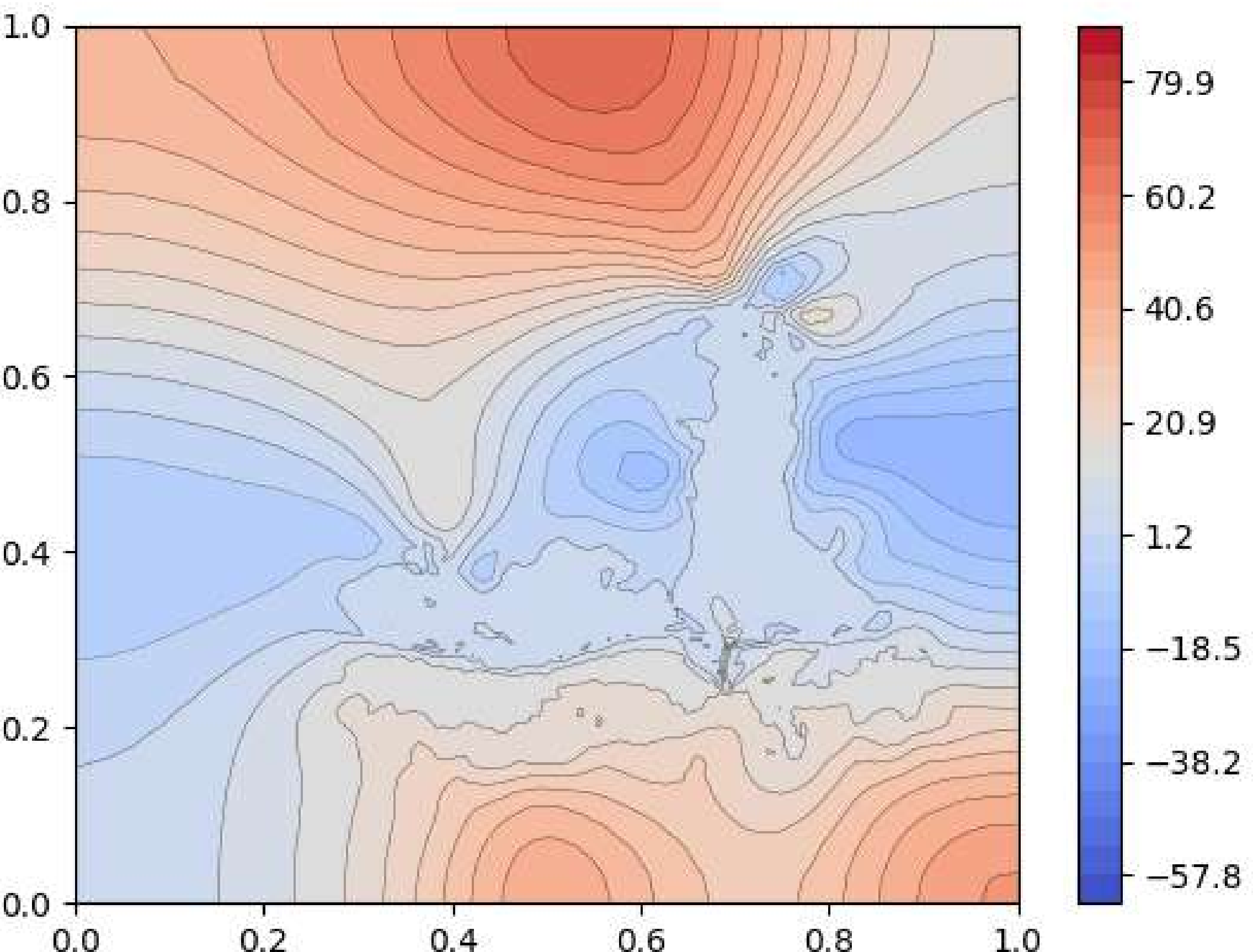}
		\caption{}
		\end{subfigure}
		\begin{subfigure}[b]{0.32\textwidth}
		\centering
		\includegraphics[width=\textwidth]{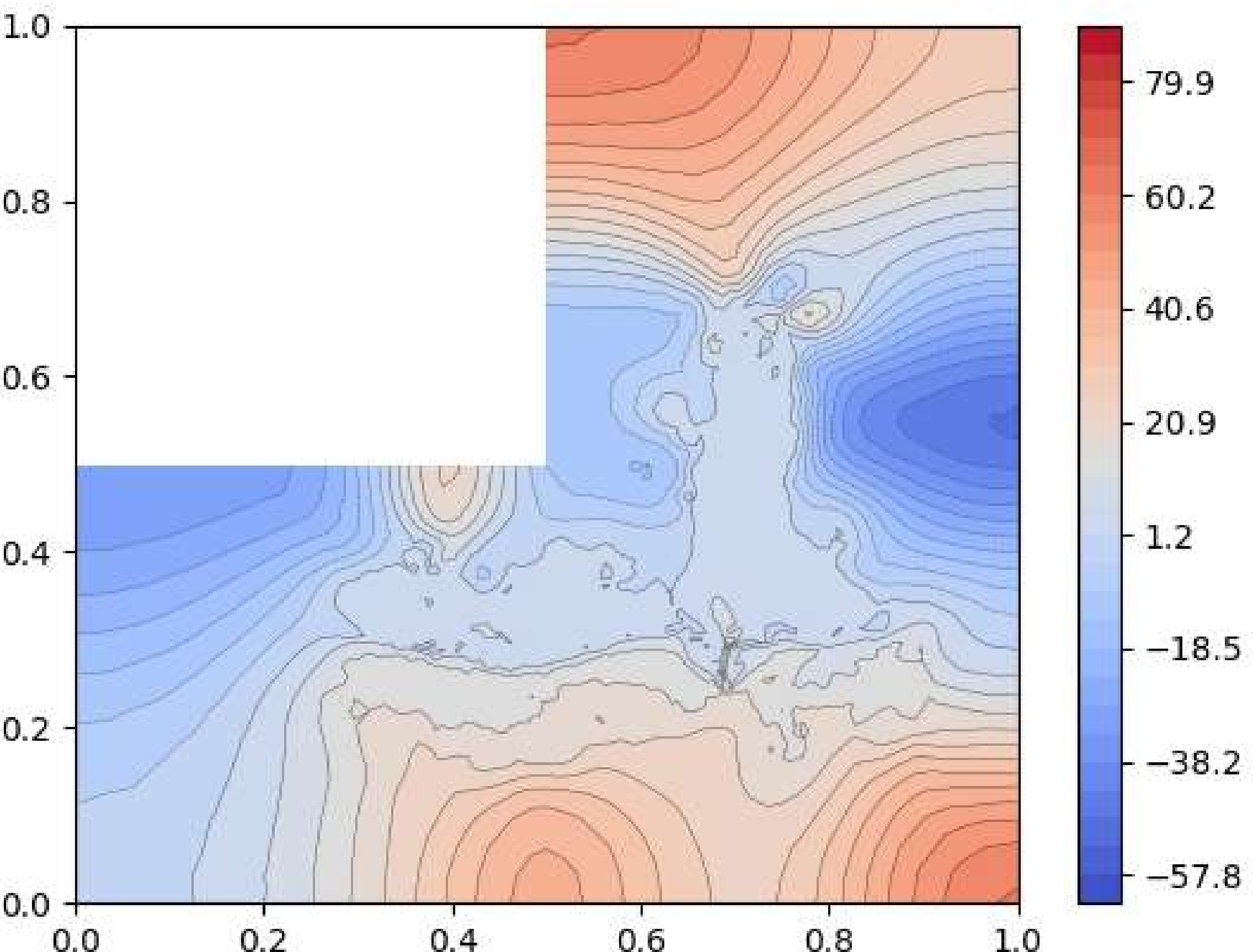}
		\caption{}
		\end{subfigure}
		\begin{subfigure}[b]{0.32\textwidth}
		\centering
		\includegraphics[width=\textwidth]{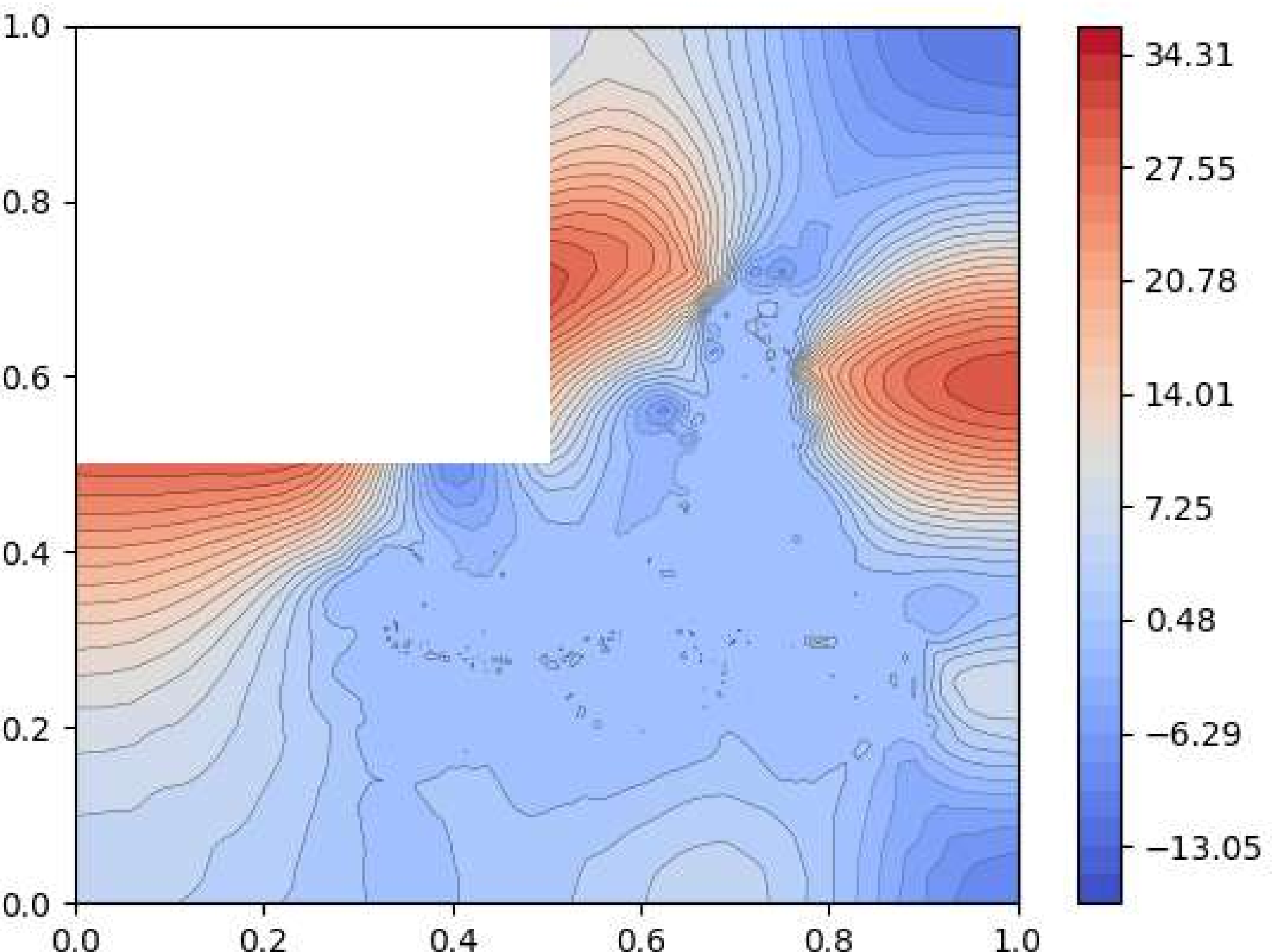}
		\caption{}
		\end{subfigure}
		\caption{Contour plots of TPSFEM using auxiliary problem error indicator in (a) square domain; and (b) L-shaped domain. Contour plot of difference between these two smoothers within the L-shaped domain in (c). Note that Figure (c) is scaled differently to show more details.}
		\label{fig:tpsfem_contour_coast}
	\end{figure}

    Two contour plots of the approximations obtained using square and L-shaped domains for the Coastal Region data are shown in Figures~\ref{fig:tpsfem_contour_coast}(a) and~\ref{fig:tpsfem_contour_coast}(b), respectively. While both smoothers are similar in regions covered by the observed data, the approximation in the L-shaped domain has a steeper descent near the points~$[0.2,0.5]$,~$[0.5,0.7]$ and~$[0.9,0.6]$ as illustrated in Figure~\ref{fig:tpsfem_contour_coast}(c). Note the scaling of Figure~\ref{fig:tpsfem_contour_coast}(c) is different from Figures~\ref{fig:tpsfem_contour_coast}(a) and~\ref{fig:tpsfem_contour_coast}(b) and was chosen to highlight the details.

\begin{figure}
		\centering
		\begin{subfigure}[b]{0.45\textwidth}
		\centering
		\includegraphics[width=\textwidth]{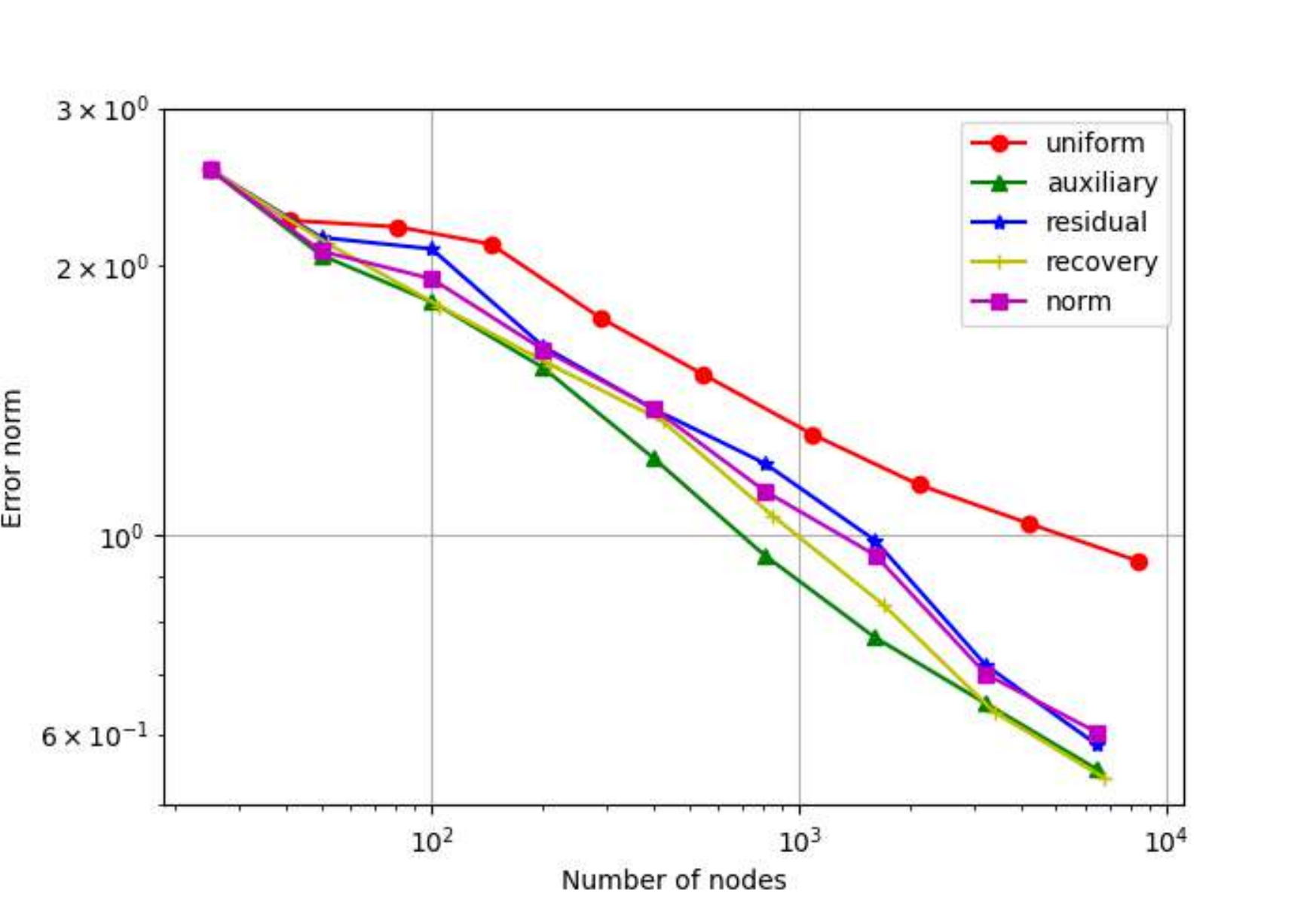}
		\caption{}
		\end{subfigure}
		\hspace{0.5cm}
		\begin{subfigure}[b]{0.45\textwidth}
		\centering
		\includegraphics[width=\textwidth]{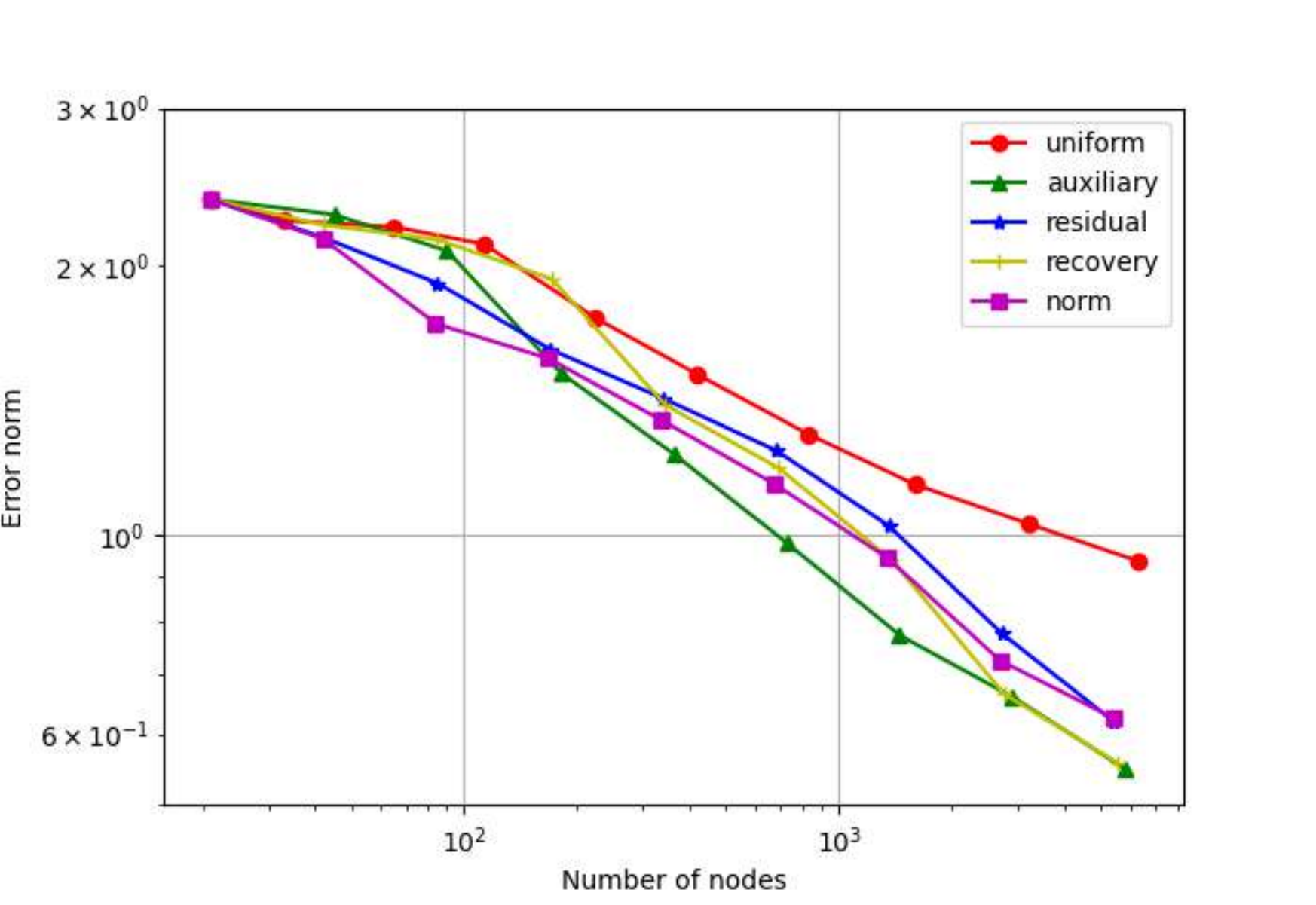}
		\caption{}
		\end{subfigure}
		\caption{RMSE of TPSFEM using uniform and adaptive grids with (a) square domain; and (b) L-shaped domain.}
		\label{fig:rmse_coast}
	\end{figure}

	The convergence of the RMSE of the TPSFEM using uniform and adaptive grids for the Coastal Region data using square and L-shaped domains is shown in Figures~\ref{fig:rmse_coast}(a) and~\ref{fig:rmse_coast}(b), respectively. The initial square and L-shaped grids contain 25 and 21 nodes, respectively. They are both refined using uniform and adaptive refinement for at most 10 and 8 iterations, respectively. Adaptive grids produced by the auxiliary problem and recovery-based error indicators achieve higher efficiency than ones of the residual-based and norm-based error indicators. Statistics of refined grids for the Coastal Region data are provided in Table~\ref{tab:metric_cost}. The adaptive grids have significantly higher error convergence rates compared to the uniform grid in both square and L-shaped domains. They achieve about $60.23\%$ of RMSE compared to the uniform grids using about $39.18\%$ of the number of nodes, $47.78\%$ of solve time and $41.47\%$ of build time for square domains; and about $66.29\%$ of RMSE using about $44.21\%$ of the number of nodes, $58.33\%$ of solve time and $43.60\%$ of build time for L-shaped domains as shown in Table~\ref{tab:metric_cost}. Moreover, the L-shaped domain improves the efficiency of adaptive grids for both uniform and adaptive refinement.

	\begin{table}
		\centering
		\caption{Regression metrics of the TPSFEM for Costal Region data}
		\label{tab:metric_cost}
		\begin{tabular}{lllllll}
		\hline\noalign{\smallskip}
		Domain & Metric & Uniform & Auxiliary & Residual & Recovery & Norm \\
		\noalign{\smallskip}\hline\noalign{\smallskip}
		Square & RMSE & 0.88 & 0.58 & 0.55 & 0.54 & 0.55 \\
		 & RMSPE & 0.030 & 0.019 & 0.020 & 0.018 & 0.020 \\
		 & MAX & 12.51 & 7.07 & 7.12 & 6.17 & 7.07 \\
		 & \# nodes & 16641 & 6416 & 6432 & 6786 & 6448 \\
		 & Solve & 0.45 & 0.21 & 0.21 & 0.22 & 0.22 \\
          & Build & 25.05 & 10.31 & 10.11 & 11.26 & 9.87 \\
          & Indicator &  & 310.403 & 263.77 & 125.10 & 83.63 \\
		\noalign{\smallskip}\hline\noalign{\smallskip}
		L-shaped & RMSE & 0.89 & 0.55 & 0.62 & 0.56 & 0.63 \\
		 & RMSPE & 0.030 & 0.019 & 0.021 & 0.019 & 0.021 \\
		 & MAX & 12.60 & 7.02 & 7.08 & 6.22 & 6.36 \\
		 & \# nodes & 12545 & 5818 & 5440 & 5528 & 5397 \\
          & Solve & 0.33 & 0.20 & 0.19 & 0.19 & 0.19 \\
          & Build & 19.61 & 8.60 & 9.17 & 8.27 & 8.16 \\
          & Indicator &  & 254.40 & 222.95 & 104.27 & 72.30 \\
		\noalign{\smallskip}\hline
		\end{tabular}
	\end{table}

	\begin{figure}
		\centering
		\begin{subfigure}[b]{0.43\textwidth}
		\includegraphics[width=\textwidth]{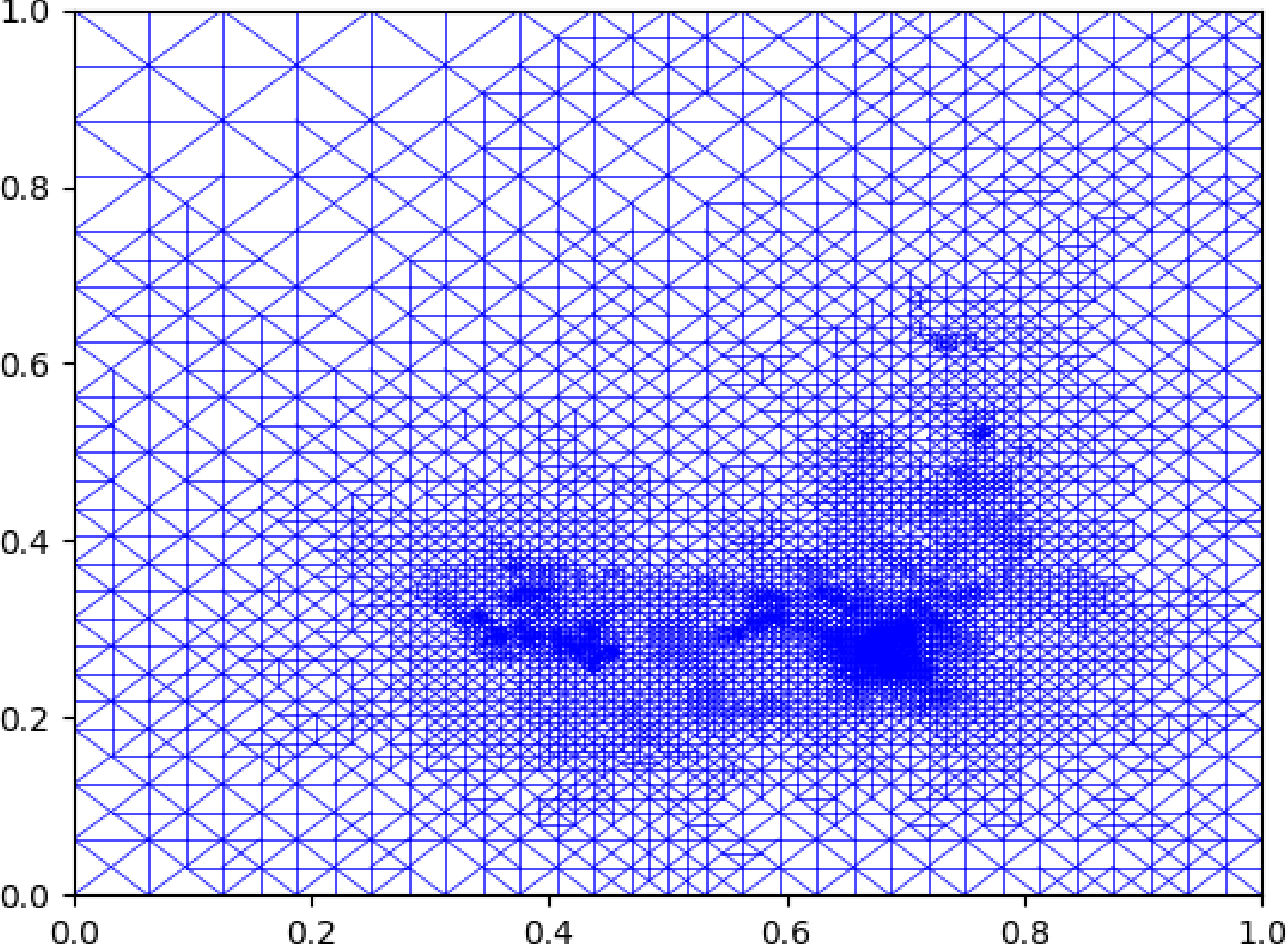}
		\caption{}
		\end{subfigure}
		\hspace{0.5cm}
		\centering
		\begin{subfigure}[b]{0.43\textwidth}
		\includegraphics[width=\textwidth]{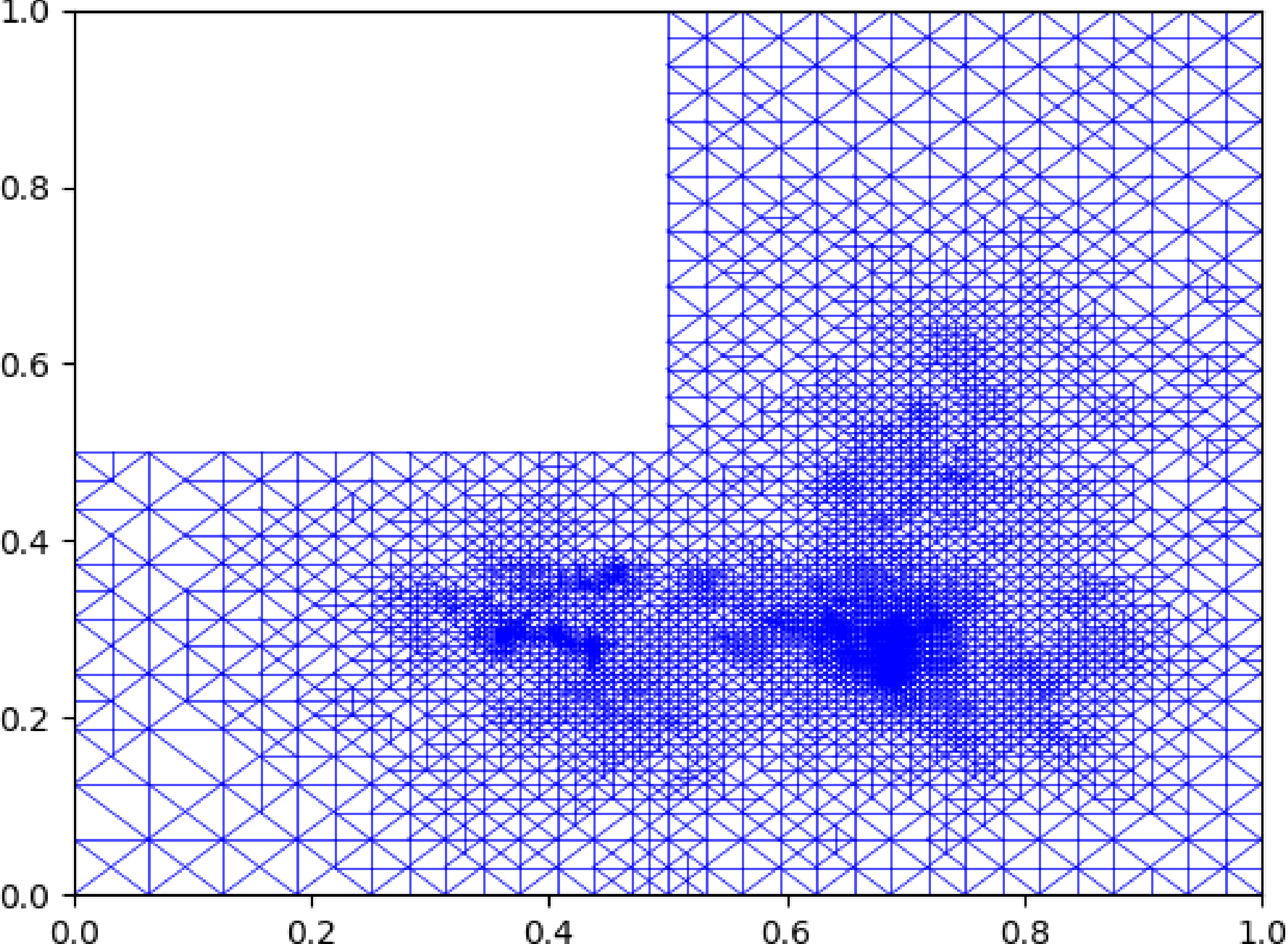}
		\caption{}
		\end{subfigure}
		\caption{Adaptive grid using recovery-based error indicator with (a) square domain; and (b) L-shaped domain.}
		\label{fig:adaptive_grid_coast}
	\end{figure}

	Two example adaptive grids produced using the recovery-based error indicator on the square and L-shaped domains are shown in figures~\ref{fig:adaptive_grid_coast}(a) and~\ref{fig:adaptive_grid_coast}(b), respectively. Refinement is concentrated at the bay area illustrated in Figure~\ref{fig:grid_data_coast}.

\section{Future directions}
\label{sec:future}

The number of nodes in the FEM grid affects the computational costs of many components in Algorithm~\ref{alg:tpsfem_adaptive}, including building and solving the system in Equation~\eqref{eqn:system} and the error indicators. We have provided the runtime of these three components too reflect the impacts of the number of nodes.
The overall runtime for the three data sets using uniform or adaptive refinement is shown in Table~\ref{tab:total}. The TPSFEM for these three cases was built using square domains for consistency. Additionally, the TPSFEM for the peaks function and Crater Lake data used Dirichlet boundaries and the Coastal Region data used Neumann boundaries. Tables \ref{tab:metric_model},
\ref{tab:metric} and \ref{tab:metric_cost} show that adaptive refinement achieved similar RMSE as uniform refinement for the peaks function and Crater Lake data and lower RMSE for the Coastal Region data. The overall runtime includes all the components of Algorithm~\ref{alg:tpsfem_adaptive} and additional components for analysis, including file handling and the GCV, which can be computationally expensive. For instance, the GCV took about 585 seconds for the peaks function using the norm-based error indicator, which occupied about $39.13\%$ of the total runtime.

	\begin{table}
		\centering
		\caption{Total runtime of TPSFEM (measured in seconds)}
		\label{tab:total}
		\begin{tabular}{llllll}
		\hline\noalign{\smallskip}
		Data & Uniform & Auxiliary & Residual & Recovery & Norm   \\
		\noalign{\smallskip}\hline\noalign{\smallskip}
		peaks & 6487 & 2772 & 2429 & 1796 & 1495 \\
		Crater & 9772 & 59111 & 44738 & 4561 & 4248  \\
        Coastal & 7812 & 2850 & 2115 & 1677 & 1572 \\
		\noalign{\smallskip}\hline
		\end{tabular}
	\end{table}

Table~\ref{tab:total} shows that the auxiliary problem and residual-based error indicators are markedly more computationally expensive than the recovery-based and norm-based error indicators. Since the former two directly access data to indicate errors, it takes longer especially when the data size is large. The Crater Lake data consists of 12,936,068 data points, which is much higher than those of the peaks function (62,500 points) and Coastal Region data (48,905 points). Consequently, the runtime of the former two error indicators is high, especially in the first few iterations with coarse grids, where the number of local data points is high. While the runtime of most adaptive refinement cases is lower than their corresponding uniform refinement, the auxiliary problem and residual-based error indicators took more time compared to uniform refinement for the Crater Lake data. Thus, in terms of accuracy and efficiency in the numerical experiments, the recovery-based error indicator is the best out of the four error indicators. 

Many PDE-based error indicators have been developed and achieved superior efficiency compared to uniform refinement. Thus, we adapted some of them and successfully improved the efficiency of the TPSFEM, which has a data-based formulation. While the efficiency of the two data-dependent error indicators deteriorated for large data sets, we have some ideas on how to fix them. For example, the auxiliary problem error indicator can build local approximations with a subset of local data points. This will significantly reduce runtime for large data sets.

We did not optimise our implementations and include more statistics in the article because the runtime is affected by many factors, including the number of nodes, the number of data points and their distribution patterns. It is not the focus of this article and a comprehensive analysis will be too complicated. In addition, we plan to work on the parallel implementation and use of iterative solvers for the TPSFEM. We believe that being able to parallelise the implementation using our local basis functions is another advantage compared to the global basis functions. It will also help to reduce the costs of error indicators, which are highly parallelisable. In future work, we will include a more thorough discussion of the timings of the various components. Nevertheless, the initial results we discuss here clearly demonstrate that adaptive refinement works very well, especially for the PDE-based error indicators.

\section{Conclusion}
\label{sec:conclusion}

	We formulated an adaptive refinement process for a FEM approximation of the thin plate spline. The FEM approximation is called the TPSFEM. We demonstrated this process needs to take into account not just the spacing of the FEM grid, as is standard with adaptive refinement techniques, but also the data distribution as well as the smoothing parameter in the thin plate spline. Given the estimated errors of the FEM approximation do not solely depend on the grid spacing, the usual stopping criteria will not work. Consequently, we provided a detailed description of a new iterative adaptive refinement process. This process adopts stopping criteria that are based on the rate of change of the root-mean-square error (RMSE). Also, a low-cost technique is used to update the smoothing parameter by using the parameter obtained in the previous iteration of the refinement process.

	We studied five error indicators. One error indicator was based on the RMSE. We showed this error indicator should not be used as it is sensitive to noise in the data and leads to over-refinement. We concluded this approach is not worth considering further. The remaining four error indicators are based on examples developed for PDEs. They are the auxiliary problem error indicator, residual-based error indicator, recover-based error indicator and norm-based error indicator. 
	
	We conducted experiments using both Dirichlet and Neumann boundary conditions and the error indicators showed similar performance irrespective of the boundary conditions in most cases. The one exception was the auxiliary problem error indicator which overrefined near one of the boundaries of the FEM grid when using Neumann boundary conditions. Note that this problem did go away when we looked at another set of test problems with a finer initial coarse grid. The other three error indicators worked well with Neumann boundary conditions when the FEM grid is extended well past the position of the data points.

To evaluate the error indicators, we studied one model problem and two bathymetric surveys with varied boundaries and domains. Adaptive refinement outperformed uniform refinement for all three problems. The four PDE-based error indicators gave similar error convergence for the first model problem, which was markedly higher than that obtained when using uniform refinement. For the first bathymetric survey carried out on a crater lake, the four error indicators again gave similar performance and were more efficient than uniform refinement when using Dirichlet boundary conditions. However, the auxiliary problem error indicator did not work well with Neumann boundary conditions while others did. All previous experiments were carried out on a square domain. For the second bathymetric survey of a coastal region, both a square and an L-shaped domain were studied and the L-shaped domain improved the efficiency of the TPSFEM as the FEM grid contained fewer regions without data. The auxiliary problem and recovery-based error indicators performed relatively better than the others in these numerical experiments. However, the auxiliary problem error indicator uses local data and becomes computationally expensive for large data sets. It may also over-refine near oscillatory Neumann boundary surfaces. Overall, the recovery-based error indicator is preferred for its effectiveness and stability.




\section*{Acknowledgments}

\subsection*{Author contributions}

	All authors contributed to the study's conception and design. Material preparation, data collection and analysis were performed by Lishan Fang. The first draft of the manuscript was written by Lishan Fang and all authors commented on previous versions of the manuscript. All authors read and approved the final manuscript.

\subsection*{Financial disclosure}

	There are no financial conflicts of interest to disclose.

\subsection*{Conflict of interest}

	The authors did not receive support from any organization for the submitted work.

\bibliographystyle{elsarticle-num}
\bibliography{elsarticle_tpsfem_fang}

\end{document}